%% file: paper.tex
\documentclass{siamart}
\usepackage{bbm,graphicx,color,euscript,mathtools}
\usepackage{accents,cite}
\usepackage{algorithmic,enumerate}
\reversemarginpar

\Crefname{subsection}{Section}{Subsections}

\input macros

\title{Contour Integral Methods\\ for Nonlinear Eigenvalue Problems: \\ A Systems Theoretic Approach%
\thanks{This work was funded by the U.S. National Science Foundation under grant DMS-1720257.}}
\author{Michael C. Brennan\footnotemark[2] 
        \and Mark Embree\footnotemark[3]
        \and Serkan Gugercin\footnotemark[3]}
\headers{Contour Integral Methods for Nonlinear Eigenvalue Problems}{M. Brennan, M. Embree, and S. Gugercin}

\begin{document}
\maketitle

\renewcommand{\thefootnote}{\fnsymbol{footnote}}
\footnotetext[2]{Center for Computational Engineering, Massachusetts Institute of Technology, Cambridge, MA 02139-4307
(\email{mcbrenn@mit.edu})}
\footnotetext[3]{Department of Mathematics and
Division of Computational Modeling and Data Analytics,
Academy of Integrated Science,
Virginia Tech, Blacksburg, VA 24061
(\email{embree@vt.edu}, \email{gugercin@vt.edu}).}
\renewcommand{\thefootnote}{\arabic{footnote}}

\input abstract
\input intro

\input hankel

\input hankelsys
\input singlept 
\input loewner 
\input rationalforTinv

\input filterfunctions
\input modaltruncation

\input conclude

\section*{Acknowledgements}  
We thank Thanos Antoulas, Jonathan Baker, 
Joe Ball, Alex Grimm, and John Rossi 
for helpful discussions about aspects of this work.

%%%%%%%%%%%%%%%%%%%%%%%%%%%%%%%%%%%%%%%%%%%%%%%%%%%%%%%%%%%%%%%%%%%%%%%%%%%%%%%%

%%%%%%%%%%%%%%%%%%%%%%%%%%%%%%%%%%%%%%%%%%%%%%%%%%%%%%%%%%%%%%%%%%%%%%%%%%%%%%%%

\end{document}

%% file: macros.tex
\newsiamremark{plainremark}{Remark}

\usepackage{mathdots}
\usepackage[bb=boondox, bbscaled=1.06,cal=euler,scr=boondoxo,scrscaled=1.04,frak=euler,frakscaled=1.05]{mathalpha}
\usepackage{sansmath}

\def\cal{\EuScript}

\def\wh{\widehat}
\def\dop{{\rm d}}

\def\iop{{\rm i}}

\def\R{\mathfrak{R}}
\def\C{\mathfrak{C}}

\def\H{\mathbb{H}}
\def\L{\mathbb{L}}

\def\rd{\mathbb{c}} % left-data
\def\ld{\mathbb{b}}  % right-data
\def\rdm{\mathbb{C}} % left-data
\def\ldm{\mathbb{B}}  % right-data

\def\Hsamp{{\bf H}_{\bf LR}}
\def\OO{\mathbcal{O}}
\def\RR{\mathbcal{R}}

\def\bB{\mathbb{B}}
\def\bC{\mathbb{C}}

\def\QL{\mathbb{Q}_{L}}
\def\QR{\mathbb{Q}_{R}}

\usepackage{dsfont}

\def\lp{\theta}
\def\rp{\sigma}

\def\qn{z} % quadrature notes
\def\qw{w} % quadrature weights

\def\pow{k} % the power used in the filter function section 

\def\Ared{\wh{\BA}}
\def\Ered{\wh{\BE}}
\def\Bred{\wh{\BB}}
\def\Cred{\wh{\BC}}

\def\Res{\operatorname{Res}}
\def\rank{{\rm rank}}

\def\Cn{\C^n}

\def\Cmm{\C^{m\times m}}

\def\Cnn{\C^{n\times n}}

\def\Rmm{\R^{m\times m}}

\def\BA{{\bf A}}  
\def\Bb{{\bf b}}\def\BB{{\bf B}}  
\def\Bc{{\bf c}}\def\BC{{\bf C}}  
  
\def\Be{{\bf e}}\def\BE{{\bf E}}  
  
\def\BG{{\bf G}}  
\def\BH{{\bf H}}  
\def\BI{{\bf I}}  
  
\def\BK{{\bf K}}  
\def\BL{{\bf L}}  
\def\BM{{\bf M}}  
\def\BN{{\bf N}}  
  
\def\BP{{\bf P}}  
\def\BQ{{\bf Q}}  
\def\Br{{\boldsymbol{r}}}\def\BR{{\bf R}}  
\def\Bs{{\bf s}}\def\BS{{\bf S}}  
\def\BT{{\bf T}}  
\def\Bu{{\bf u}}  \def\CU{{\cal U}} \def\CBU{{\boldsymbol{\CU}}}
\def\Bv{{\bf v}}\def\BV{{\bf V}}  
\def\Bw{{\bf w}}\def\BW{{\bf W}}  
\def\Bx{{\bf x}}\def\BX{{\bf X}}  
\def\By{{\bf y}}\def\BY{{\bf Y}}  \def\CY{{\cal Y}} \def\CBY{{\boldsymbol{\CY}}}
\def\Bz{{\bf z}}  

\def\Bzero{\boldsymbol{0}}

\def\BLambda{\boldsymbol{\Lambda}}
\def\BLambdas#1{\boldsymbol{\Lambda}^{\kern-1.25pt#1\kern1pt}}   % to better accept a power

\def\BPsi{\boldsymbol{\Psi}}
\def\BSigma{\boldsymbol{\Sigma}}

\def\BAs{\BA{\kern-1.5pt}}
\def\BVs#1{\BV_{\kern-1.25pt #1}}
\def\BWs#1{\BW_{\kern-1.25pt #1}}
\def\BLs#1{\BL^{\kern-1.25pt #1}}

\def\Bell{{\boldsymbol{\ell}}}

\def\Cmod{\BC}
\def\Bmod{\BB}
\def\Amod{\BA}
\def\TF{\BG}
\def\TFr{\BG_{\rm r}}
\def\TFt{\BG_{\rm t}}
\def\nin{n_{\rm i}}
\def\nout{n_{\rm o}}
\def\ns{n}
\def\nred{m}
\def\pole{\lambda}
\def\lrd{\Bc}
\def\rrd{\Bb}
\def\cS{\mathcal{S}}

\mathcode`\@="8000
{\catcode`\@=\active \gdef@{\mkern1mu}}
\def\mydate{\number\day\ {\ifcase\month \or January\or February\or
              March\or April\or May\or June\or July\or August\or
              September\or October\or November\or December\fi}
\number\year}
\DeclareMathSymbol{\widehatsym}{\mathord}{largesymbols}{"65}
\newcommand\lowerwidehatsym{%
  \text{\smash{\raisebox{-1.3ex}{%
    $\widehatsym$}}}}

\def\wt#1{\accentset{\textstyle\lowerwidehatsym}{#1}}

\def\SigmaPts{\mbox{\sansmath $\mathrm{\Sigma}$}}
\def\ThetaPts{\mbox{\sansmath $\mathrm{\Theta}$}}

%% file: abstract.tex
%!TEX root = paper.tex
\begin{abstract}
Contour integral methods for nonlinear eigenvalue problems seek to
compute a subset of the spectrum in a bounded region of the complex plane.
We briefly survey this class of algorithms, establishing a relationship to
system realization techniques in control theory.  This connection motivates
a new general framework for contour integral methods (for linear and nonlinear
eigenvalue problems), building on recent developments in multi-point
rational interpolation of dynamical systems.  
These new techniques, which replace the usual Hankel matrices with 
Loewner matrix pencils,  incorporate general interpolation schemes 
and permit ready recovery of eigenvectors. 
Because the main computations (the solution of linear systems 
associated with contour integration) are identical for these Loewner 
methods and the traditional Hankel approach, a variety of new eigenvalue 
approximations can be explored with modest additional work.
Numerical examples illustrate the potential of this approach.
We also discuss how the concept of filter functions can be employed 
in this new framework, and show how contour methods enable a
data-driven modal truncation method for model reduction.
\end{abstract}

\begin{keywords}
Nonlinear eigenvalue problem, contour integrals, rational interpolation, system realization, Loewner matrices, model reduction, filter function, modal truncation
\end{keywords}

\begin{AMS}
15A18, 65F15, 93B20, 93B30
\end{AMS}

%% file: intro.tex
%!TEX root = paper.tex 
\section{Introduction}
Let $\BT(z):\C\to\C^{n\times n}$ denote an analytic matrix-valued function.
The nonlinear eigenvalue problem (NLEVP) seeks $\lambda\in\C$ and nonzero $\Bv\in\Cn$ such 
that $\BT(\lambda)\Bv=\Bzero$.  
NLEVPs are typically much more 
challenging than standard eigenvalue problems, with infinitely many eigenvalues
possible for problems of finite dimension, $n$; moreover, eigenvectors 
associated with distinct eigenvalues need not be linearly independent.
Many algorithms have been proposed to solve NLEVPs,
ranging from Newton methods that compute one eigenvalue at a time
to linearization algorithms based on local polynomial or rational approximations to $\BT(\lambda)$;
the recent survey of G\"uttel and Tisseur provides a comprehensive overview of 
theory and algorithms~\cite{GT17}, 
complementing earlier surveys by Mehrmann and Voss~\cite{MV04,Vos14}.  
While no one algorithm has yet emerged as a definitive method of choice
for all NLEVPs, a class of algorithms based on contour integration 
of $\BT(z)^{-1}$, initiated by Asakura et al.~\cite{ASTIK09} and Beyn~\cite{Beyn12},
shows much promise.
Inspired by contour methods for 
linear eigenvalue problems~\cite{Pol09,SS03},
these methods provide a basis for new black-box software for NLEVPs~\cite{TP19}.

We seek all the eigenvalues in a bounded (open) domain $\Omega\subset \C$.\ \ Contour 
integral methods build on a fundamental result of Keldysh from the 1950s~\cite{Kel51,Keld71}
that decomposes $\BT(z)^{-1}$ into the sum of a resolvent for a linear operator
corresponding to the eigenvalues of $\BT(z)$ in $\Omega$, 
and a nonlinear remainder with no poles in $\Omega$.
We state a simplified version.

%%%%%%%%%%%%%%%%%%%%%%%%%%%%%%%%%%%%%%%%%%%%%%%%%%%%%%%%%%%%%%%%%%%%%%%%%%%%%%%%
\begin{theorem}[Keldysh]
   \label{thm:Keldysh}
   Suppose $\BT(z):\C\to\C^{n\times n}$ has $m$ eigenvalues
   $\lambda_1, \ldots, \lambda_m$ (counting multiplicity) in the bounded domain $\Omega\subset \C$,
   all semi-simple $($i.e., $\lambda_j$ is a first-order pole of $\BT(z)^{-1}$ for {$j=1,\ldots, m$}$)$.  Then one can write
   \begin{equation} \label{eq:keldysh}
   \BT(z)^{-1} = \BV(z\BI-\BLambda)^{-1}\BW^* + \BN(z),
   \end{equation}
for $\BLambda=\operatorname{diag}(\lambda_1,\ldots,\lambda_m)$, 
     $\BV = [\Bv_1\ \cdots\ \Bv_m]\in\C^{n\times m}$ 
and  $\BW = [\Bw_1\ \cdots\ \Bw_m]\in \C^{n\times m}$.
The matrices $\BV$ and $\BW$ contain the right and left eigenvectors 
$\Bv_j$ and $\Bw_j$ for the eigenvalue $\lambda_j$,
normalized so $\Bw_j^*\BT'(\lambda_j^{})\Bv_j^{} = 1$;
the matrix-valued function $\BN(z)$ is analytic in $\Omega$.
\end{theorem}
%%%%%%%%%%%%%%%%%%%%%%%%%%%%%%%%%%%%%%%%%%%%%%%%%%%%%%%%%%%%%%%%%%%%%%%%%%%%%%%%

For a more general version of Keldysh's theorem that handles 
defective eigenvalues (i.e., Jordan blocks), see~\cite[Thm.~2.8]{GT17}.
Notice that the columns of $\BV$ and $\BW$ need not be linearly independent; indeed, 
while we often envision $m\ll n$, the theorem permits $m>n$, 
in which case the dimensions require 
$\BV$ and $\BW$ to have linearly dependent columns.
\Cref{fig:schematic} provides a schematic illustration of~\cref{thm:Keldysh}.
Note that $\lambda\in\C$ is an eigenvalue if and only if $\BT(\lambda)$ is not invertible, 
i.e., $\lambda$ is a pole of at least one entry of the
matrix-valued function $\BT(z)^{-1}$.
By~\cref{eq:keldysh}, $\BT(\lambda)$ is not invertible at $\lambda\in\Omega$ 
precisely when $\lambda\BI-\BLambda$ is not invertible, 
i.e., $\lambda$ is a conventional eigenvalue of the (diagonal) matrix $\BLambda$.\ \ 
Supposing the number of eigenvalues $m$ inside $\Omega$ is much smaller than the dimension of the problem, $n$, 
the Keldysh decomposition gives a low-dimensional linear component within the larger nonlinear problem.
While \cref{thm:Keldysh} ensures the \emph{existence} of this decomposition, 
it does not explicitly \emph{reveal} the critical linear component 
$\BH(z) := \BV(z\BI-\BLambda)^{-1}\BW^*$ from which we could extract the
eigenvalues in~$\Omega$.
Contour integral techniques give access to this linear part by enabling calculation
of samples $\BH(\sigma)$ for $\sigma\in\C\setminus \overline{\Omega}$, where
$\overline{\Omega}$ denotes the closure of the domain $\Omega$.
All the methods we discuss apply naturally to linear eigenvalue problems too,
in which case $\BN(z)$ is a rational function having poles outside $\overline{\Omega}$.

%%%%%%%%%%%%%%%%%%%%%%%%%%%%%%%%%%%%%%%%%%%%%%%%%%%%%%%%%%%%%%%%%%%%%%%%%%%%%%%%
\begin{figure}[t!]
	\input schematic

\vspace*{-5pt}
	\caption{\label{fig:schematic}
		Schematic illustration of Keldysh's Theorem, showing the $n>m$ case. 
    The number of distinct poles in $\BH(z)$ corresponds to the number of
    distinct eigenvalues in $\Omega$.
	}
\end{figure}
%%%%%%%%%%%%%%%%%%%%%%%%%%%%%%%%%%%%%%%%%%%%%%%%%%%%%%%%%%%%%%%%%%%%%%%%%%%%%%%%

In this paper, we show that contour integral algorithms for eigenvalue problems
are closely related to data-driven system identification techniques from control theory
that use samples of $\BH(\rp)$ at points $\rp\in\C$ 
(or even just \emph{tangential samples} like $\Bell^*\BH(\rp)$,
$\BH(\rp)\Br$, or $\Bell^*\BH(\rp)\Br$ {with $\Bell,\Br \in \C^n$})
to recover the matrices $\BV$, $\BLambda$, and $\BW$ comprising $\BH(z)$.
This connection provides a new perspective on existing contour methods,
and suggests new eigenvalue algorithms based on rational interpolation
and Loewner matrices.  
We begin by identifying $\BH(z) = \BV(z\BI-\BLambda)^{-1}\BW^*$ 
as the \emph{transfer function} for the linear 
$n$-input, $n$-output time-invariant dynamical system
\begin{align}
\Bx'(t) &= \BLambda \Bx(t) + \BW^* \Bu(t) \label{eqn:LTI} \\
\By(t) &= \BV \Bx(t), \nonumber
\end{align}
where $\Bx(t) \in \C^m$, $\Bu(t) \in \C^n$, and $\By(t) \in \C^n $ are the 
\emph{states}, \emph{inputs}, and \emph{outputs} 
of~\cref{eqn:LTI} and $\Bx(0)=\Bzero$; see, e.g., \cite{Ant05b,Kai80}. 
Let ${\CBU}(z)$ and ${\CBY}(z)$ denote the Laplace transforms of $\Bu(t)$ and $\By(t)$. 
The transfer function $\BH(z)$ maps $\CBU(z)$ to $\CBY(z)$, 
i.e., $\CBY(z)= \BH(z)\CBU(z)$.\ \  
To expose this $\BH(z)$ term in the Keldysh decomposition
of $\BT(z)^{-1}$, let $f$ be any function analytic on $\Omega$, 
and compute the contour integral
\begin{equation} \label{eq:int1}
 {1\over 2\pi@\iop} \int_{\partial \Omega} f(z) \BT(z)^{-1}\,\dop z
   = {1\over 2\pi@\iop} \int_{\partial \Omega} f(z) \BH(z)\,\dop z
     + 
     {1\over 2\pi@\iop} \int_{\partial \Omega} f(z)\BN(z)\,\dop z
\end{equation}
about the boundary $\partial\Omega$ of $\Omega$.\ \ 
Since $\BN(z)$ is analytic on $\Omega$, Cauchy's theorem gives
\[ {1\over 2\pi@\iop} \int_{\partial \Omega} f(z)\BN(z) \,\dop z = \Bzero,\]
thus reducing~\cref{eq:int1} to 
\begin{align} \label{eqn:cauchy}
 {1\over 2\pi@\iop} \int_{\partial \Omega} f(z) \BT(z)^{-1}\,\dop z
   &= {1\over 2\pi@\iop} \int_{\partial \Omega} f(z) \BH(z)\,\dop z \\[.5em]
   &= \BV \biggl( {1\over 2\pi@\iop} \int_{\partial \Omega} f(z) (z\BI-\BLambda)^{-1} \,\dop z\biggr)\BW^* 
   = \BV f(\BLambda)\BW^*. \nonumber
\end{align}
(The last step just uses the Cauchy integral formula for functions of matrices~\cite{Hig08}.)
Thus, given only access to $\BT(z)$, one can compute $\BV f(\BLambda)\BW^*$ 
via the contour integral on the left-hand side of~\cref{eqn:cauchy}.
We seek to combine different choices of $f$ to reveal 
the eigenvalues of $\BT(z)$ in~$\Omega$ as efficiently and reliably 
as possible.  

Existing algorithms use $f(z) \equiv 1$, $f(z) = z$, 
or more generally $f(z) = z^k$ for nonnegative integers $k$;
with such choices the integral computes $\BV\BLambdas{k}\BW^*\in\Cnn$.
In systems theory these matrices, called \emph{Markov parameters}, 
play a central role in \emph{realization algorithms} that determine
the transfer function $\BH(z)$ from measurements of $\BV\BLambdas{k}\BW^*$.
We elaborate upon this connection in \Cref{sec:hankel}, 
then apply it later to design a new class 
of eigenvalue algorithms that use rational functions for~$f$.

The integrand in~\cref{eq:int1} involves the inverse $\BT(z)^{-1}$.
For all but the smallest $n$ one prefers to avoid such inversion,
reducing the dimension by applying left and right \emph{probing matrices} 
$\BL\in\C^{n\times \ell}$ and $\BR\in\C^{n\times r}$:
\begin{equation} \label{eq:fprobe}
 \BL^*\BV f(\BLambda)\BW^*\BR = {1\over 2\pi@\iop} \int_{\partial \Omega} f(z)\, \BL^*\BT(z)^{-1}\BR\,\dop z \in \C^{\ell\times r}.
\end{equation}
For example, 
$\BT^{-1}(z)\BR \in \C^{n\times r}$ in the integrand 
can be computed (for a fixed $z$) by solving $\BT(z)\BPsi(z) = \BR$ for $\BPsi(z)$.
(The use of probing matrices resembles \emph{sketching} techniques
in randomized numerical linear algebra~\cite{TYUC17,Woo14}.)
We shall show how the $\BL$ and $\BR$ matrices relate to
\emph{tangential interpolation directions} in the context of system 
realization~\cite{ABG10,MA07}.
In contrast to~\cref{eqn:LTI}, 
we identify 
\[ \Hsamp(z) := \BL^*\BV(z\BI-\BLambda)^{-1}\BW^*\BR\]
as the transfer function of the $r$-input, $\ell$-output 
\emph{sampled} dynamical system
\begin{align} 
\Bx'(t) &= \BLambda \Bx(t) + \BW^*\BR \Bu(t) \label{eq:LRsys} \\
\By(t) &= \BL^*\BV \Bx(t),  \nonumber
\end{align}
which, subject to mild conditions on $\BL$ and $\BR$,
has the same poles as $\BV(z\BI-\BLambda)^{-1}\BW^*$:
the eigenvalues of $\BT(z)$ in $\Omega$.

%%%%%%%%%%%%%%%%%%%%%%%%%%%%%%%%%%%%%%%%%%%%%%%%%%%%%%%%%%%%%%%%%%%%%%%%%%%%%%%%
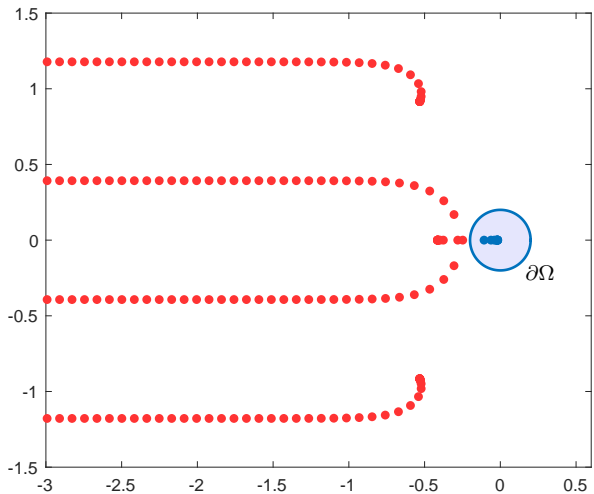
\begin{figure}[b!]
\begin{center}
\includegraphics[scale=0.5]{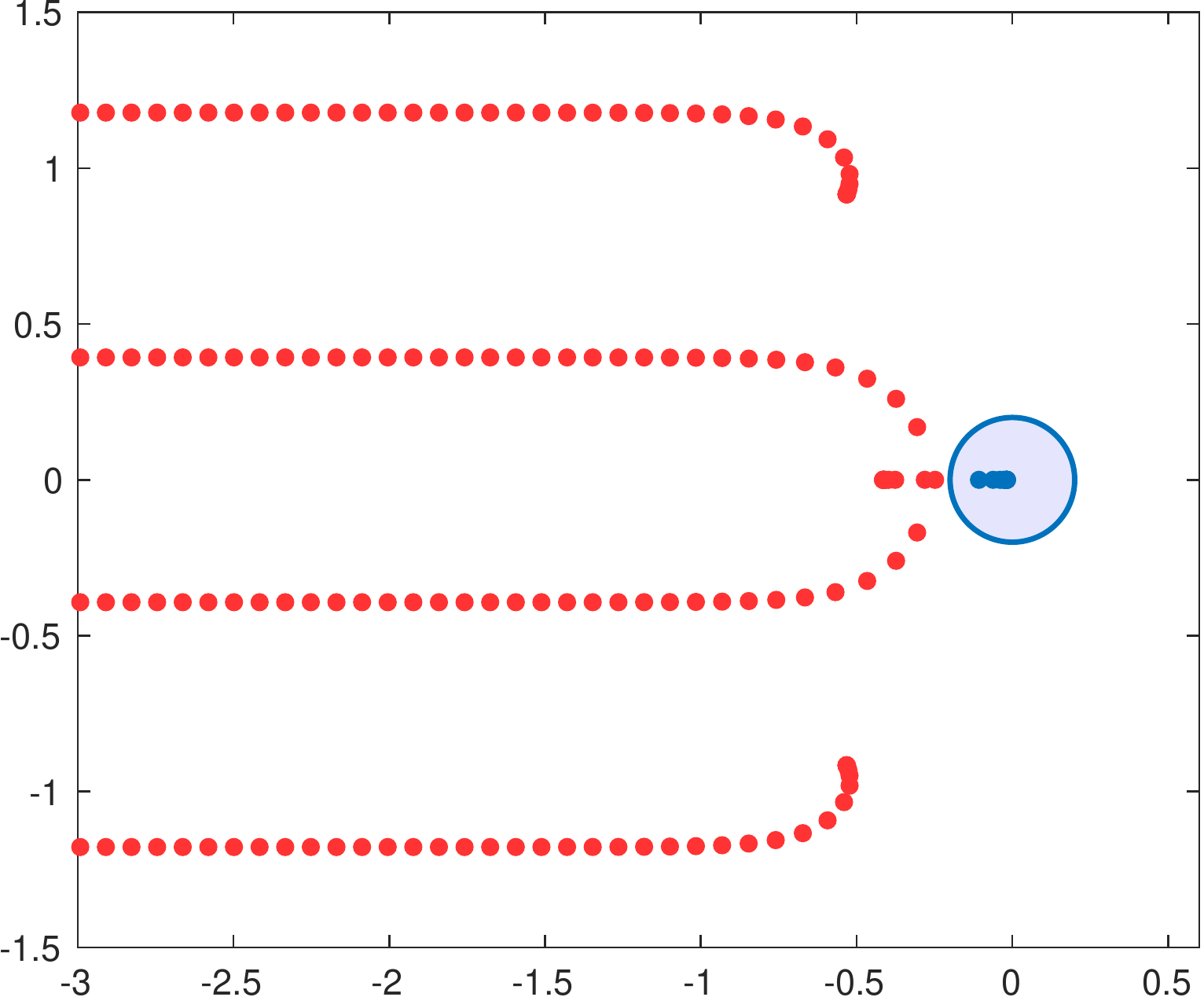}

\begin{picture}(0,0)
\put(85,92){\footnotesize $\partial\Omega$}
\end{picture}
\end{center}

\vspace*{-1.5em}
\caption{\label{fig:contour_sketch}
Sketch of the set-up in the complex plane: computing 11~eigenvalues 
(cluster of blue dots)
of an NLEVP for stability analysis of a delay differential equation.  
Integrate around $\partial\Omega$, the boundary of the blue region,
to compute the eigenvalues inside.  The red dots show some of the
infinitely many other eigenvalues of $\BT(z)$.  
While $\Omega$ is a circle, other shapes are easy to implement.
}
\end{figure}
%%%%%%%%%%%%%%%%%%%%%%%%%%%%%%%%%%%%%%%%%%%%%%%%%%%%%%%%%%%%%%%%%%%%%%%%%%%%%%%%

In practice, the integral~\cref{eqn:cauchy} is approximated 
via a quadrature rule, giving
\begin{equation} \label{eq:trap}
 \BL^*\BV f(\BLambda)\BW^*\BR 
   \ \approx\ 
   \sum_{k=1}^N \qw_k@  f(\qn_k)\, \BL^*\BT(\qn_k)^{-1}\BR
\end{equation}
for quadrature weights $\{\qw_k\}_{k=1}^N$ and nodes {$\{\qn_k\}_{k=1}^N$}. 
Since the individual $\BL^*\BT(\qn_k)^{-1}\BR$ values in \cref{eq:trap} are independent of
one another, they can be computed in parallel, a major
appeal of these contour integration algorithms.
Moreover, we emphasize that the computation of
$\BL^*\BT(\qn_k)^{-1}\BR$ for each quadrature node $\qn_k$ dominates the computational complexity;
with those quantities in hand, it is easy to evaluate~(\ref{eq:trap}) 
with various different $f$.

The trapezoidal rule~\cite{TW14} is often the method of choice in \cref{eq:trap},
 although one can tailor the quadrature rule to the problem through the use of 
\emph{rational filter functions}~\cite{BK16}. 
The term \emph{rational} in that setting is quite different from the concept
of (multi-point) rational interpolation that is the main theme of this paper. 
This distinction will be clarified in \Cref{section:filterfunctions},
after we establish our proposed framework.

%%%%%%%%%%%%%%%%%%%%%%%%%%%%%%%%%%%%%%%%%%%%%%%%%%%%%%%%%%%%%%%%%%%%%%%%%%%%%%%%

%%%%%%%%%%%%%%%%%%%%%%%%%%%%%%%%%%%%%%%%%%%%%%%%%%%%%%%%%%%%%%%%%%%%%%%%%%%%%%%%
\begin{figure}[t!]
\hfill 
\includegraphics[scale=0.36]{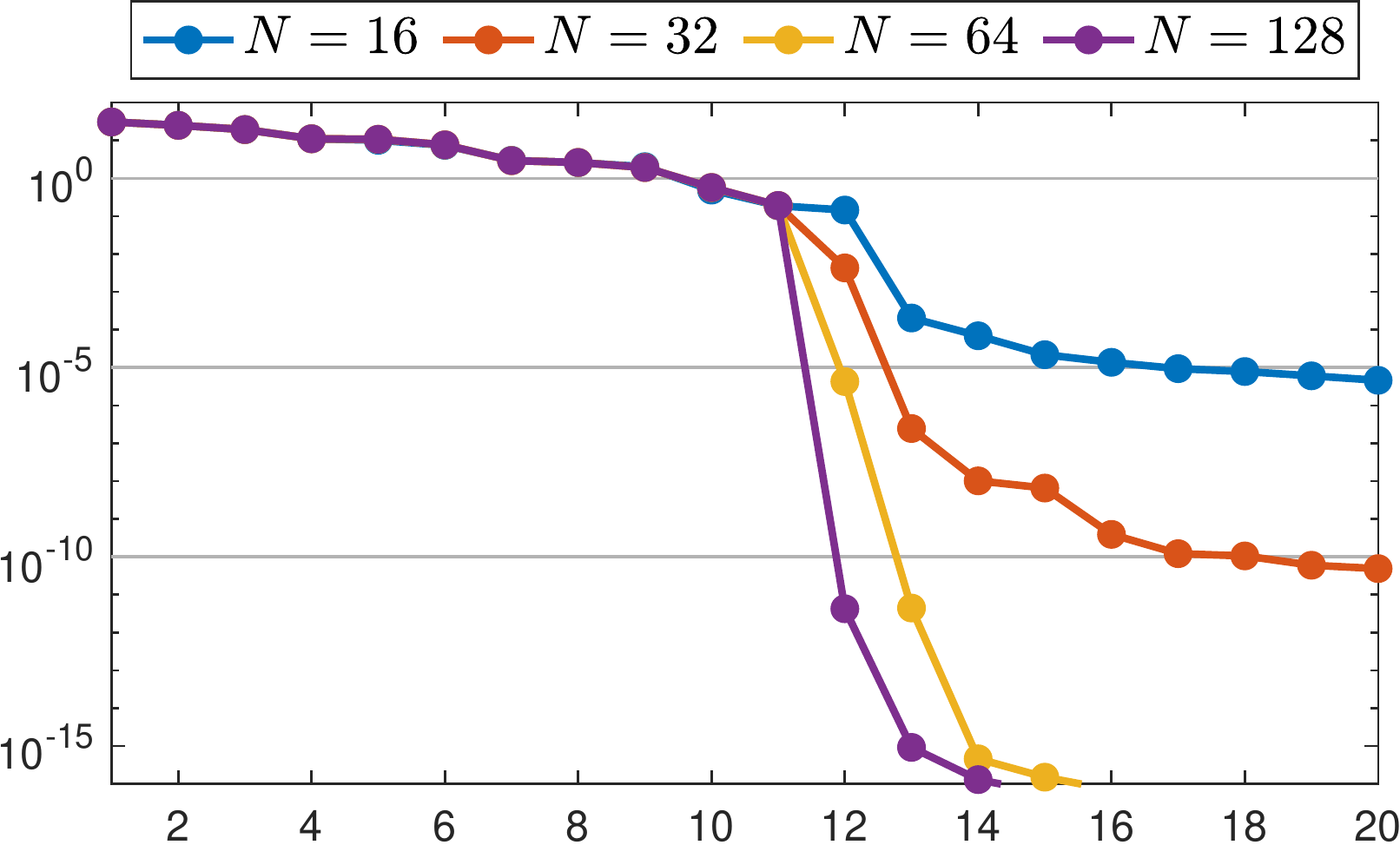}\hspace*{1.75em}
\includegraphics[scale=0.36]{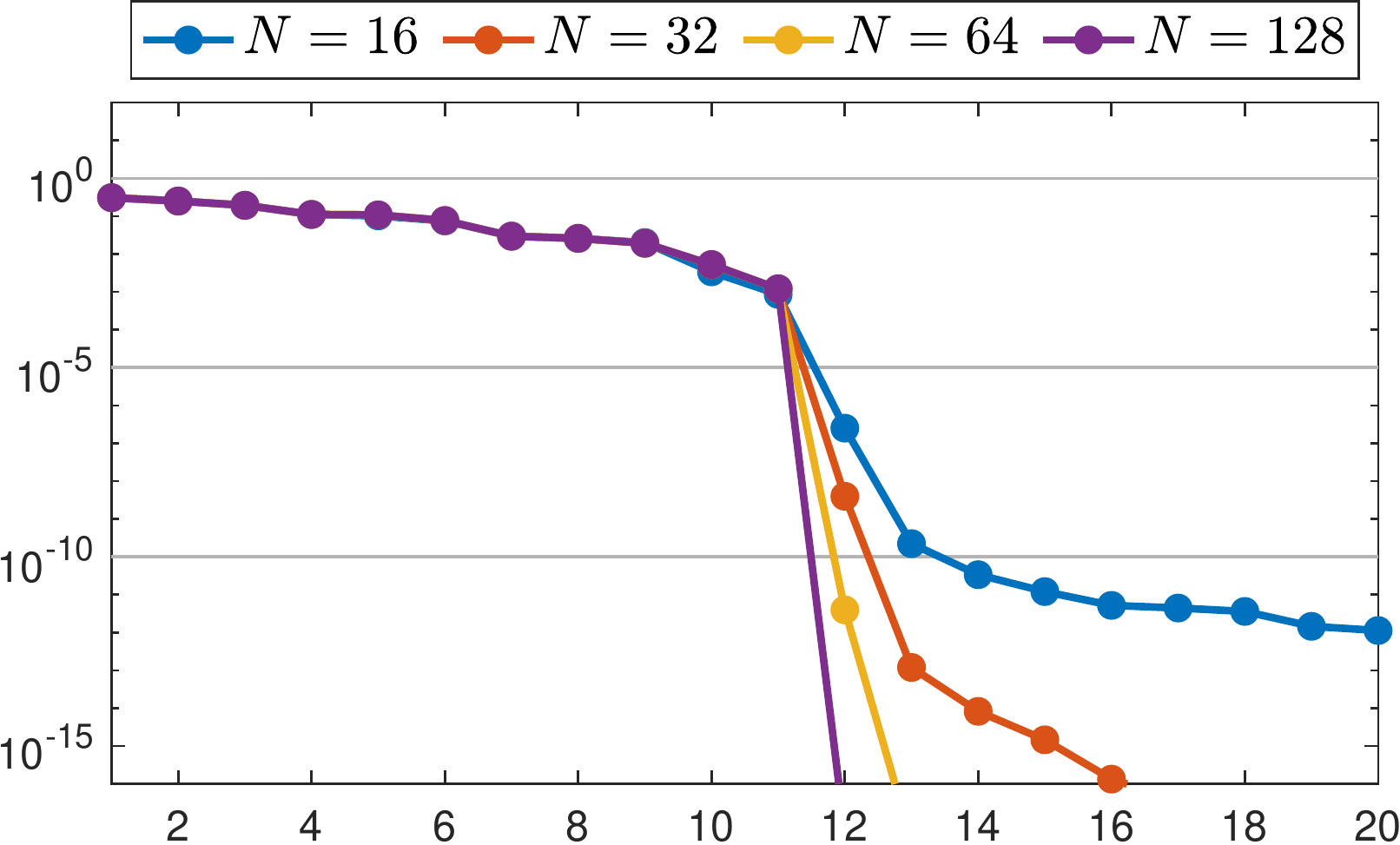}

\begin{center}
\begin{picture}(0,0)
\put(-183,18){\rotatebox{90}{\footnotesize $k$th singular value of $\H$}}
\put(7,18){\rotatebox{90}{\footnotesize $k$th singular value of $\L$}}
\put(-153,25){\footnotesize Hankel}
\put(37,35){\footnotesize single point }
\put(37,25){\footnotesize Loewner, $\sigma=10$}
\put(-32,4){\footnotesize $k$}
\put(157,4){\footnotesize $k$}
\end{picture}
\end{center}

\vspace*{-1em}
\caption{\label{fig:contour_trap}
Singular values of the block Hankel matrix $\H$ (established method; 
\Cref{sec:hankel,sec:hankelsys}) and the single-shift
block Loewner matrix $\L$ (new method; expansion about $\sigma=10$; 
\Cref{sec:singlepoint}) as a function of the number of quadrature points, $N$,
for the NLEVP in \Cref{fig:contour_sketch} using 
$\ell=r=11$ probing directions and $K=5$ blocks.
The matrices with exact data both have rank $m=11$,  
which in this case is clearer at small $N$ for the method on the right.
}
\end{figure}
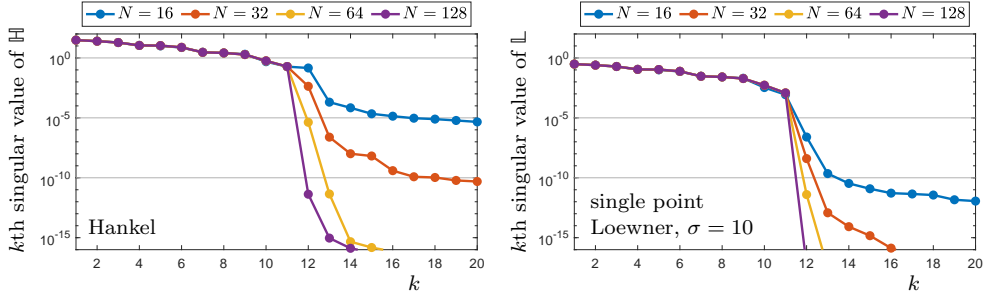
%%%%%%%%%%%%%%%%%%%%%%%%%%%%%%%%%%%%%%%%%%%%%%%%%%%%%%%%%%%%%%%%%%%%%%%%%%%%%%%%

\Cref{fig:contour_sketch} shows a typical setting.
While a user of the algorithm specifies the target region $\Omega$
(e.g., to look for rightmost eigenvalues for stability analysis), 
practical implementations should determine the number $m$ of 
eigenvalues in $\Omega$.\ \ 
In the established algorithm, $m$ is revealed by the rank 
of a block Hankel matrix composed from contour integrals 
of the form~\cref{eq:fprobe} (\Cref{sec:hankel}); 
the methods we develop in \Cref{sec:singlepoint} and~\Cref{sec:multipoint}
use block Loewner matrices.
With exact data, these methods all precisely reveal $m$;
in practice, the quadrature approximations~\cref{eq:trap}
give inexact data, so $m$ must be discerned from the 
decaying singular values of the block Hankel or Loewner matrix.  
\Cref{fig:contour_trap} illustrates this challenge using the 
Hankel approach (left) and the single-point Loewner method (\Cref{sec:singlepoint}),
applied to the problem in \Cref{fig:contour_sketch}.
We seek the $m=11$ eigenvalues in the disk $\Omega$, using $\ell=r=11$ random
probing directions and $K=5$ block rows and columns (explained in the following sections).  
With exact data these matrices both have rank $m$, but with trapezoid-rule
approximations, the singular values behave differently, and the rank $m$ 
can be more easily approximated with modest $N$ using one method or the other.

%% file: schematic.tex
\begin{picture}(0,150)(5,100)
\put(60,155){\textcolor{blue!10!white}{\rule[-3pt]{41.5pt}{96pt}}}
\put(60,155){\framebox{\makebox(35,90)[c]{$\BV$}}}
\put(107,210){\textcolor{blue!10!white}{\rule[-3pt]{41.5pt}{41pt}}}
\put(107,210){\framebox{\makebox(35,35)[c]{\,\footnotesize$(z\BI\kern-1.4pt-\kern-1.4pt\BLambda)^{-1}$}}}
\put(154,210){\textcolor{blue!10!white}{\rule[-3pt]{96.5pt}{41pt}}}
\put(154,210){\framebox{\makebox(90,35)[c]{$\BW^*$}}}
 \put(260,200){\makebox(0,0)[c]{$+$}}
\put(270,155){\textcolor{red!10!white}{\rule[-3pt]{96.5pt}{96pt}}}
\put(270,155){\framebox{\makebox(90,90)[c]{$\BN(z)$}}}
\put(35,200){\makebox(0,0)[c]{$\BT(z)^{-1} = {}$}}
\put(160,120){\makebox(0,0)[c]{\shortstack[c]{{\small$\BH(z)\coloneqq\BV(z\BI-\BLambda)^{-1}\BW^*$} \\ \emph{$n\times n$ rational function, order ${}\le m$:} \\ \emph{${}\le m$ poles in $\Omega$}}}}
\put(317,120){\makebox(0,0)[c]{\shortstack[c]{{\small $\BN(z)$} \\ \emph{$n\times n$ nonlinear function,}\\ \emph{no poles in $\Omega$}}}}
\end{picture}

%% file: hankel.tex
%!TEX root = paper.tex
\section{Algorithms based on Hankel matrices} \label{sec:hankel}
Established contour integration methods for NLEVPs are based on
Hankel matrix techniques~\cite{ASTIK09,Beyn12},
which we briefly summarize in this section. 
Consider the analytic matrix-valued function 
$\BT: \Omega \subset \C \to \Cnn$ 
that defines the NLEVP
\[
\BT(\lambda)\Bv = \Bzero.
\]
Assume that $\Omega \subset \C$ is an open domain and $\BT$ has $m$ eigenvalues $\lambda_1, \ldots, \lambda_m$ (distinctly labeling each copy of a multiple eigenvalue)
in $\Omega$.\ \ \Cref{thm:Keldysh} ensures that 
\[
\BT(z)^{-1} = \BV(z\BI-\BLambda)^{-1}\BW^* + \BN(z).
\]
With $f(z)=z^k$ and probing matrices $\BL\in\C^{n\times \ell}$ and $\BR\in\C^{n\times r}$ 
in \cref{eq:fprobe}, define
\begin{equation} \label{eq:ApLR}
\BA_k :=  {1\over 2\pi\iop} \int_{\partial \Omega} z^k\ \BLs{*}\BT(z)^{-1} \BR \,\dop z\in \C^{\ell\times r}.
\end{equation}
Using \cref{eqn:cauchy},  we can factor
\begin{equation}
\label{eqn:Ap}
\BA_k = \BLs{*}\BV \BLambdas{k} \BW^*\BR;
\end{equation}
in particular,
\begin{equation} \label{eq:A0A1}
\BA_0 = \BLs*\BV \BW^* \BR, \qquad
\BA_1 = \BLs*\BV \BLambda \BW^* \BR.
\end{equation}
At this stage, one could  analyze the (possibly singular) rectangular matrix pencil
$\BA_1 - z\BA_0 = \BLs*\BV(\BLambda - z@\BI)\BW^*\BR$ of dimension 
$\ell\times r$~\cite{BEGM05,HMP,Ste94,WT02}.
Instead, established nonlinear eigenvalue contour algorithms 
reduce $\BA_1-z\BA_0$ to a square pencil.

Note that $\rank(\BA_0)$ and $\rank(\BA_1)$ 
depend on the probing matrices $\BL$ and $\BR$,
the eigenvectors in $\BV$ and $\BW$, and 
the number $m$ of eigenvalues in $\Omega$.
Suppose for now that we have at least as many left and right
probing directions as desired eigenvalues ($\ell, r \ge m$) and 
the eigenvectors that form the columns of $\BV\in\C^{n\times m}$ 
and $\BW\in\C^{n\times m}$ are linearly independent.
Then for generic choices of 
$\BL\in\C^{n\times \ell}$ and $\BR \in \C^{n\times r}$,
\begin{equation}
\label{eq:rank_cond}
 \rank(\BLs*\BV) = \rank(\BW^*\BR) = m.
\end{equation}
In practice,  the probing matrices $\BL$ and $\BR$ 
are constructed with random entries,
so the rank conditions~\cref{eq:rank_cond} hold 
with high probability~\cite{CD05};
cf.~\cite[sect.~3]{Beyn12}.

Via Sylvester's rank inequality (see, e.g., \cite[Sect.~0.4.5]{HJ85}),
the condition~\cref{eq:rank_cond} implies
$\rank(\BA_0) = m$.
Take the (reduced) singular value decomposition (SVD)
\begin{equation} \label{eq:A0svd}
\BA_0^{} = \BX @ \BSigma @ \BY^*,
\end{equation}
with $\BX \in \C^{\ell\times m}$ and
$\BY \in \C^{r\times m}$ having orthonormal columns,
and invertible $\BSigma\in\Rmm$.
From these ingredients we can expose $\BLambda$
and extract the $m$ eigenpairs of $\BT$ in $\Omega$
by solving a standard (linear) matrix eigenvalue problem.

%%%%%%%%%%%%%%%%%%%%%%%%%%%%%%%%%%%%%%%%%%%%%%%%%%%%%%%%%%%%%%%%%%%%%%%%%%%%%%%%
\begin{theorem}
	\label{thm:few_eigs}
Assume $\Omega$ contains $m$~eigenvalues
$\lambda_1,\ldots, \lambda_m$
with linearly independent right and left eigenvectors stored
in $\BV, \BW \in\C^{n\times m}$, and the probing matrices 
$\BL \in \C^{n\times \ell}$ and $\BR \in \C^{n\times r}$ 
satisfy the rank conditions~\cref{eq:rank_cond}
(implying $\ell, r \ge m$).  
Then there exists an invertible $\BS  \in\Cmm$
such that $\BLs*\BV = \BX\BS$, 
and
	\[
	\BB := \BX^* @{\BA_1} \kern-1.5pt \BY \BSigma^{-1} 
        \in \Cmm
	\] 
can be diagonalized as $\BB=\BS \BLambda \BS^{-1}$,
with $\BLambda = {\rm diag}(\lambda_1,\ldots,\lambda_m)$ containing the 
$m$ eigenvalues of the nonlinear eigenvalue problem
$\BT(\lambda)\Bv=\Bzero$ in~$\Omega$.
\end{theorem}
%%%%%%%%%%%%%%%%%%%%%%%%%%%%%%%%%%%%%%%%%%%%%%%%%%%%%%%%%%%%%%%%%%%%%%%%%%%%%%%%

\begin{proof}
The rank conditions~\cref{eq:rank_cond} imply that
$\BLs*\BV$ and $\BA_0 = \BLs*\BV\BW^*\BR$ have the same 
column space, which hence also agrees 
with the column space of $\BX$ in~(\ref{eq:A0svd}).
Thus there exists a unique invertible 
$\BS \in \Cmm$ such that $\BLs*\BV = \BX@\BS$.\ \ 
Equate the expressions for $\BA_0$ in~\cref{eq:A0A1} and~\cref{eq:A0svd}
to see
\[\BLs*\BV\BW^*\BR = \BX@\BSigma @\BY^*.\]
Substitute $\BLs*\BV = \BX@\BS$ on the left to obtain
$\BX@\BS@\BW^*\BR = \BX \BSigma\BY^*$.
Premultiply by $\BX^*$ and postmultiply by $\BY$:
since the columns of $\BX$ and $\BY$ are orthonormal,
\begin{equation} \label{eq:Sig0}
 \BS@\BW^*\BR\BY = \BSigma.
\end{equation}
From the expression \cref{eq:A0A1} for $\BA_1$, 
we find
\begin{align*}
\BB := \BX^*{\BA_1} \kern-1pt \BY @\BSigma^{-1}
     &= \BX^*(\BLs*\BV\BLambda\BW^*\BR)\BY\BSigma^{-1} \\
     &= \BX^*(\BLs*\BV)\BLambda(\BS^{-1}\BS)(\BW^*\BR)\BY\BSigma^{-1}  \\
     &= \BX^*(\BX@\BS)\BLambda\BS^{-1}(\BS\BW^*\BR\BY)\BSigma^{-1} 
         = \BS\BLambda\BS^{-1},
\end{align*}
where the last equality follows from the expression
for $\BSigma$ in~\cref{eq:Sig0}.
\end{proof}
%%%%%%%%%%%%%%%%%%%%%%%%%%%%%%%%%%%%%%%%%%%%%%%%%%%%%%%%%%%%%%%%%%%%%%%%%%%%%%%%

The matrix $\BB$ in \cref{thm:few_eigs} does not reveal
the \emph{eigenvectors} of the 
nonlinear eigenvalue problem $\BT(\lambda)\Bv=\Bzero$.
In systems theory terms, 
\Cref{thm:few_eigs} uses the \emph{two-sided samples} $\BLs*\BV f(\BLambda)\BW^*\BR$ 
to \emph{realize} the \emph{sampled} transfer function 
$\Hsamp(z) = \BLs*\BH(z) \BR = \BLs*\BV (z\BI-\BLambda)^{-1}\BW^*\BR$
associated with the sampled dynamical system~\cref{eq:LRsys}:
the process recovers $\BLs*\BV$ and $\BW^*\BR$, instead of $\BV$ and $\BW^*$.
If the left probing is trivial, $\BL=\BI$ (so $\ell=n$), partition 
$\BS = [\Bs_1,\ldots, \Bs_m]$ by columns to reveal the eigenvectors
of the nonlinear eigenvalue problem, stored in the columns of $\BV$:
the $m$ eigenpairs associated with $\Omega$ are 
$(\lambda_j, \BX@\Bs_j)$ for $j=1,\ldots, m$.

\medskip
If the eigenvectors that form the columns of 
$\BV$ and $\BW$ are not linearly independent 
(certainly the case when $m>n$, but possible in other situations),
or if we have too few probing directions ($\ell, r < m$),
then the rank conditions~\cref{eq:rank_cond} cannot be satisfied,
regardless of $\BL$ and $\BR$.\ \ 
Such situations can be handled by forming, for suitable 
$K\ge 1$, the \emph{block Hankel matrices}
\begin{equation}  
\label{eqn:H}
	\H = \left[
	\begin{array}{cccc}
		\BA_0 & \BA_1 & \!\!\cdots\!\! & \BA_{K-1} \\
		\BA_1 & \BA_2 & \!\!\cdots\!\! & \BA_{K}\\
		\vdots & \vdots & \!\!\iddots\! & \vdots \\
		\BA_{K-1} & \BA_{K} & \!\!\cdots\!\! & \BA_{2K-2}
	\end{array}
	\right], \qquad
	\H_s = \left[
	\begin{array}{cccc}
		\BA_1 & \!\!\BA_2 & \!\!\cdots\!\! & \BA_{K} \\
		\BA_2 & \!\!\BA_3 & \!\!\cdots\!\! & \BA_{K+1}\\
		\vdots & \!\!\vdots & \!\!\iddots\!\! & \vdots \\
		\BA_{K} & \!\!\BA_{K+1} & \!\!\cdots\!\! & \BA_{2K-1}
	\end{array}
	\right]. 
\end{equation}
(We use double-struck characters like $\H$ to represent matrices
constructed from data obtained via contour integration.)
The submatrices $\BA_j\in\C^{\ell\times r}$ have the form~\cref{eqn:Ap},
so that $\H, \H_s\in \C^{\ell K\times r K}$.
(The $K = 1$ case gives $\H = \BA_0$ and $\H_s = \BA_1$.)
All our approaches will utilize convenient (theoretical)
factorizations of structured block matrices, 
$\H$ and $\H_s$ in this case.
Use the decomposition~\cref{eqn:Ap} for $\BA_k$
to factor
\begin{equation} \label{eqn:Hp_decomp} 
\H = \OO@\RR, \qquad \H_s = \OO\BLambda \RR,
\end{equation}
where $\OO \in \C^{\ell K \times m}$ 
and $\RR \in \C^{m\times r K}$ have the form
\begin{equation} \label{eqn:OORR}
\OO = \left[ 
\begin{array}{c}
\BL^*\BV \\ 
\BL^*\BV \BLambda \\
\vdots \\
\BL^*\BV \BLambda^{K-1}
\end{array} 
\right], \quad 
\RR = \left[\BW^* \BR \ \ \  \BLambda \BW^* \BR \ \ \ \cdots \ \ \ \BLambda^{K-1}\BW^* \BR \right].
\end{equation}
In systems theory, these matrices are called the 
\emph{observability} and \emph{reachability} matrices
for the sampled dynamical system~\cref{eq:LRsys}.
If $\rank(\OO) = m$, the system is \emph{observable};
if $\rank(\RR) = m$, the system is \emph{reachable}.

Assume $\BL$, $\BR$, and $K\ge 1$ have been chosen so that
\[\rank(\H) = m,\] 
which in turn implies that
$ \rank(\OO) = \rank(\RR) = m$.
(When the eigenvectors in  $\BV$ and $\BW$ are linearly independent,
it suffices to take $\ell K, r K \ge m$
for generic $\BL$ and $\BR$; cf.~\cite[sect.~5]{Beyn12}.)
Take the reduced SVD
\begin{equation}
\label{eqn:H SVD}
\H = \BX @ \BSigma@  \BY^*,
\end{equation}
where
$\BSigma = {\rm diag}(\sigma_1,\dots,\sigma_m)\in\R^{m\times m}$
contains the nonzero singular values of $\H$,
and $\BX \in \C^{\ell K \times m}$ and $\BY \in \C^{r K\times m}$
both have orthonormal columns.
Just as with the $K=1$ case covered in \cref{thm:few_eigs},
this set-up enables us to compute the $m$ eigenvalues of $\BT(z)$
in $\Omega$ by computing the eigenvalues of a matrix.

%%%%%%%%%%%%%%%%%%%%%%%%%%%%%%%%%%%%%%%%%%%%%%%%%%%%%%%%%%%%%%%%%%%%%%%%%%%%%%%%
\begin{theorem} \label{thm:many_eigs}
	With the notation above, assume $\rank(\H)=m$.
There exists an invertible $\BS\in\C^{m\times m}$ such that
$\OO = \BX\BS$, and the matrix
\[
	\BB := \BX^* @@\H_s \kern-1pt \BY \BSigma^{-1} \in \C^{m\times m}
\]
can be diagonalized as $\BB=\BS\BLambda\BS^{-1}$,
with $\BLambda = {\rm diag}(\lambda_1,\ldots,\lambda_m)$ containing the 
$m$ eigenvalues of the nonlinear eigenvalue problem
$\BT(\lambda)\Bv=\Bzero$ in~$\Omega$.
\end{theorem}
%%%%%%%%%%%%%%%%%%%%%%%%%%%%%%%%%%%%%%%%%%%%%%%%%%%%%%%%%%%%%%%%%%%%%%%%%%%%%%%%

The proof of \cref{thm:many_eigs} is an immediate generalization 
of the proof of \cref{thm:few_eigs}.
The existence of $\BS\in\Cmm$ follows from the assumption that
$\rank(\H)=m$, which, given the decomposition~\cref{eqn:Hp_decomp},
implies that $\rank(\OO)=\rank(\H)=\rank(\BX)$.

As in the $K=1$ case,
the eigenvectors of the nonlinear eigenvalue problem in $\BV$
are obscured by the probing matrix $\BL^*$.\ \  
Recovery of these eigenvectors will be among the issues 
that will be naturally resolved with our systems theory perspective.

%% file: hankelsys.tex
%!TEX root = paper.tex

%%%%%%%%%%%%%%%%%%%%%%%%%%%%%%%%%%%%%%%%%%%%%%%%%%%%%%%%%%%%%%%%%%%%%%%%%%%%%%%%
\section{Systems theory perspective on Hankel contour integration} \label{sec:hankelsys}
%%%%%%%%%%%%%%%%%%%%%%%%%%%%%%%%%%%%%%%%%%%%%%%%%%%%%%%%%%%%%%%%%%%%%%%%%%%%%%%%
As observed in the introduction,
we view the critical term $\BH(z) =  \BV(z\BI-\BLambda)^{-1}\BW^*$
in the  Keldysh decomposition \cref{eq:keldysh} for $\BT(z)^{-1}$
as the transfer function of the (unsampled)
linear dynamical system~\cref{eqn:LTI}.
Use a Neumann series to expand $\BH(z)$ around $z= \infty$:
\begin{align}
\BH(z)  = \BV(z\BI-\BLambda)^{-1}\BW^* 
       & = z^{-1}\BV\Bigl(\BI-z^{-1}\BLambda\Bigr)^{-1}\BW^* 
         \label{eqn:infseries}\\
       & = z^{-1} \BV\biggl(\sum_{k=0}^\infty z^{-k}\BLambdas{k}\biggr)\BW^*
         = \sum_{k=0}^\infty   \BV \BLambdas{k} \BW^* z^{-(k+1)}. \nonumber
\end{align}
In systems theory,
the coefficients $\BM_k := \BV \BLambdas{k} \BW^*$ of $z^{-(k+1)}$
in this  expansion are called \emph{Markov parameters}
of $\BH(z)$ \cite[Sect.~4.1]{Ant05b}.
Therefore, the contour integral~\cref{eqn:cauchy} with $f(z) = z^k$
\emph{computes the Markov parameters of $\BH(z)$}.
When left and right probing matrices are included in the integral~\cref{eq:ApLR},
the ``two-sided samples'' $\BL^*\BM_k\BR$ are the 
\emph{Markov parameters of $\Hsamp(z)$ for the sampled system}.
From the systems theory perspective, the Hankel contour integral methods
described in \cref{thm:few_eigs} and \cref{thm:many_eigs} 
address a \emph{realization problem}~\cite[Sect.~4.4]{Ant05b}.
Provided $\rank(\H) = m$, the poles of the rational function $\Hsamp(z)$ 
reveal the eigenvalues of $\BT(z)$ in $\Omega$.
%%%%%%%%%%%%%%%%%%%%%%%%%%%%%%%%%%%%%%%%%%%%%%%%%%%%%%%%%%%%%%%%%%%%%%%%%%%%%%%%

\medskip
\begin{center}
\framebox{\smallskip\ \begin{minipage}{4.7in}
\vspace*{2pt}
\textbf{Realization problem: Data at {\boldmath$z=\infty$}, two-sided samples}\\
\textsl{Given samples $\BA_k := \BL^*\BM_k\BR = \BL^*\BV \BLambdas{k} \BW^*\BR$
of the Markov parameters,
construct the transfer function $\Hsamp(z)$ for the sampled system~\cref{eq:LRsys}.}
\vspace*{2pt}
\end{minipage}\ }
\end{center}
\medskip
%%%%%%%%%%%%%%%%%%%%%%%%%%%%%%%%%%%%%%%%%%%%%%%%%%%%%%%%%%%%%%%%%%%%%%%%%%%%%%%%

\plainremark{We point out that the established eigenvalue algorithms 
from~\cite{ASTIK09,Beyn12}
in \cref{thm:few_eigs,thm:many_eigs} 
amount to 
applications of the \emph{Ho--Kalman algorithm}~\cite{DeS00,HK66} for realizing a dynamical system
from Markov parameters, in this case obtained from the contour integrals~\cref{eq:ApLR}.  
Ho and Kalman find $\BP$ and $\BQ$ such that
\[ \BP@\H@\BQ = \left[\begin{array}{cc} \BI & \Bzero \\ \Bzero & \Bzero\end{array}\right],\]
and then observe that the upper-left $m\times m$ submatrix  $(\BP\H_s\BQ)_{1:m,1:m}$ recovers $\BLambda$ 
(up to a coordinate transformation).
If $\H = \wh{\BX}\wh{\BSigma}\wh{\BY}^*$ is a full SVD,
then taking $\BP=\wh{\BX}^*$ and $\BQ = \wh{\BY}\wh{\BSigma}@^+$ will
recover
the algorithm in \cref{thm:many_eigs}, 
where $\BB = \BX^*\H_s\BY\BSigma^{-1} = (\wh{\BX}^*\H_s\wh{\BY}\wh{\BSigma}@^+)_{1:m,1:m} = (\BP\H_s\BQ)_{1:m,1:m}$,
and $\BSigma^+$ denotes the pseudoinverse of $\BSigma$.
The \emph{Silverman realization algorithm}~\cite[Sect~4.4.1]{Ant05b}, \cite{Silv71} provides
an alternative technique for realizing a system using full-rank $m\times m$ submatrices of $\H$ and $\H_s$.
}

\medskip
In this spirit, note that \emph{any system realization algorithm can be applied to solve 
eigenvalue problems}, provided one can use contour integrals like~\cref{eqn:cauchy} and~\cref{eq:fprobe}
to obtain whatever sample data is required by the realization algorithm.
In this and the following sections, we will show several ways in which the 
\emph{Loewner modeling framework} 
allows us to recover the full transfer function $\BH(z)$ using ``one-sided samples''
that only apply probing vectors on the left or right of $\BV f(\BLambda)\BW^*$.
We begin by showing how \emph{one-sided samples} of the form $\BL^*\BM_k$ and $\BM_k\BR$ can be used
to recover not just $\BLambda$, but the transfer function $\BH(z)$ of the \emph{full system}
(and hence the eigenvector matrices $\BV, \BW \in \C^{n\times m}$).
In subsequent sections we will repeat this same template, using different data about the NLEVP to
recover $\BH(z)$.

%%%%%%%%%%%%%%%%%%%%%%%%%%%%%%%%%%%%%%%%%%%%%%%%%%%%%%%%%%%%%%%%%%%%%%%%%%%%%%%%
\medskip
\begin{center}
\framebox{\smallskip\ \begin{minipage}{4.9in}
\vspace*{2pt}
\textbf{Realization problem: Data at {\boldmath$z=\infty$, one-sided samples}}\\
\textsl{Given samples $\BL^*\BM_k$ and $\BM_k\BR$ of the Markov parameters
$\BM_k := \BV \BLambdas{k}\BW^*$ along the directions $\BL$ and $\BR$, 
construct the transfer function $\BH(z)$ for the full system~\cref{eqn:LTI}.}

\vspace*{2pt}
\end{minipage}\ }
\end{center}
\medskip
%%%%%%%%%%%%%%%%%%%%%%%%%%%%%%%%%%%%%%%%%%%%%%%%%%%%%%%%%%%%%%%%%%%%%%%%%%%%%%%%

Start with the factorizations of the Hankel and shifted Hankel matrices in~\cref{eqn:Hp_decomp}:
\[ \H = \OO@\RR, \quad \H_s = \OO\BLambda \RR \in \C^{\ell K \times r K}\]
for $\OO\in \C^{\ell K\times m}$ and $\RR\in\C^{m\times r K}$ given in~\cref{eqn:OORR}.
As before, assume that $\H$ has rank~$m$, with reduced SVD 
$\H = \BX\BSigma\BY^*$.
The $m\times m$ matrices $\BX^*\OO$  and $\RR@\BY$ must also have rank~$m$,
since $(\BX^*\OO)(\RR@\BY) = \BX^*\H\BY = \BSigma$ has rank~$m$.

Recall now the \emph{full} dynamical system~\cref{eqn:LTI}.  Since $\RR@\BY$ is invertible,
we can change the coordinate system for the state variable $\Bx(t)$ in~\cref{eqn:LTI} to 
$\wh{\Bx}(t) := (\RR@\BY)^{-1}\Bx(t)$, giving the equivalent system
\begin{align*}
(\RR@\BY)\, \wh{\Bx}@'(t) &= \BLambda(\RR@\BY)\, \wh{\Bx}(t) + \BW^* \Bu(t) \\[.25em]
\By(t) &= \BV (\RR@\BY)\,\wh{\Bx}(t). 
\end{align*}
Premultiply the first equation by the invertible matrix $\BX^*\OO$ to get
another equivalent system:
\begin{align}
(\BX^*\OO)(\RR@\BY)\, \wh{\Bx}@'(t) 
         &= (\BX^*\OO) \BLambda(\RR@\BY)\, \wh{\Bx}(t) 
          + (\BX^*\OO) \BW^* \Bu(t) \label{eqn:LTItrans}\\[.25em]
\By(t) &= \BV (\RR@\BY)\,\wh{\Bx}(t).  \nonumber
\end{align}
Note that $\BX^*\OO\RR@\BY = \BX^*\H@\BY = \BSigma$ 
and $\BX^*\OO\BLambda\RR@\BY = \BX^*\H_s\BY$.  
Collect the one-sided data $\{\BL^*\BM_k\}_{k=0}^{K-1}$ and $\{\BM_k\BR\}_{k=0}^{K-1}$ into matrices
\begin{subequations} \label{eqn:BC}
\begin{align} 
   \bB :=& \left[\! \begin{array}{c}
   \BL^{*}\BM_{0} \\
   \BL^{*}\BM_{1} \\
   \vdots \\
   \BL^{*}\BM_{K-1}
   \end{array}\!\right] = 
 \left[\! \begin{array}{c}
   \BL^{*}\BV \\
   \BL^{*}\BV\BLambda \\
   \vdots \\
   \BL^{*}\BV\BLambda^{K-1}
   \end{array}\!\right] \BW^*  = \OO\BW^*, \\[.75em]
\bC :=& \left[\,\BM_{0} \BR ~~\BM_{1} \BR ~~ \cdots ~~ \BM_{K-1} \BR \,\right]  \\
     =& \left[\,\BV \BW^* \BR ~~ \BV\BLambda\BW^*\BR ~~ \cdots ~~ \BV\BLambda^{K-1}\BW^*\BR \,\right] 
= \BV\RR.  \nonumber
\end{align}
\end{subequations}
With this notation, the system~\cref{eqn:LTItrans} becomes
\begin{align}
\BSigma\, \wh{\Bx}@'(t) 
         &= (\BX^*\H_s\BY)\, \wh{\Bx}(t) + (\BX^*\bB)\,\Bu(t) \label{eqn:LTIinf}\\[.25em]
\By(t) &= (\bC\BY)\,\wh{\Bx}(t).  \nonumber
\end{align}
Since the full system~\cref{eqn:LTI} and this transformed system~\cref{eqn:LTIinf}
only differ in the coordinate system for the internal state variables, they describe the
same input--output map, and hence share the common transfer function
\begin{align*}
   \BH(z) &= \BV (z\BI-\BLambda)^{-1}\BW^* \\[.25em]
          &= \bC\BY(z \BSigma - \BX^*\H_s\BY)^{-1}\BX^*\bB.
\end{align*}
The eigenvalues of the NLEVP in $\Omega$, which are the poles of $\BH(z)$,
are thus eigenvalues of the $m \times m$ pencil $(\BX^*\H_s\BY,\BSigma)$,
and so the eigenvalues of $\BX^*\H_s\BY\BSigma^{-1}$,
the key matrix in~\cref{thm:many_eigs}.
Moreover, since this pencil corresponds to an equivalent realization
of $\BH(z)$, its eigenvalues have the same multiplicities and
indices as those of $\BLambda$; since we assumed the eigenvalues of 
the NLEVP in $\Omega$ are semi-simple, so too must be the eigenvalues
of  $(\BX^*\H_s\BY,\BSigma)$.
If we diagonalize $\BX^*\H_s\BY\BSigma^{-1} = \BS\BLambda\BS^{-1}$, then
\begin{align*}
 \BV(z\BI-\BLambda)^{-1}\BW^* 
      &= \bC\BY(z \BSigma - \BX^*\H_s\BY)^{-1}\BX^*\bB \\[.25em]
      &= (\bC\BY\BSigma^{-1}\BS) (z \BI - \BLambda)^{-1}(\BS^{-1}\BX^*\bB).
\end{align*}
We thus identify the $j$th column of $\bC\BY\BSigma^{-1}\BS$ and $j$th row of $\BS^{-1}\BX^*\bB$
as right and left eigenvectors of $\BT(z)$ associated with $\lambda_j$. 
Using one-sided data permits recovery of the full system, and hence recovery of 
eigenvalues and eigenvectors.
(Though we restrict to semi-simple eigenvalues, 
these ideas directly extend to defective eigenvalues and Jordan blocks,
since the construction obtains an equivalent realization of $\BH(z)$.)

\plainremark{Two steps led to the transformed dynamical
system in~\cref{eqn:LTItrans}: (a)~the change of variables $\Bx(t) = \RR\BY\wh{\Bx}(t)$,
and (b)~the premultiplication by $\BX^*\OO$.\ \ 
These steps parallel projection-based model reduction algorithms;
see, e.g., \cite[sect.~3.2]{AntBG20}.
In that context, step~(a) is replaced by an approximation that restricts the state vector
to a lower-dimensional subspace, and step~(b) imposes a Petrov--Galerkin
condition to close the system.  In the NLEVP setting, the number $m$ of eigenvalues in 
$\Omega$ will often be unknown, to be estimated from the singular values of a block
Hankel matrix $\H$ populated with inexact (quadrature) data.
If $m$ is \emph{underestimated}, the process described in this section
develops a reduced-order model of $\BH(z)$,
and the computed eigenvalues will only be approximations to the
true eigenvalues of $\BT(z)$ in $\Omega$.
}

%% file: singlept.tex
%!TEX root = paper.tex

\section{Single-point Loewner algorithm for NLEVPs} \label{sec:singlepoint}

The Hankel techniques described in the last two sections use samples of
the Markov parameters $\BV\BLambdas{k}\BW^*$, the coefficients in the series 
expansion of $\BH(z)$ at $z=\infty$~\cref{eqn:infseries}.
What if we instead expand
$\BH(z)$ about some other point $\rp\in\C$ outside the domain $\Omega@$?\ \ 
With exact data (contour integrals), this approach yields an equivalent 
realization of $\BH(z) = \BV(z\BI-\BLambda)^{-1}\BW^*$;
for inexact data (quadrature), the accuracy of the computed eigenvalues will differ.
As we will show in this section, 
the most computationally-intensive work required in this new approach 
is precisely the same as for the Hankel case: computing the terms 
$\BL^*\BT(\qn_k)^{-1}\BR$ in the quadrature approximations~\cref{eq:trap}
to the exact integrals. Therefore, one can explore
this method (indeed, with multiple values of $\rp$) at little extra cost.

For $\rp \notin \overline{\Omega}$, expand $\BH(z)$ in the Taylor series
\[ \BH(z) =  \BV(z\BI-\BLambda)^{-1}\BW^*  
= \sum_{k=0}^\infty \Big({1\over k!} \BH^{(k)}(\rp)\Big) (z-\rp)^k
=: \sum_{k=0}^\infty \BM_{k} (z-\rp)^k, \]
where we now define
\begin{equation} \label{eq:Hderiv}
 \BM_k := {1\over k!} \BH^{(k)}(\sigma) = (-1)^k \, \BV(\sigma\BI-\BLambda)^{-(k+1)}\BW^*.
\end{equation}
The coefficients $\BM_k\in\Cnn$ can be obtained by
evaluating the integrals~\cref{eqn:cauchy} and~\cref{eq:ApLR} with 
$f(z) = (-1)^k/(\rp-z)^{k+1}$
at essentially the same expense as the $f(z) = z^k$ used in the Hankel method;
\Cref{sec:contsinglepoint} provides details.
The matrices $\BM_k$ are \emph{moments of the transfer function $\BH(z)$ about $z=\rp$}; 
see, e.g., \cite[p.~109]{Ant05b}. 
Akin to the last section, we
seek to recover the full transfer function $\BH(z)$ from
one-sided (``tangential'') measurements $\BL^*\BM_k$ and $\BM_k\BR$ of $\BM_k$, 
given 
$\BL \in \C^{n \times \ell}$ and $\BR \in \C^{n \times r}$ ($\ell, r \le n$).

%%%%%%%%%%%%%%%%%%%%%%%%%%%%%%%%%%%%%%%%%%%%%%%%%%%%%%%%%%%%%%%%%%%%%%%%%%%%%%%%
\medskip
\begin{center}
\framebox{\smallskip\ \begin{minipage}{4.8in}
\vspace*{2pt}
\textbf{Realization problem: Data at {\boldmath$z=\rp\not \in \overline{\Omega}$, one-sided samples}}\\
\textsl{Given samples of the moments $\BM_k := (-1)^k \BV (\sigma\BI-\BLambda)^{-(k+1)} \BW^*$
along the directions $\BL$ and $\BR$, i.e., $\BL^*\BM_k$ and $\BM_k\BR$,
construct the transfer function $\BH(z)$ for the full system~\cref{eqn:LTI}.}

\vspace*{2pt}
\end{minipage}\ }
\end{center}
\medskip
%%%%%%%%%%%%%%%%%%%%%%%%%%%%%%%%%%%%%%%%%%%%%%%%%%%%%%%%%%%%%%%%%%%%%%%%%%%%%%%%
Inspired by the rational bi-tangential interpolation model reduction 
framework of Mayo and Antoulas~\cite[sect.~6]{MA07},
construct the $\ell K \times r K$ matrices 
\begin{align} 
\L &=
     \left[
	\begin{array}{cccc}
		\BLs*\BM_{1} \BR& \BLs*\BM_{2}\BR & \cdots & \BLs*\BM_{K} \BR\\
	\BLs*\BM_{2}\BR & \BLs*\BM_{3} \BR& \cdots & \BLs*\BM_{K+1}\BR\\
		\vdots & \vdots & \iddots & \vdots \\
		\BLs*\BM_{K}\BR & \BLs*\BM_{K+1} \BR& \cdots & \BLs*\BM_{2K-1}\BR
	\end{array}
	\right],   \label{eqn:L} \\[10pt]
	\L_0 &=  
\left[
	\begin{array}{cccc}
		\BL^*\BM_{0}\BR & \BL^*\BM_{1}\BR & \cdots & \BL^*\BM_{K-1} \BR\\
		\BL^*\BM_{1}\BR & \BL^*\BM_{2}\BR & \cdots & \BL^*\BM_{K}\BR\\
		\vdots & \vdots & \iddots & \vdots \\
		\BL^*\BM_{K-1}\BR & \BL^*\BM_{K}\BR & \cdots & \BL^*\BM_{2K-2}\BR
	\end{array}
	\right],  
	\label{eqn:L0} \\[10pt]
\L_s &= \rp@\L + \L_0.
	\label{eqn:Ls}
\end{align}
Following Mayo and Antoulas, we call $\L$ and $\L_s$ 
the \emph{Loewner and shifted Loewner matrices corresponding to $\BH(z)$ 
for the single point $\rp$}. 
Analogous to the Hankel development, define
\begin{align*}
 \OO &:= \left[\begin{array}{c}
                 \BL^*\BV(\sigma\BI-\BLambda)^{-1}  \\
                -\BL^*\BV(\sigma\BI-\BLambda)^{-2}  \\
                \vdots  \\
                (-1)^{K+1}\BL^*\BV(\sigma\BI-\BLambda)^{-K}  
         \end{array}\right] \in \C^{\ell K\times m}, \\[.5em]
 \RR &:= \left[\begin{array}{ccc}
               (\sigma\BI-\BLambda)^{-1}\BW^*\BR 
                & \cdots
                & (-1)^{K+1}(\sigma\BI-\BLambda)^{-K}\BW^*\BR 
         \end{array}\right] \in \C^{m\times rK}.
\end{align*}
We do not have access to these matrices, as they require 
the eigenvalues and eigenvectors;
they merely support the derivation.
In particular, we can express
\begin{subequations} \label{eq:LoewnerOR}
\begin{align}
  \L_0 &= \OO(\sigma\BI-\BLambda)\RR, \\[.25em]
  \L_{\phantom{0}}   &= -\OO\RR, \\[.25em]
  \L_s &= \sigma\L + \L_0 = -\OO\BLambda \RR.
\end{align}
\end{subequations}
Collect the left and right samples in
\begin{subequations} \label{eq:YZ}
\begin{align} 
	\bB :=& \left[\! \begin{array}{c}
	\BL^{*}\BM_{0} \\
	\BL^{*}\BM_{1} \\
	\vdots \\
	\BL^{*}\BM_{K-1} 
	\end{array}\! \right] \,=\, \OO\BW^* \in \C^{\ell K  \times n}, \\[.75em]
\bC :=& \left[\,\BM_{0} \BR ~~\BM_{1} \BR ~~ \cdots ~~ \BM_{K-1} \BR \,\right] 
    \,=\, \BV\RR \in \C^{n \times r K}.
\end{align}
\end{subequations}
Now compute the reduced (rank-$m$) SVD
\begin{equation} \label{eqn:shiftsvd}
 \L = \BX \BSigma\BY^*,
\end{equation}
with $\BX \in  \C^{\ell K \times m}$, $\BSigma\in\R^{m\times m}$,
and $\BY\in \C^{r K \times m}$.
Since
\begin{equation} \label{eq:SigXORY}
  \BSigma = \BX^*\L\BY = -\BX^*\OO\RR\BY
\end{equation}
has rank~$m$, the $m\times m$ matrices $\BX^*\OO$ and $\RR\BY$ are invertible.

Just as in the last section, we recover the full transfer function $\BH(z)$
by appropriately transforming the original system~\cref{eqn:LTI}.
Change variables according to $\Bx(t) =: \RR\BY\wh{\Bx}(t)$, giving
\begin{align*}
   \RR\BY@\wh{\Bx}'(t) &= \BLambda \RR\BY@\wh{\Bx}(t) + \BW^* \Bu(t) \\[.25em]
                  \By(t) &= \BV\RR\BY@\wh{\Bx}(t).
\end{align*}
Premultiply the first equation by $\BX^*\OO$ to get
\begin{align*}
   \BX^*\OO\RR\BY@\wh{\Bx}'(t) &= \BX^*\OO\BLambda \RR\BY@\wh{\Bx}(t) + \BX^*\OO\BW^* \Bu(t) \\[.25em]
                  \By(t) &= \BV\RR\BY@\wh{\Bx}(t).
\end{align*}
Identifying, via~\cref{eq:LoewnerOR}, \cref{eq:YZ}, and \cref{eq:SigXORY},
\[ 
\BX^*\OO\RR\BY = -\BSigma, 
\quad
\BX^*\OO\BLambda\RR\BY = -\BX^*\L_s\BY, 
\quad
\BX^*\OO\BW^* = \BX^*\bB,
\quad
\BV\RR\BY = \bC\BY,\]
we build an equivalent realization of the transfer function, using only data matrices
and the SVD~\cref{eqn:shiftsvd} of the data matrix $\L$:
\begin{align} 
 \BH(z) &= \big(\BV\RR\BY\big)\big(z(\BX^*\OO\RR\BY) - (\BX^*\OO\BLambda \RR\BY)\big)^{-1}\big(\BX^*\OO\BW^*\big) \nonumber \\[.25em]
        &= \bC\BY(\BX^*\L_s\BY - z@\BSigma)^{-1} \BX^*\bB.  \label{eqn:Hwithmu}
\end{align}

%%%%%%%%%%%%%%%%%%%%%%%%%%%%%%%%%%%%%%%%%%%%%%%%%%%%%%%%%%%%%%%%%%%%%%%%%%%%%%%%
\begin{theorem}  	\label{thm:singleshift}
Suppose $\BT(z)$ has $m$ eigenvalues in $\Omega$, all semi-simple, and suppose $\rp\in\C$ is different from these eigenvalues.
Let $\L$ and $\L_s$ be given as in \cref{eqn:L} and \cref{eqn:Ls}.
Suppose ${\rm rank}(\L)=m$, and take the reduced  SVD
   $\L = \BX\BSigma\BY^*$.
	Define $\BB_\rp\in \C^{m \times m}$ and its eigenvalue decomposition:
	\begin{equation} \label{eqn:Bsigma}
	\BB_\rp :=   \BSigma^{-1} \BX^* \L_s\BY = \BS \BLambda\BS^{-1}.
	\end{equation}
The matrix $\BLambda = {\rm diag}(\lambda_1,\ldots, \lambda_m)$ reveals the 
$m$ eigenvalues of $\BT(z)$ in $\Omega$.\ Let $\Bs_j$ denote
the $j$th column of $\BS$.  Then $(\lambda_j,  \bC \BY \Bs_j)$ is an eigenpair of $\BT(z)$.
\end{theorem}
\marginpar[\raggedright Must $\BB_\sigma$ be diagonalizable?]{}
If we only seek eigenvalues (and not eigenvectors), we can 
construct $\BB_\rp$ entirely from two-sided samples of the form $\BL^*\BM_k\BR$,
analogous to the construction in \cref{thm:few_eigs}.

%%%%%%%%%%%%%%%%%%%%%%%%%%%%%%%%%%%%%%%%%%%%%%%%%%%%%%%%%%%%%%%%%%%%%%%%%%%%%%%%

%%%%%%%%%%%%%%%%%%%%%%%%%%%%%%%%%%%%%%%%%%%%%%%%%%%%%%%%%%%%%%%%%%%%%%%%%%%%%%%%
\subsection{Computing the moments $\BM_k$}  \label{sec:contsinglepoint}

The realization algorithm just described requires the
tangentially probed {one-sided samples $\BL^*\BM_k$ and $\BM_k\BR$.}\ \ 
\Cref{thm:rationalfunctionCI} shows how $\BM_k$ 
can be computed via a contour integral.
To compute the probed matrices {$\BL^*\BM_k$  and $\BM_k\BR$} directly, replace 
$\BT(z)^{-1}$ in the integrand with {$\BL^*\BT(z)^{-1}$ and $\BT(z)^{-1}\BR$.}
When these integrals are computed via quadrature, each integrand
evaluation then amounts to solving a linear system involving 
$\BT(\qn_k)$ at each quadrature point $\qn_k$; see \Cref{section:filterfunctions}.

%%%%%%%%%%%%%%%%%%%%%%%%%%%%%%%%%%%%%%%%%%%%%%%%%%%%%%%%%%%%%%%%%%%%%%%%%%%%%%%%
 \begin{theorem}
	\label{thm:rationalfunctionCI}	
	Suppose $\BT(z)$ has a finite number of eigenvalues in the domain $\Omega$,
   all semi-simple,
   and let $\BT(z)^{-1} = \BH(z) + \BN(z)$ as in \Cref{thm:Keldysh} and \Cref{fig:schematic}.
   Moreover, assume no other eigenvalues of $\BT(z)$ lie on the boundary $\partial\Omega$.
	\begin{enumerate}[{\rm (a)}]
		\item If ${\rp\not\in \overline{\Omega}}$, then
		\begin{equation}  \label{eqn:Hviacontour}
		\frac{1}{2\pi@\iop} \int_{\partial \Omega} \frac{1}{(\rp-z)^{k+1}} \BT(z)^{-1}\, \dop z =
	 \frac{(-1)^k}{k!} \BH^{(k)}(\rp) = 
		 (-1)^k\,\BM_k.
		\end{equation}
		\item If {$\rp\in \Omega$} but $\rp$ is not a pole of $\BH(z)$, then
		\begin{equation} \label{eqn:Rviacontour}
		\frac{1}{2\pi@\iop} \int_{\partial \Omega} \frac{1}{(\rp- z)^{k+1}} \BT(z)^{-1} \, \dop z = 
		(-1)^{k+1}\frac{\BN^{(k)}(\rp)}{k!}.
		\end{equation}
	\end{enumerate}
\end{theorem}
%%%%%%%%%%%%%%%%%%%%%%%%%%%%%%%%%%%%%%%%%%%%%%%%%%%%%%%%%%%%%%%%%%%%%%%%%%%%%%%%

As throughout, the assumptions about simplicity of the eigenvalues 
are just a convenience to avoid technicalities with the Jordan form,
for the sake of presentation.

%%%%%%%%%%%%%%%%%%%%%%%%%%%%%%%%%%%%%%%%%%%%%%%%%%%%%%%%%%%%%%%%%%%%%%%%%%%%%%%%
\begin{proof}
Recall from~\cref{eqn:cauchy} that, for any $f$ analytic on $\overline{\Omega}$,
\[
{1\over 2\pi@\iop} \int_{\partial \Omega} f(z) \BT(z)^{-1}\,\dop z
= \BV f(\BLambda)\BW^*.
\]
If $\rp\not\in\overline{\Omega}$, then $f(z) = (\rp-z)^{-(k+1)}$ is analytic on
$\overline{\Omega}$ and the formula~(\ref{eqn:Hviacontour}) follows from
\[ 
{1\over 2\pi@\iop} \int_{\partial \Omega} {1\over (\sigma-z)^{k+1}} \BT(z)^{-1}\,\dop z 
  = \BV (\sigma\BI-\BLambda)^{-(k+1)} \BW^* 
  = {(-1)^k \over k!} \BH^{(k)}(\sigma) = (-1)^k \BM_k, 
\]
using the formula~\cref{eq:Hderiv} for $\BH^{(k)}(\rp)$.

If {$\rp\in \Omega$}, then $f(z) = (\rp-z)^{-(k+1)}$ is not analytic in $\Omega$; in fact, this $f$ adds a pole of order $k+1$ to $\BT(z)^{-1}$ at $\rp$
(which is, by assumption, distinct from the poles of $\BH(z)$). 
Then, by Keldysh's Theorem (\Cref{thm:Keldysh}),
\begin{eqnarray*}
{1\over 2\pi@\iop} \int_{\partial \Omega} \frac{1}{(\rp-z)^{k+1}}&&\kern-12pt \BT(z)^{-1}\, \dop z \\[.5em]
&=& \underbrace{{1\over 2\pi@\iop} \int_{\partial \Omega} \frac{1}{(\rp-z)^{k+1}} \BH(z) \,\dop z}_{I_1} + \underbrace{{1\over 2\pi@\iop}\int_{\partial \Omega} \frac{1}{(\rp-z)^{k+1}} \BN(z) \,\dop z}_{I_2}.
\end{eqnarray*}
Use the residue theorem to express the integral $I_1$ as the sum of the residue at $\rp$ and the residues at the eigenvalues $\lambda_1, \ldots, \lambda_m$, yielding
\[
I_1 = \Res\left(\frac{1}{(\rp- z)^{k+1}} \BH(z),\,\rp\!\right) 
   + \sum_{j = 1}^{m} \Res\left(\frac{1}{(\rp- z)^{k+1}} \BH(z),\,\lambda_j\!\right).
\]
Since $(\rp-z)^{-(k+1)}$ introduces a new pole of order $k+1$ at $\rp$,
\begin{eqnarray*}
\Res\left(\frac{1}{(\rp-z)^{k+1}} \BH(z),\, \rp\!\right) 
 &=& \frac{1}{k!}\lim_{z \to \rp} \dfrac{\dop^{k}}{\dop z^k}\left(\frac{(z - \rp)^{k+1}}{(\rp-z)^{k+1}} \BH(z) \right) \\[.25em]
&=& \frac{(-1)^{k+1}}{k!}\,\BH^{(k)}(\rp) \\[.25em]
&=& \frac{(-1)^{k+1}}{k!} \left((-1)^k\,k!\,\BV(\rp\BI-\BLambda)^{-(k+1)} \BW^*\right) \\[.25em]
&=& -\sum_{j=1}^{m} \frac{\Bv_j^{}\Bw_j^*}{(\rp - \lambda_j)^{k+1}}.
\end{eqnarray*}
Since for $j=1,\ldots, m$, the eigenvalue $\lambda_j$ is semi-simple, 
it is also a simple pole of $\BH(z)$, and so
\[
\Res\left(\frac{1}{(\rp-z)^{k+1}} \BH(z),\,\lambda_j\!\right) 
   = \lim_{z\to\lambda_j} {z-\lambda_j \over (\rp-z)^{k+1}} \BH(z) 
   = \frac{\Bv_j^{} \Bw_j^*}{(\rp-\lambda_j)^{k+1}},
\]
and so we conclude that the integral $I_1$ is zero:
\[
I_1 = -\sum_{j=1}^{m}\frac{\Bv_j^{} \Bw_j^*}{(\rp-\lambda_j)^{k+1}} + \sum_{j = 1}^{m}\frac{\Bv_j^{} \Bw_j^*}{(\rp-\lambda_j)^{k+1}} = 0.
\]
For the second integral $I_2$, the Cauchy integral formula gives
\[
I_2 = (-1)^{k+1}\frac{\BN^{(k)}(\rp)}{k!},
\]
thus yielding the formula~\eqref{eqn:Rviacontour}.
\end{proof}
%%%%%%%%%%%%%%%%%%%%%%%%%%%%%%%%%%%%%%%%%%%%%%%%%%%%%%%%%%%%%%%%%%%%%%%%%%%%%%%%

\begin{plainremark}
\Cref{thm:rationalfunctionCI} has two implications.
\begin{romannum}
\item  For {$\rp\notin \overline{\Omega}$}, the contour integral in \cref{eqn:Hviacontour} evaluates $\BH(\rp)$ (and its derivatives) using only $\BT(z)^{-1}$, 
enabling computation of the moments $\BM_k$~\cref{eq:Hderiv} needed 
for the single-point Loewner method discussed above, 
and the multi-point Loewner method we discuss in the next section.
\item Even when {$\rp\in \Omega$}, one could use \cref{eqn:Rviacontour} with $k=0$ to evaluate $\BH(z)$ via the expression $\BH(\rp) = \BT(\rp)^{-1} - \BN(\rp)$, further extending the applicability of contour integral methods. For the multi-point Loewner algorithm, this result enables sampling at multiple points inside the contour $\Omega$. In the present work we pursue the more conventional approach of taking sampling points that are outside $\overline{\Omega}$.
\end{romannum}
\end{plainremark}
	
%%%%%%%%%%%%%%%%%%%%%%%%%%%%%%%%%%%%%%%%%%%%%%%%%%%%%%%%%%%%%%%%%%%%%%%%%%%%%%%%
\subsection{Numerical illustration of the single-point Loewner algorithm}
\label{sec:singlept:numerics}

We now provide numerical results that compare the Hankel contour integration approach
of \Cref{sec:hankelsys} to the single-point Loewner method introduced in this section.

\begin{figure}[h!]
	\begin{center}
		\includegraphics[width=.48\textwidth]{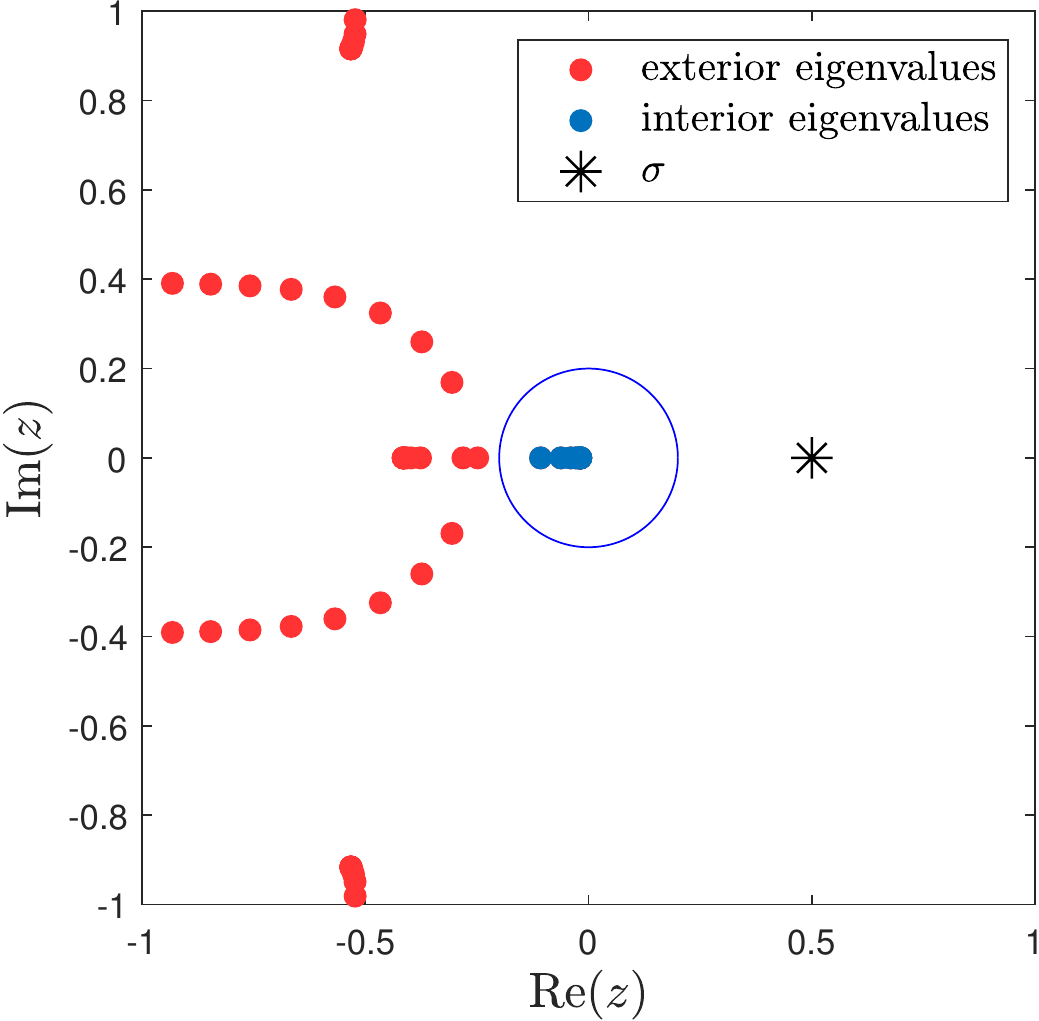}
      \quad
		\includegraphics[width=.48\textwidth]{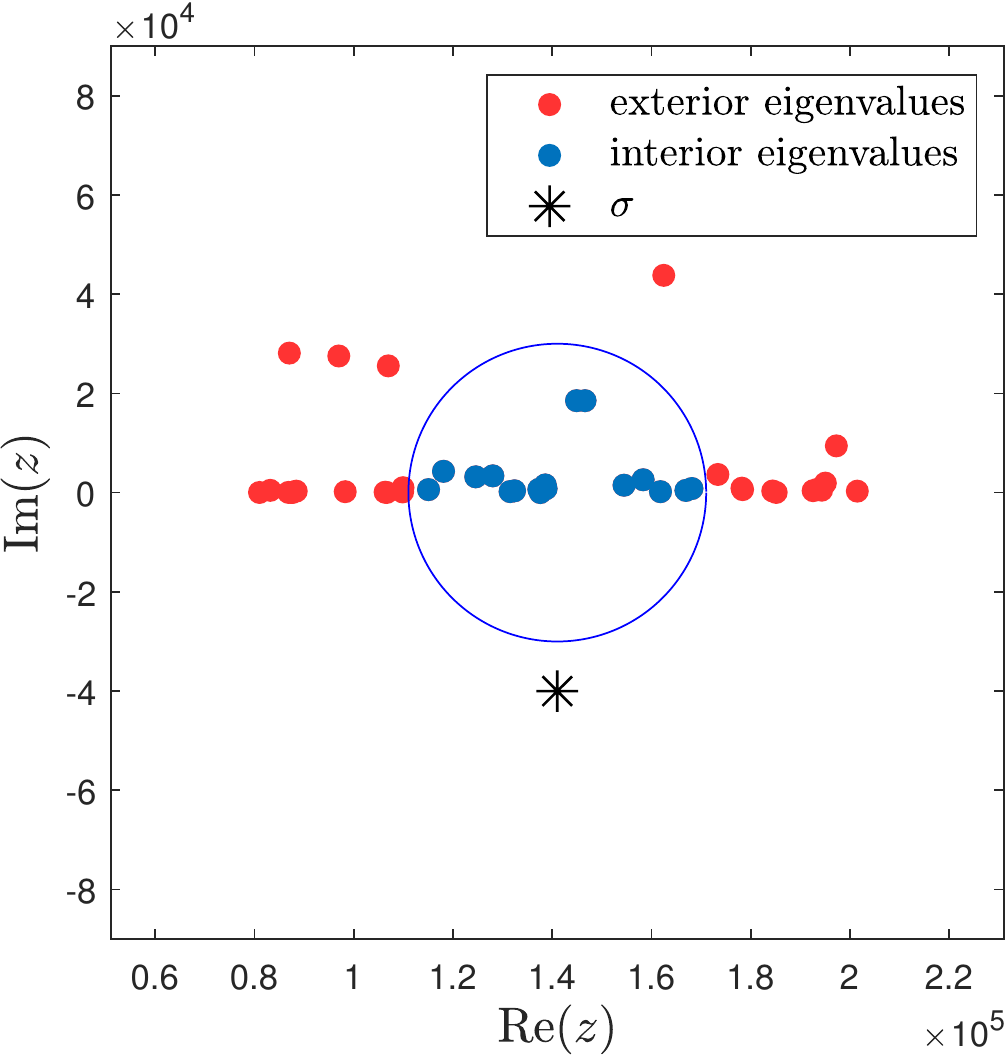}

\begin{picture}(0,0)
\put(-73,90){\footnotesize $\partial \Omega$}
\put(122,84){\footnotesize $\partial \Omega$}
\end{picture}
	\end{center}
	
	\vspace*{-2em}
	\caption{\label{fig:contour}
		The regions of interest for the delay problem (left) and {\tt gun} problem (right), 
including the contour $\partial\Omega$, eigenvalues of $\BT(z)$, 
and the interpolation point $\sigma$ $(*)$.
	}
\end{figure}
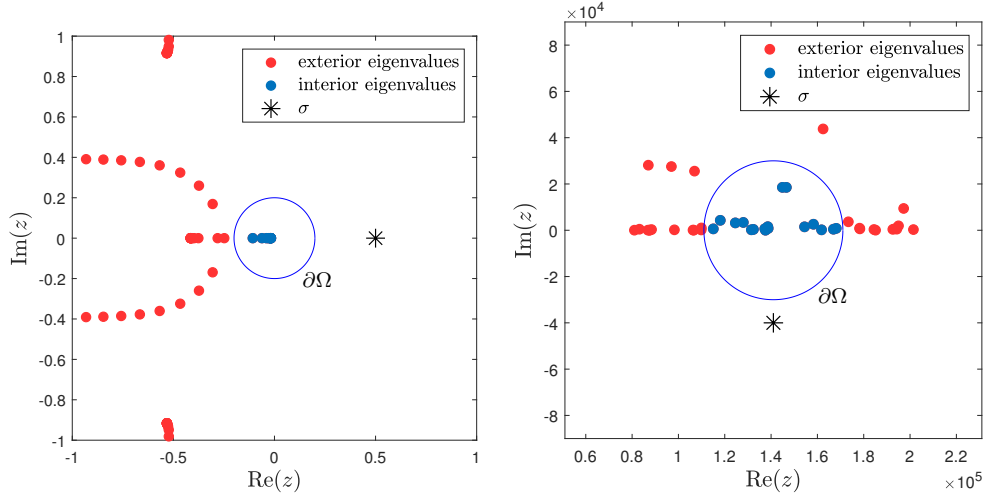

We first consider the eigenvalue problem resulting from the delay differential equation 
introduced in \Cref{fig:contour_sketch}. This NLEVP has the form
\[
\BT(z) = z\BI + c\exp(-\tau z)\BI - \BE_0,
\]
where we take the constant $c = 0.015$ and the delay length $\tau = 8$. The eigenvalues of the diagonal matrix $\BE_0\in \R^{50 \times 50}$ are chosen to be $50$ logarithmically spaced points between $-10^{10}$ and $-10^{-4}$. Given the spectrum of $\BE_0$, one can compute the eigenvalues of $\BT$ using the Lambert-W function.  (See \cite{MN14} for a comprehensive overview of delay systems,
including analysis of their corresponding eigenvalue problems.) This access to the true eigenvalues allows us to compare the performance of each method in terms of the maximum eigenvalue error, 
\begin{align}  \label{err1}
\displaystyle \max_{j=1,\ldots,m}|\lambda_j - \widetilde{\lambda}_j|,
\end{align}
where $\lambda_j$ denotes a true eigenvalue and $\widetilde{\lambda}_j$ denotes its approximation. We also consider the residual error
\begin{align} \label{err2}
 \max_{j=1,\ldots,m}\|\BT(\widetilde{\lambda}_j)\wt{\Bv}_j\|_2,
\end{align}
which measures the accuracy of the approximated eigenvector $\wt{\Bv}_j$, as well as the eigenvalue. In all cases the eigenvector is normalized so that $\|\wt{\Bv}_j\|_2 = 1$ for all $j = 1,\ldots,m$. \Cref{fig:contour} shows some of the infinitely many eigenvalues of $\BT$, the contour $\Omega$, and one choice of interpolation point, $\sigma = 1/2$. The accumulation of the eigenvalues of $\BE_0$ near the origin causes a similar concentration of the eigenvalues of $\BT$.\ \  In the context of stability analysis, one must compute the rightmost eigenvalues to good accuracy, to avoid incorrectly classifying the system as unstable.

%%%%%%%%%%%%%%%%%%%%%%%%%%%%%%%%%%%%%%%%%%%%%%%%%%%%%%%%%%%%%%%%%%%%%%%%%%%%%%%%
%%%%%%%%%%%%%%%%%%%%%%%%%%%%%%%%%%%%%%%%%%%%%%%%%%%%%%%%%%%%%%%%%%%%%%%%%%%%%%%%

%%%%%%%%%%%%%%%%%%%%%%%%%%%%%%%%%%%%%%%%%%%%%%%%%%%%%%%%%%%%%%%%%%%%%%%%%%%%%%%%
\begin{figure}[t!]
	\begin{center}
		$\overset{K = 1, \ r = 11}{\includegraphics[width=0.32\textwidth]{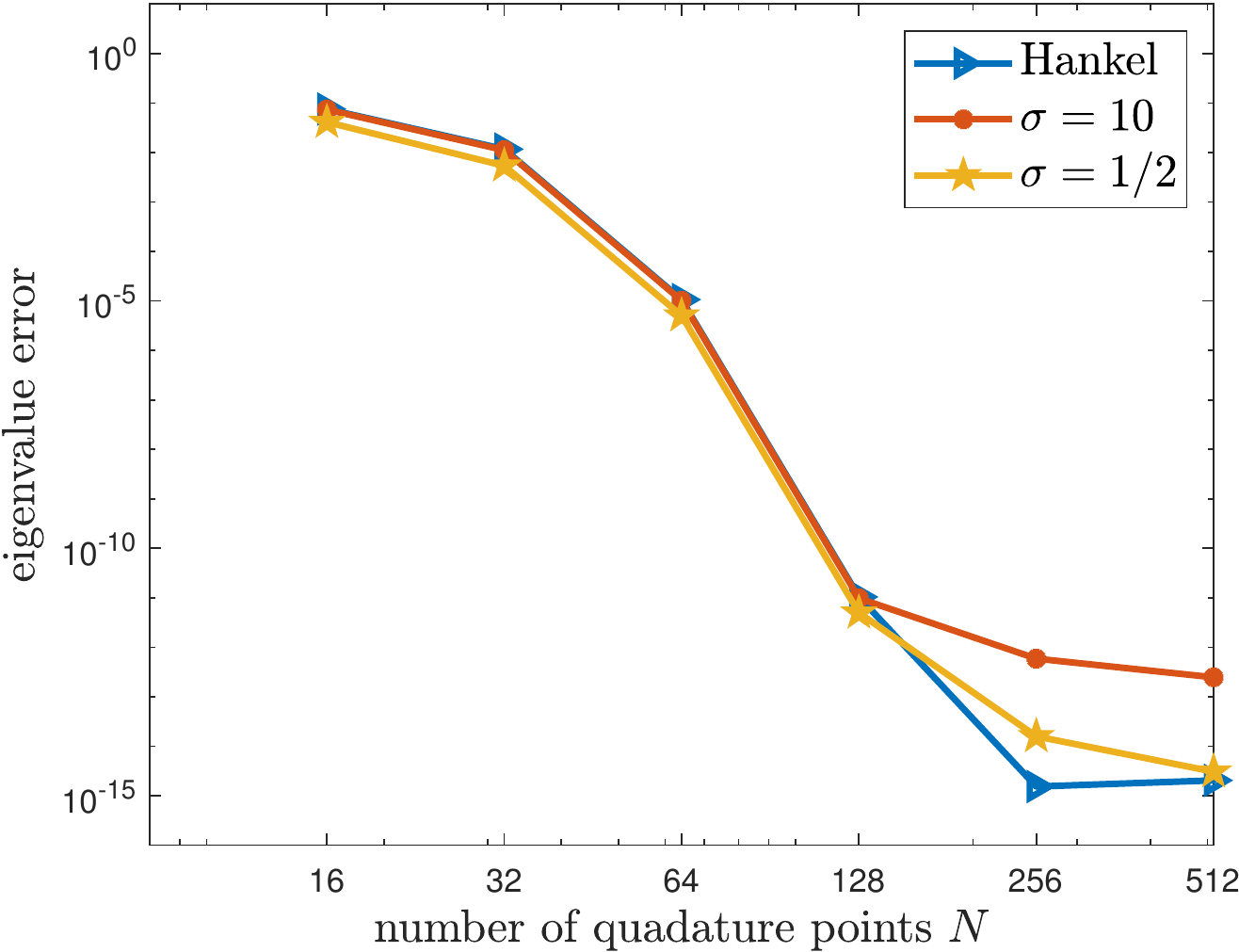}}$
		$\overset{K = 3, \ r = 11}{\includegraphics[width=0.32\textwidth]{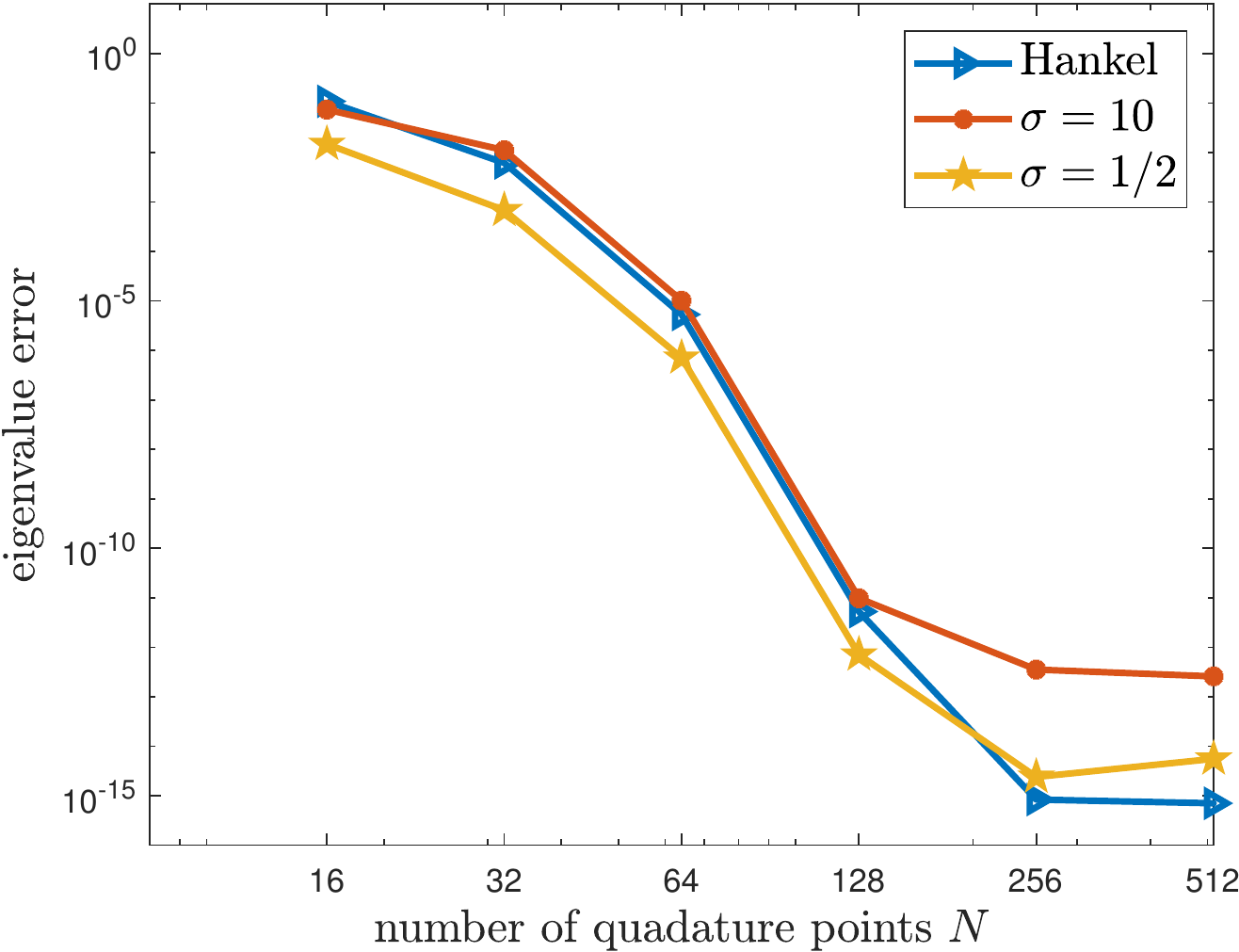}}$
		$\overset{K = 5, \ r = 11}{\includegraphics[width=0.32\textwidth]{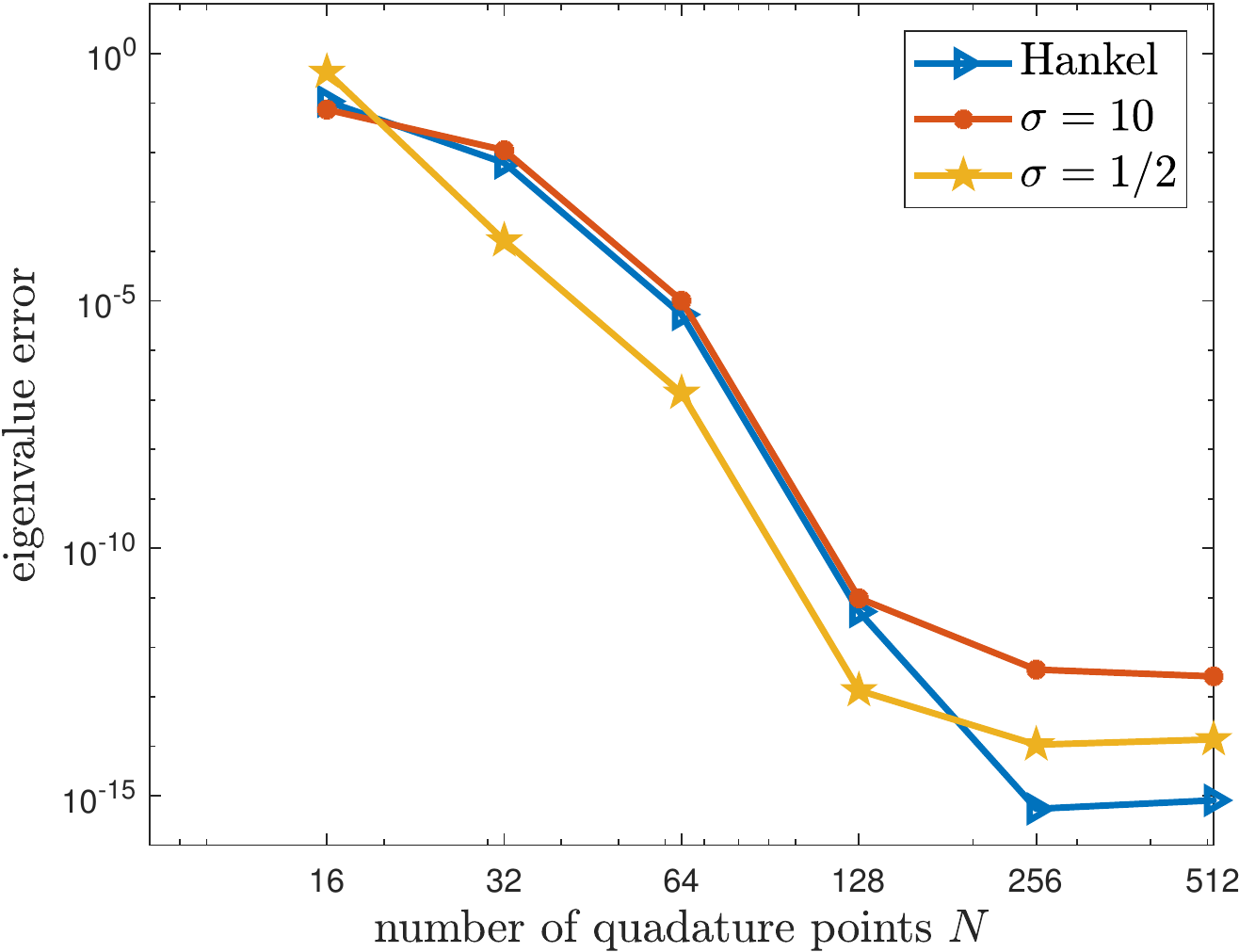}}$
		
		\includegraphics[width=0.32\textwidth]{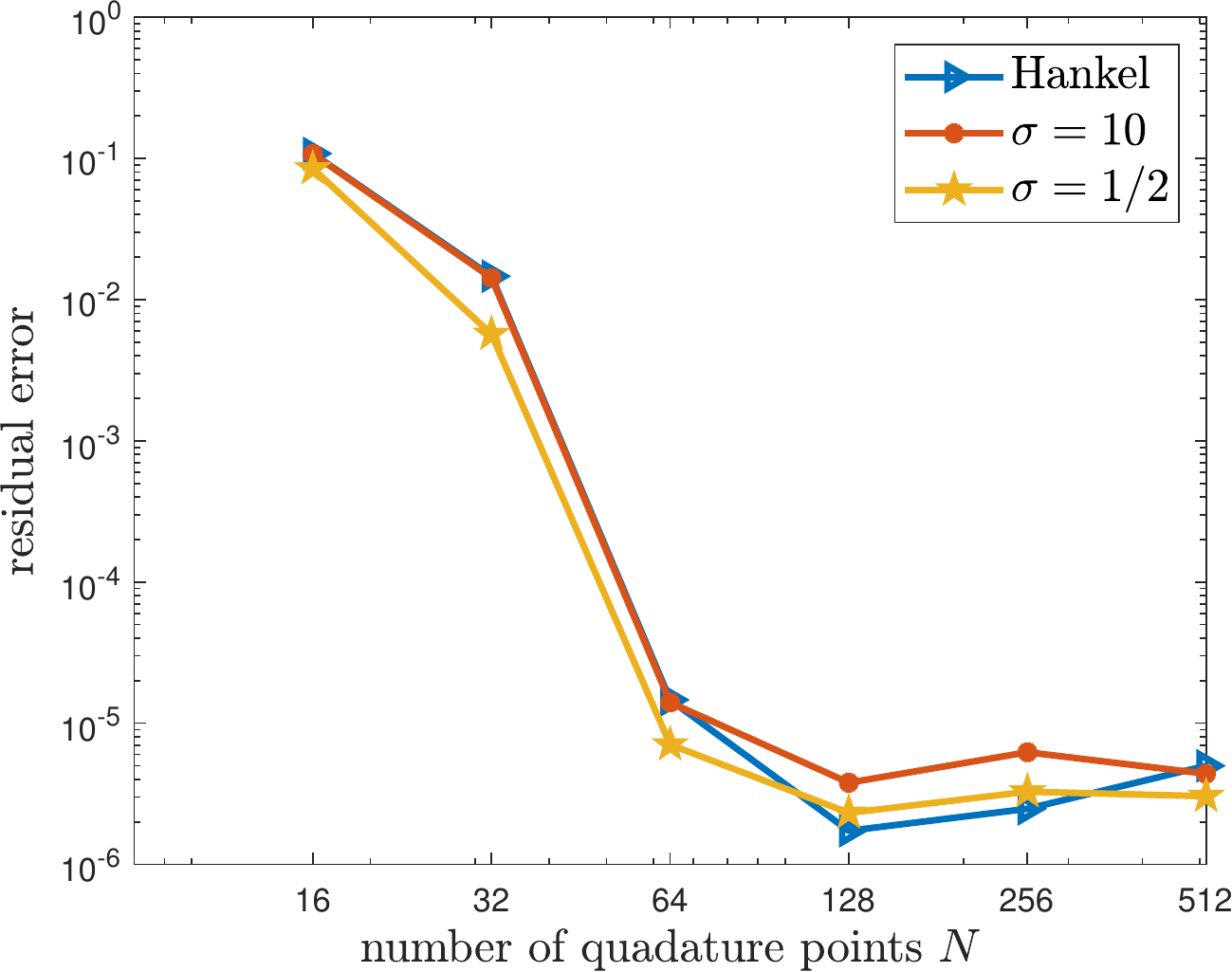}
		\includegraphics[width=0.32\textwidth]{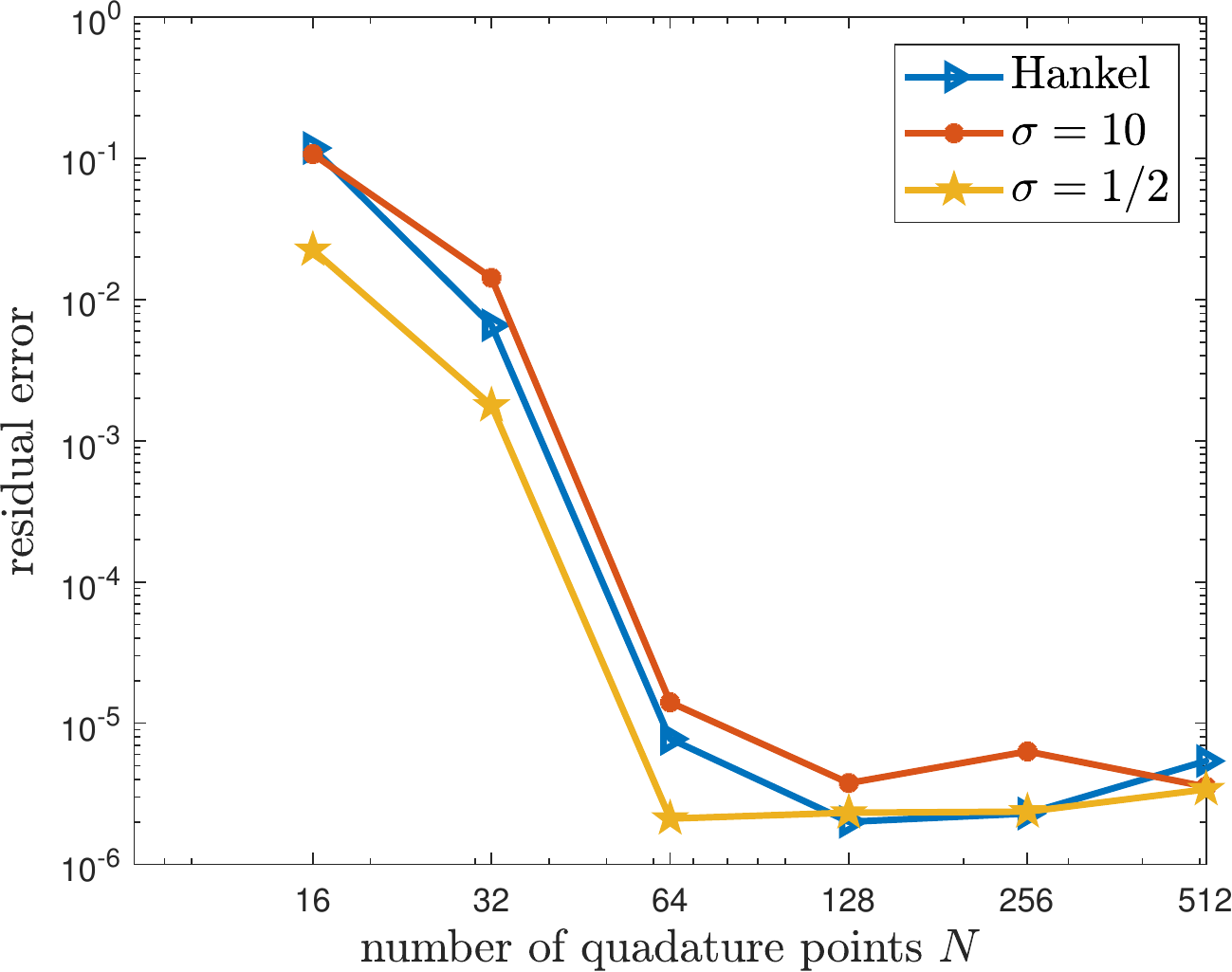}
		\includegraphics[width=0.32\textwidth]{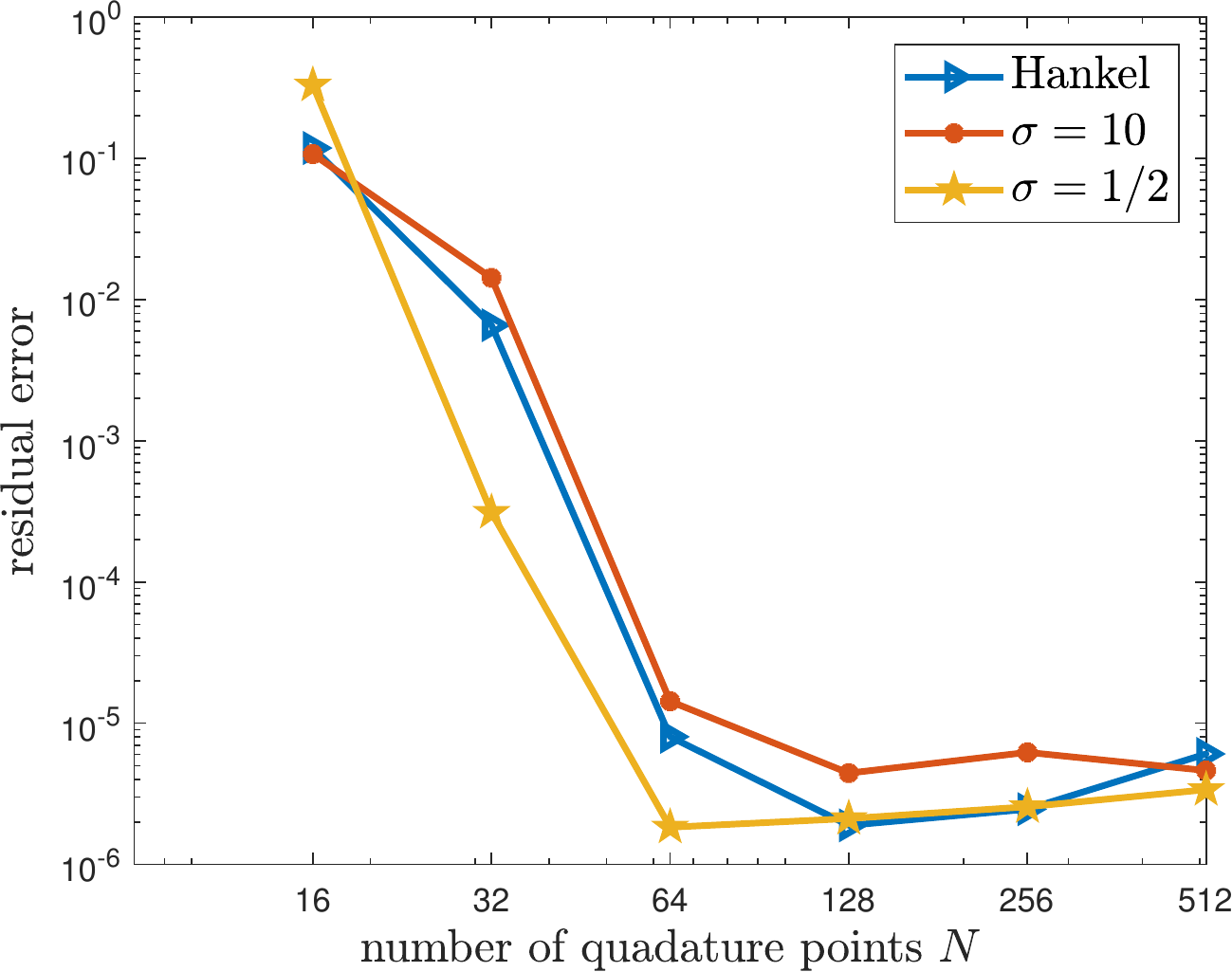}
		
		\includegraphics[width=0.32\textwidth]{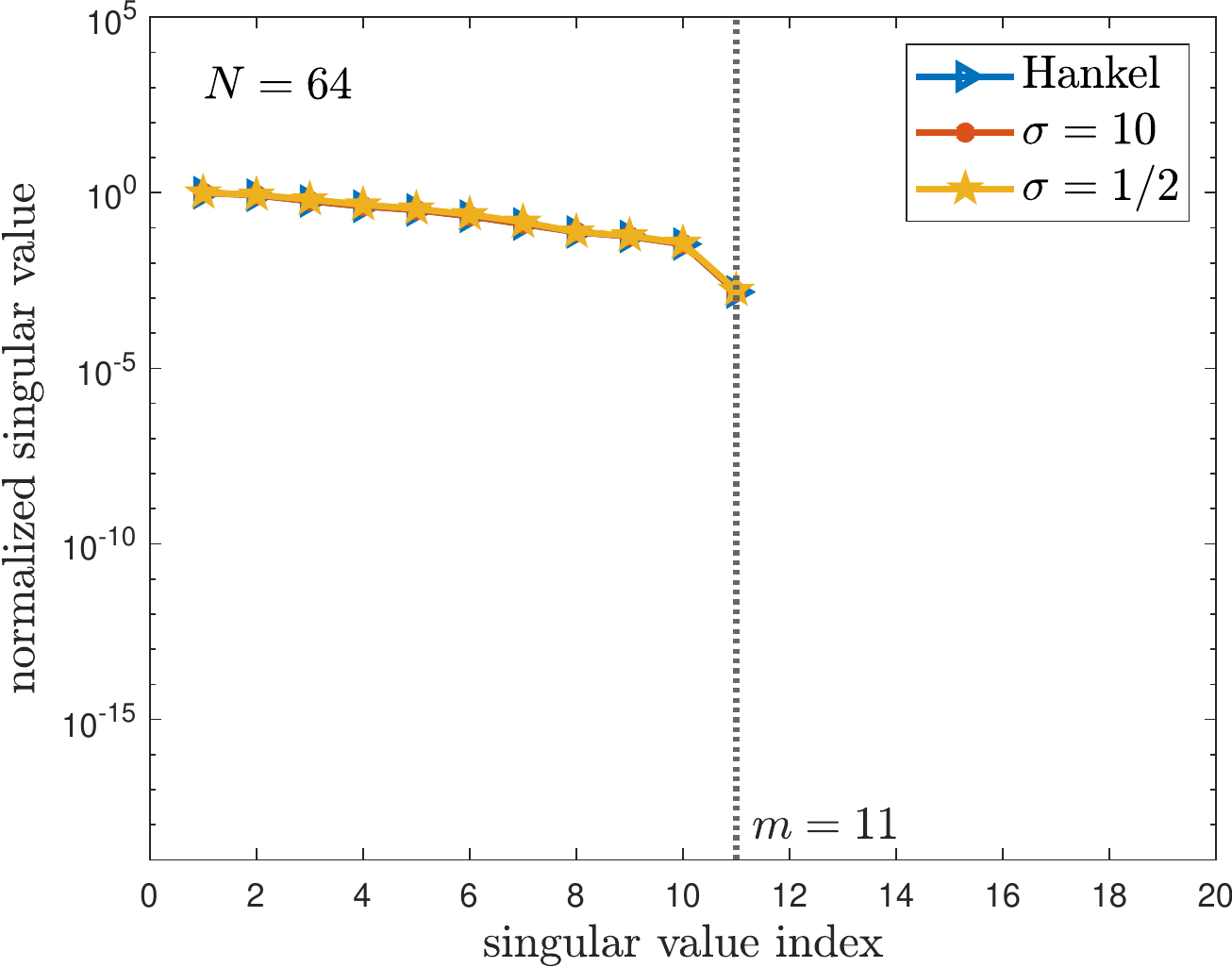}
		\includegraphics[width=0.32\textwidth]{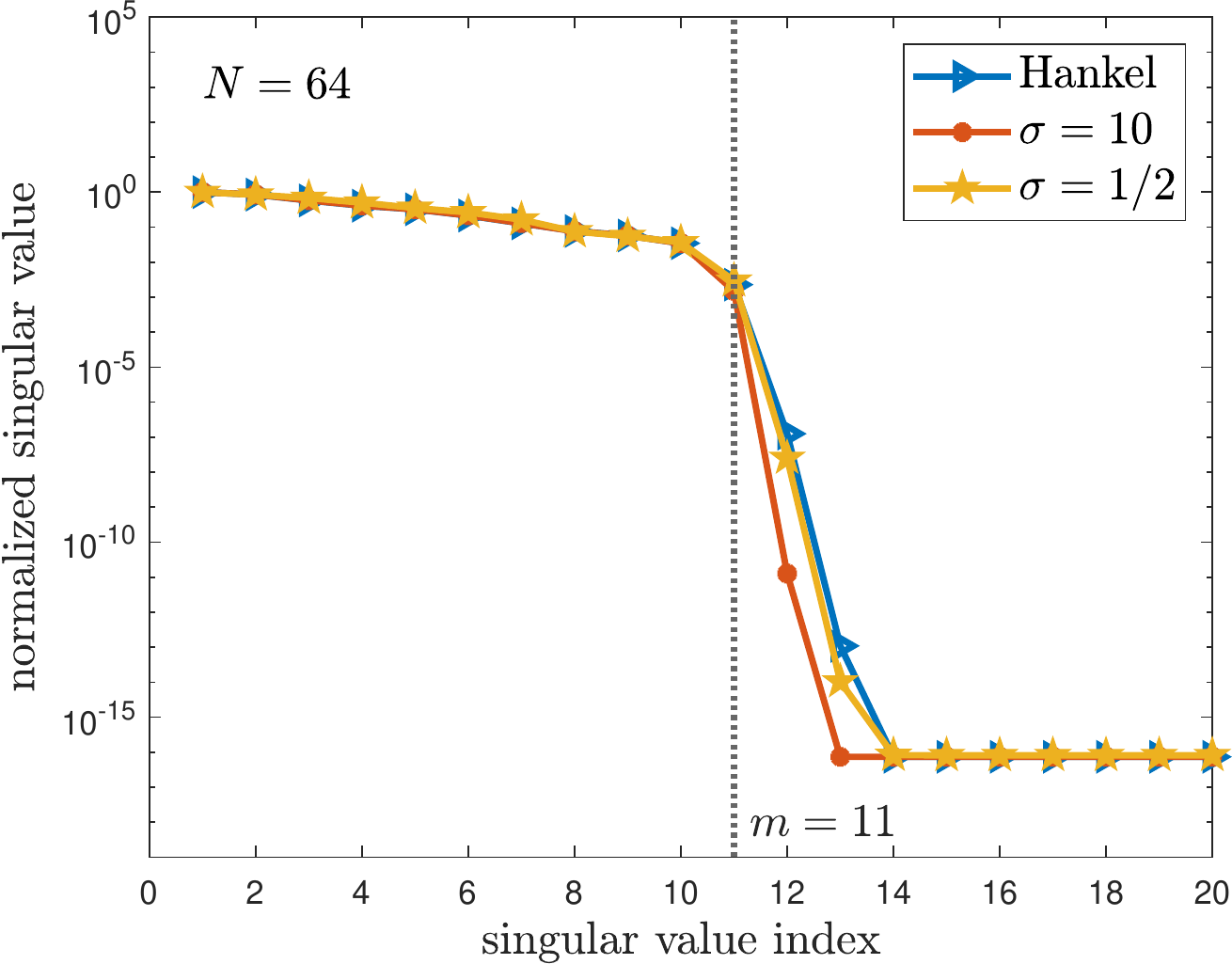}
		\includegraphics[width=0.32\textwidth]{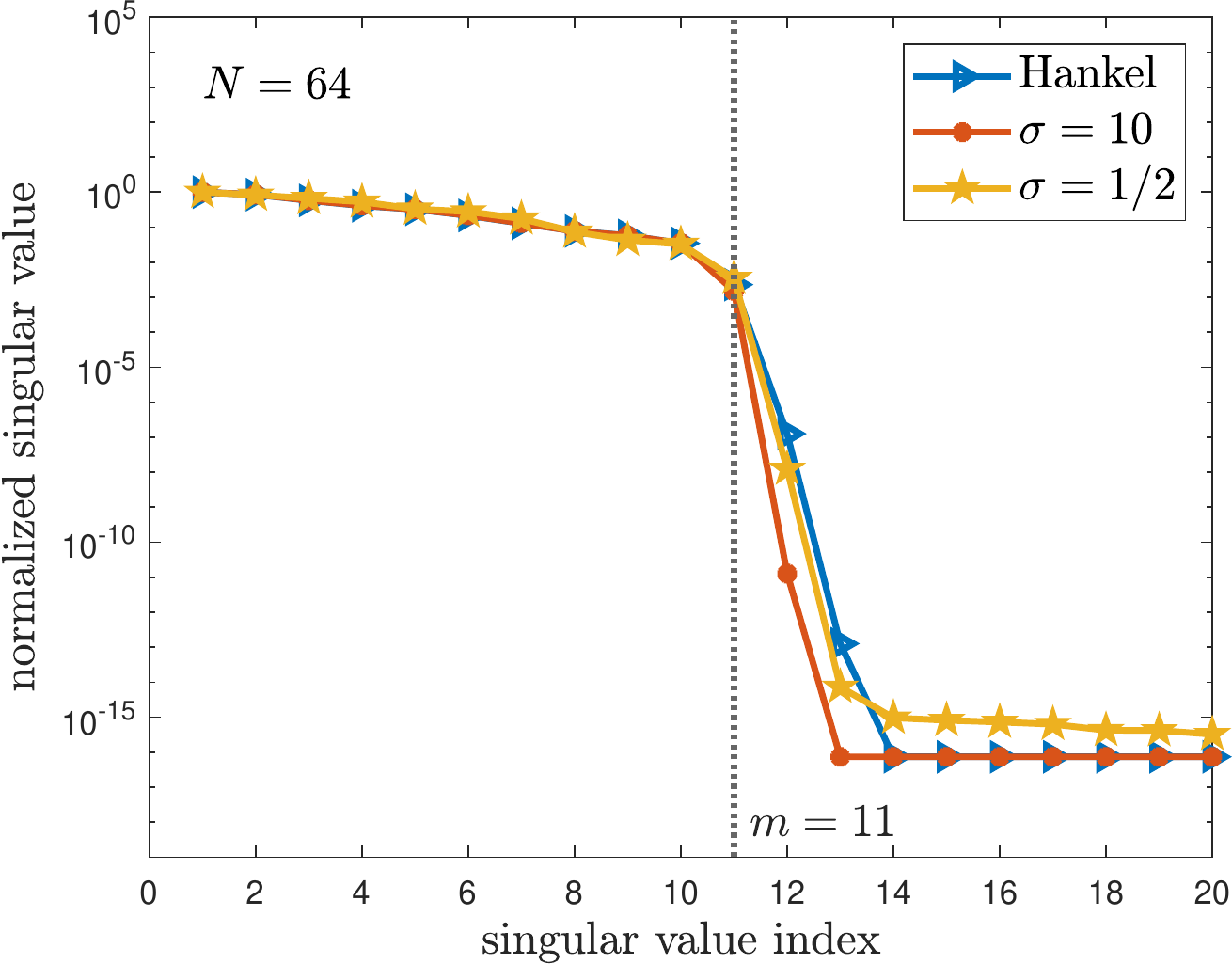}
	\end{center}
	
	\vspace*{-1em}
	\caption{\label{fig:sp_delay}
				Delay problem, $m=11$ eigenvalues in $\Omega$, Hankel and single-point Loewner methods ($\rp=10, 1/2$)
            with $\ell=r=11$ probing directions: 
            the maximum eigenvalue error (top), maximum residual error (middle), 
            and singular values of $\H$ and $\L$ with $N=64$ quadrature points (bottom).
	}
\end{figure}
%%%%%%%%%%%%%%%%%%%%%%%%%%%%%%%%%%%%%%%%%%%%%%%%%%%%%%%%%%%%%%%%%%%%%%%%%%%%%%%%

\Cref{fig:sp_delay} shows the eigenvalue and residual errors for both methods versus the number of quadrature points, $N$. The choice of interpolation point, $\sigma$, affects the accuracy of the computed eigenvalues. 
Taking $\sigma = 1/2$, the single-point Loewner method often yields results that are 
moderately more accurate in both metrics, 
\cref{err1} and \cref{err2}, compared to the Hankel based method.
While the eigenvalue error is not accessible in practice, it is possible to compute the residual error \cref{err2}. 
Therefore, at a negligible cost (compared to computing the quadrature data) one can perform the Hankel method and the single-point Loewner method with multiple choices of $\sigma$, then choose the result that yields the lowest residual error.

\begin{figure}[h!]
	\begin{center}
		$\overset{K = 1, \ r = 32}{\includegraphics[width=0.32\textwidth]{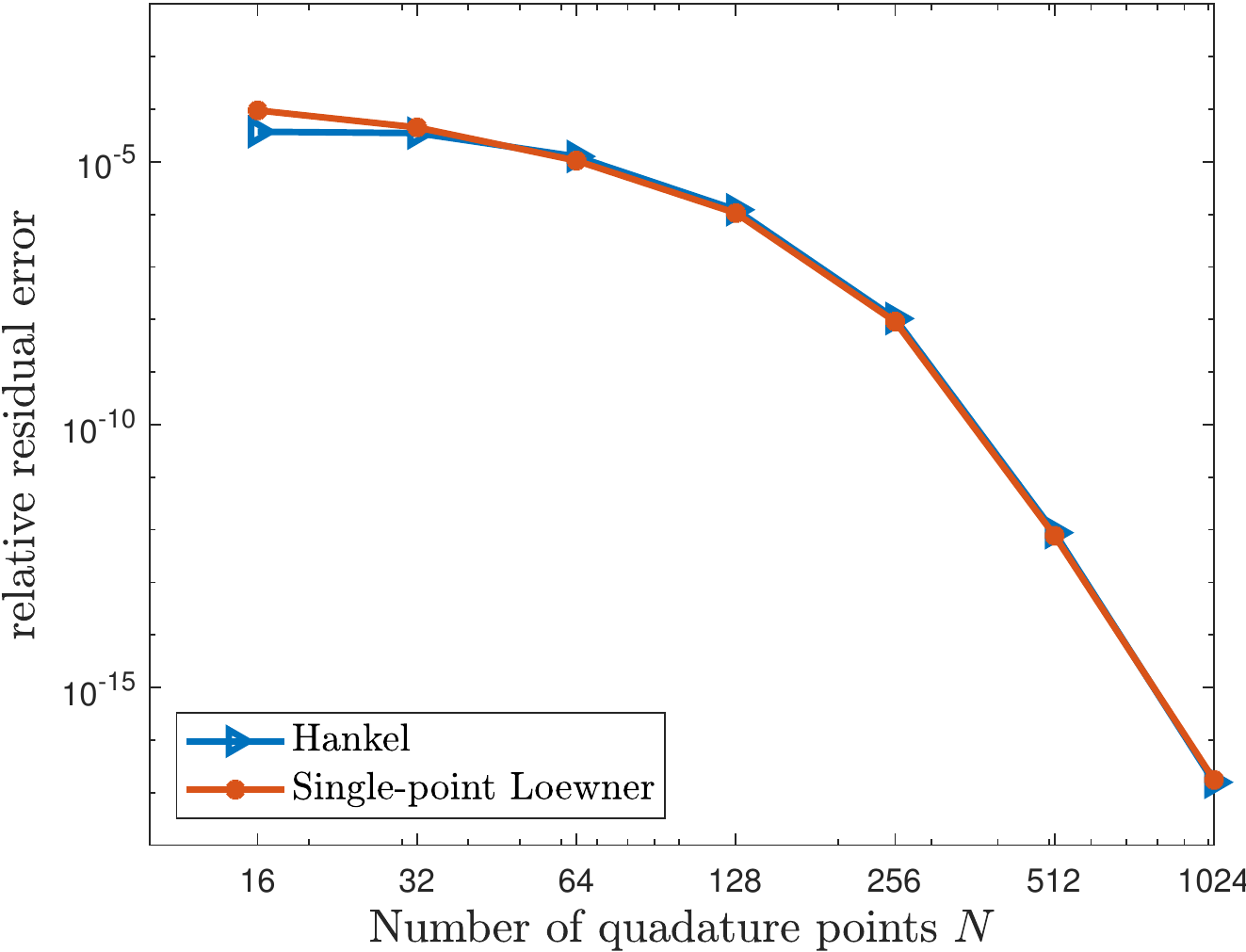}}$
		$\overset{K = 2, \ r = 16}{\includegraphics[width=0.32\textwidth]{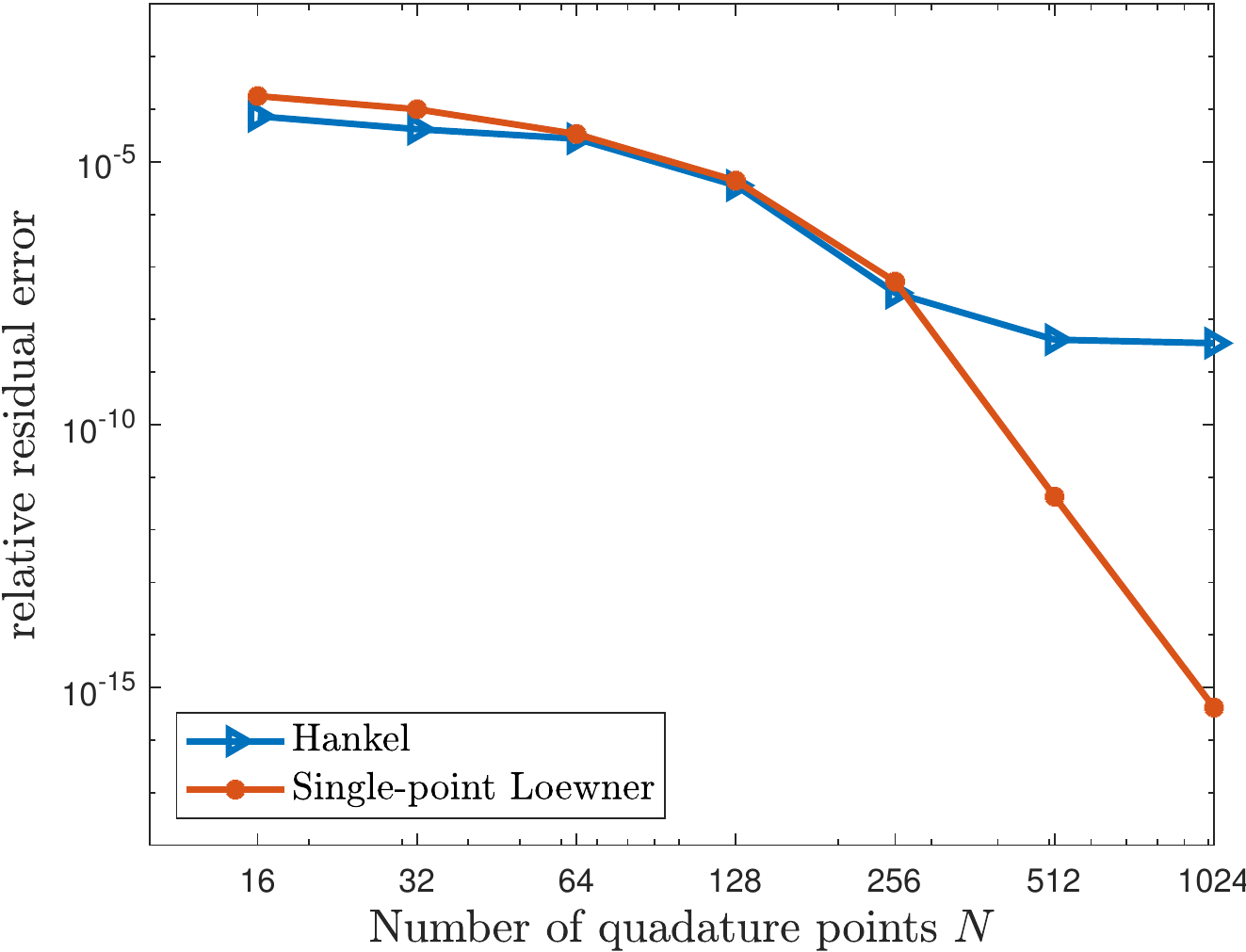}}$
		$\overset{K = 4, \ r = 8}{\includegraphics[width=0.32\textwidth]{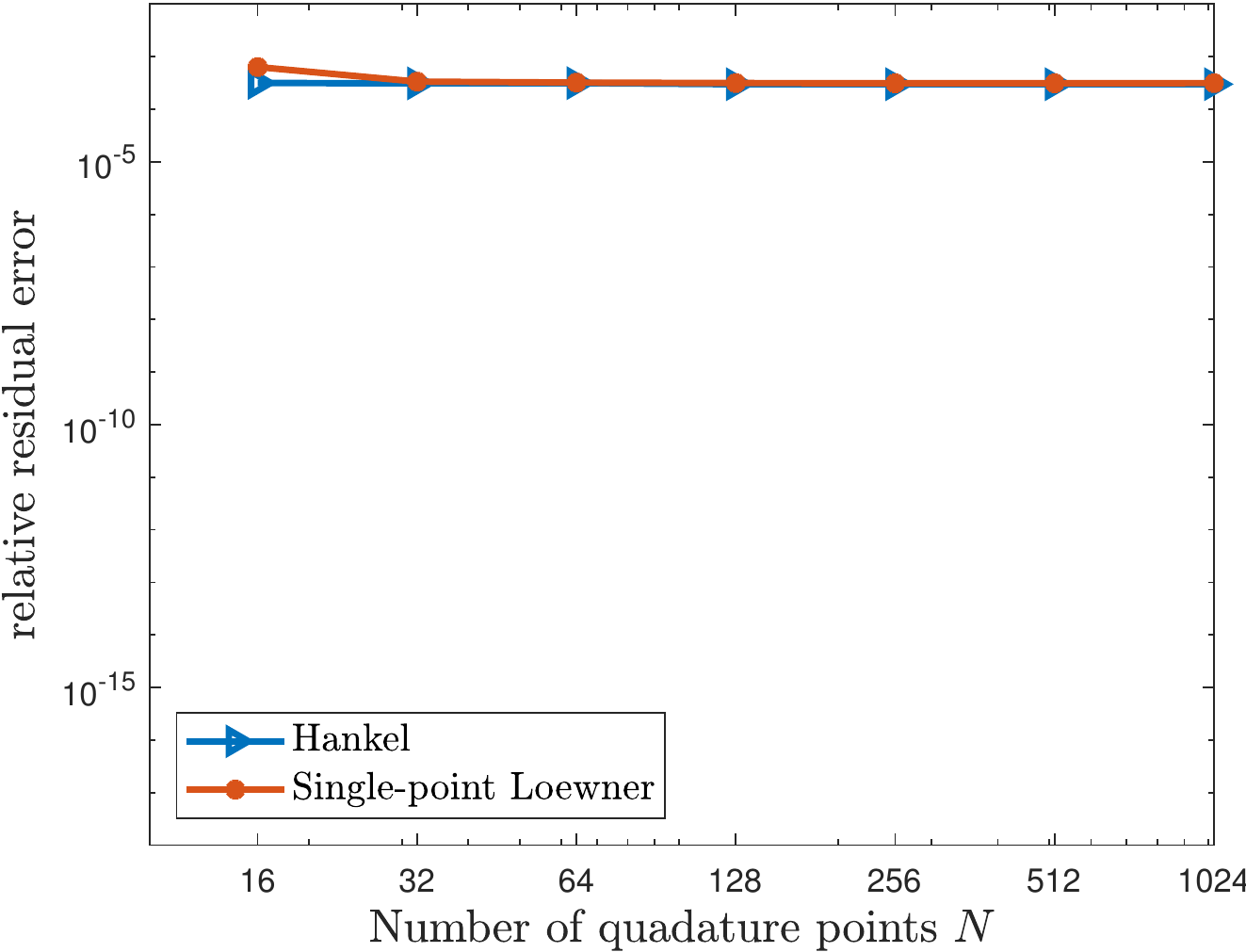}}$
		
		\includegraphics[width=0.32\textwidth]{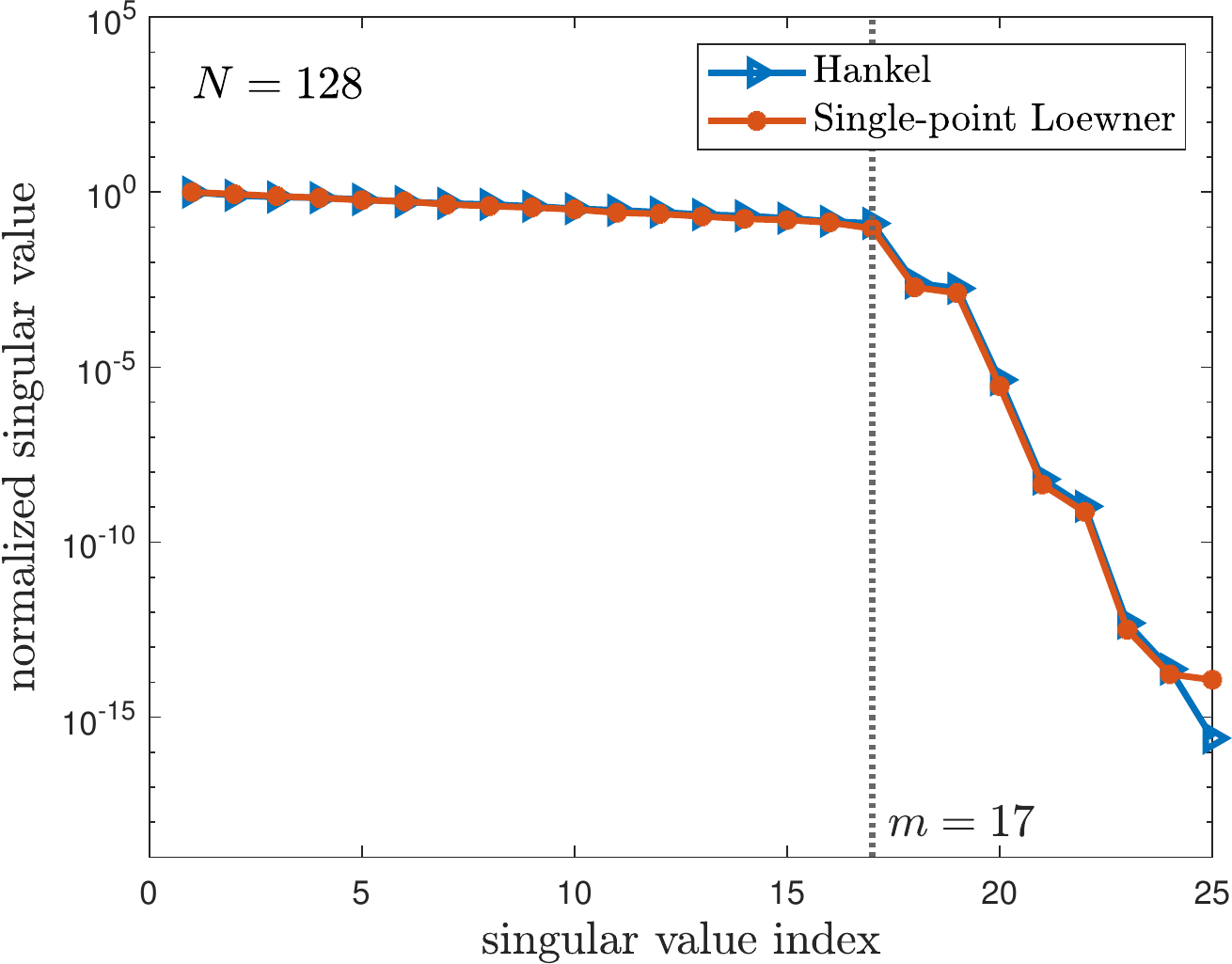}
		\includegraphics[width=0.32\textwidth]{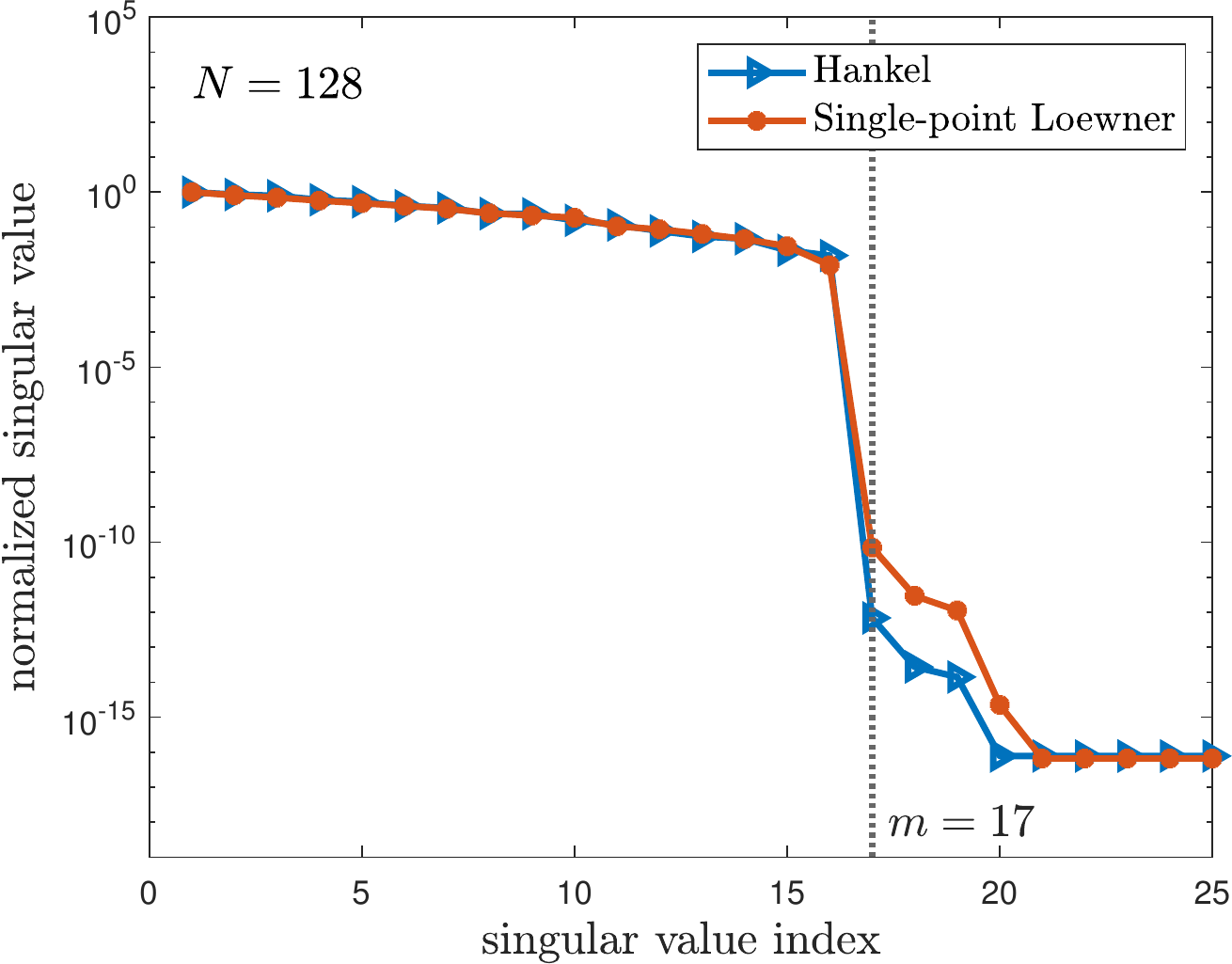}
		\includegraphics[width=0.32\textwidth]{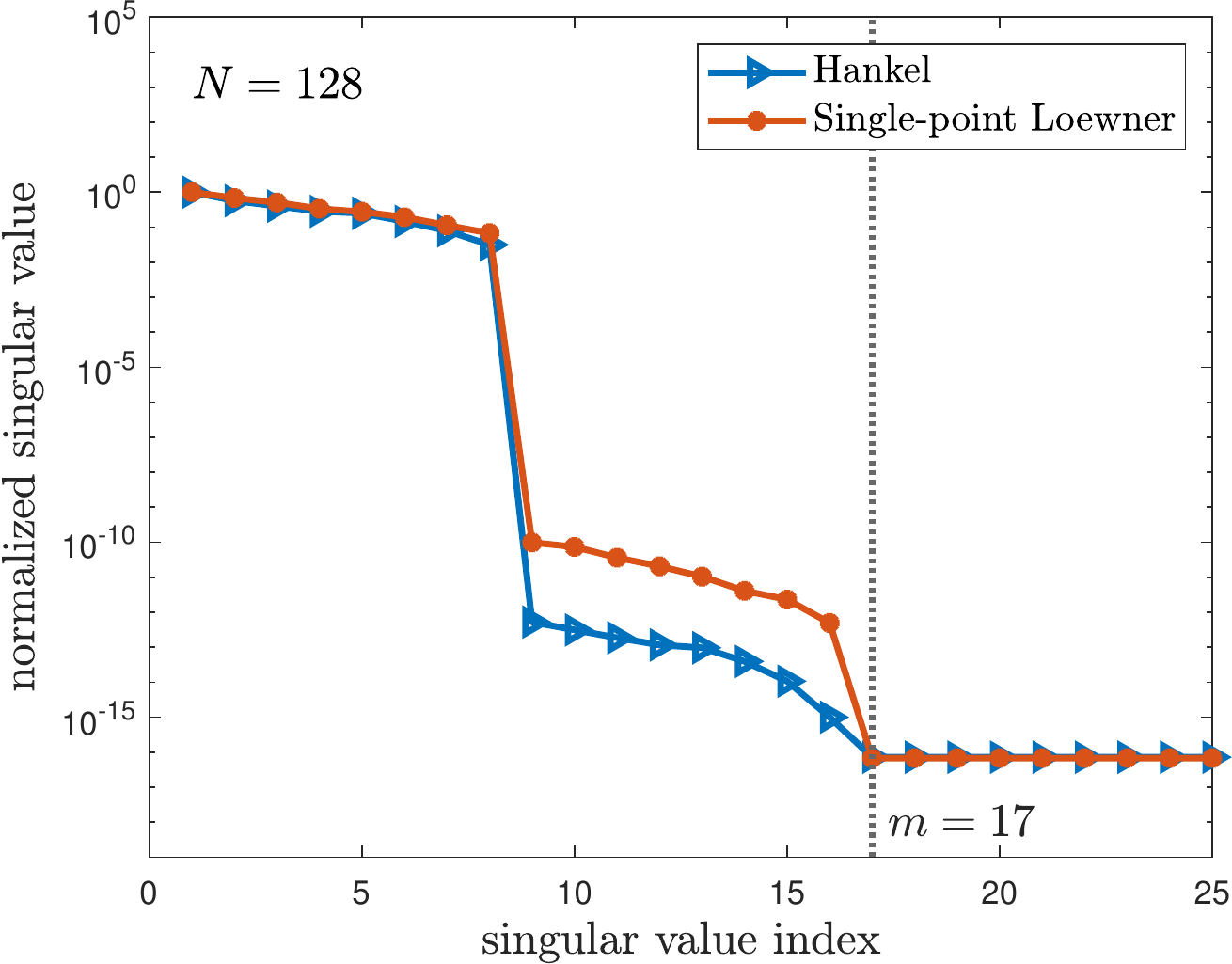}
		
		\includegraphics[width=0.32\textwidth]{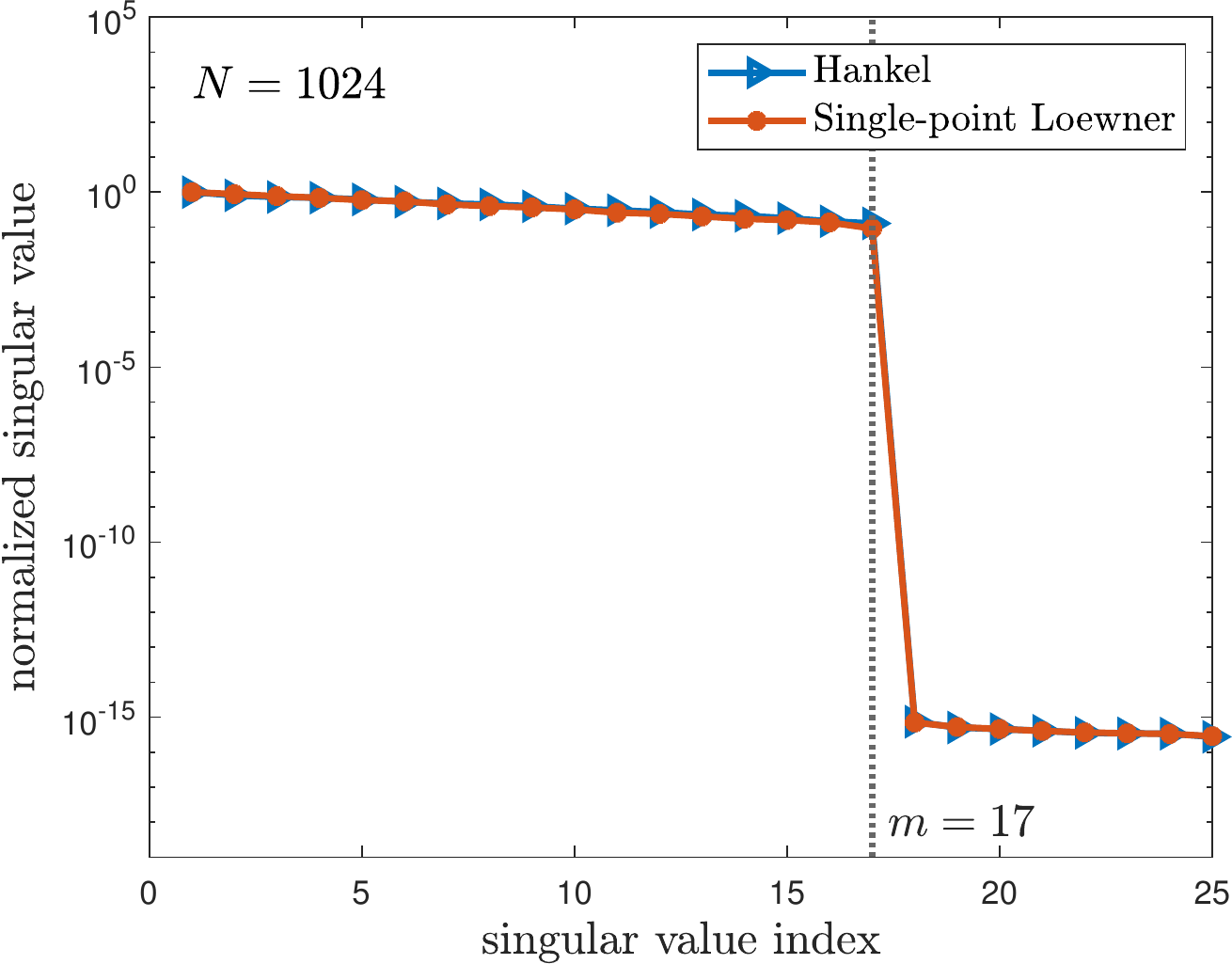}
		\includegraphics[width=0.32\textwidth]{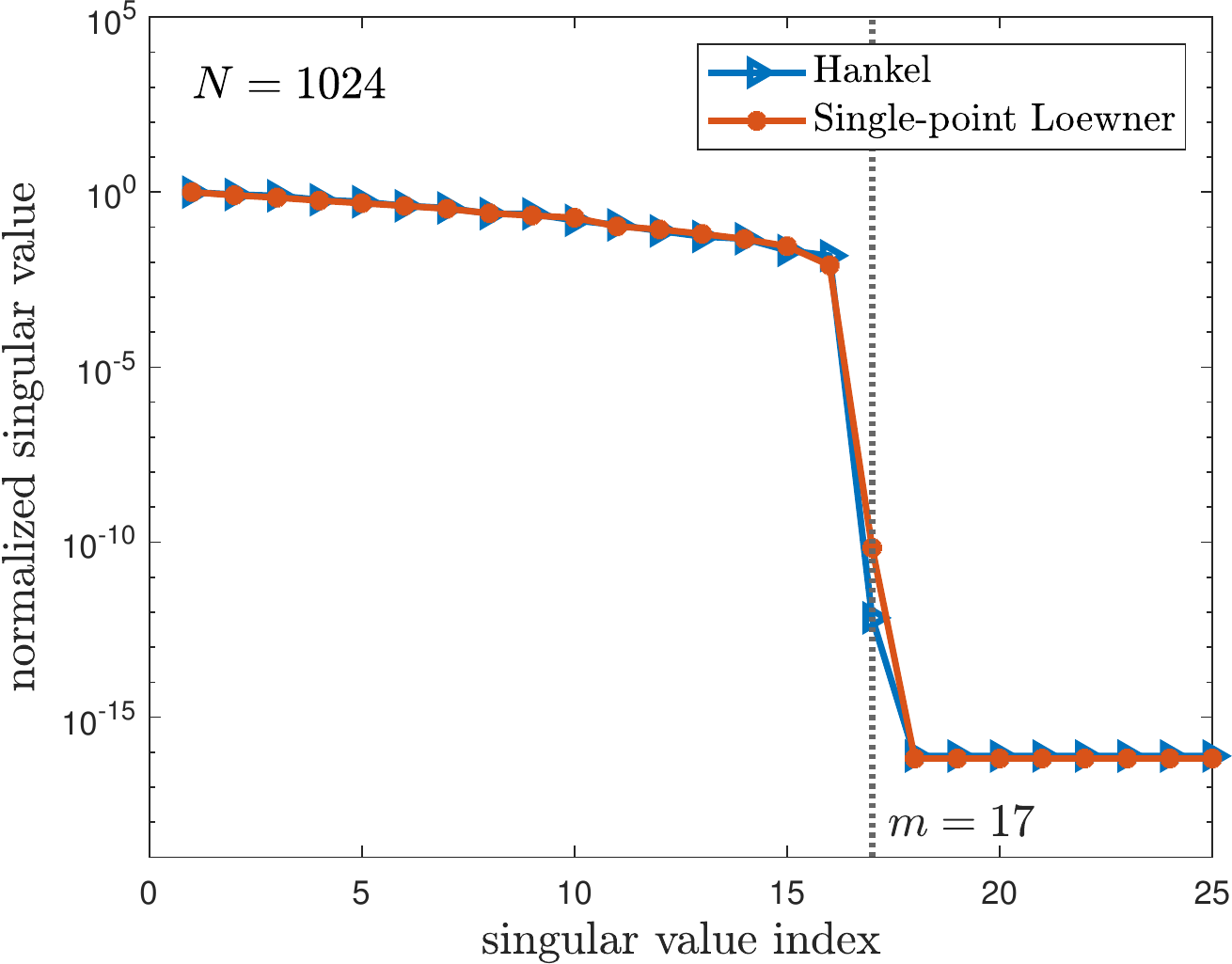}
		\includegraphics[width=0.32\textwidth]{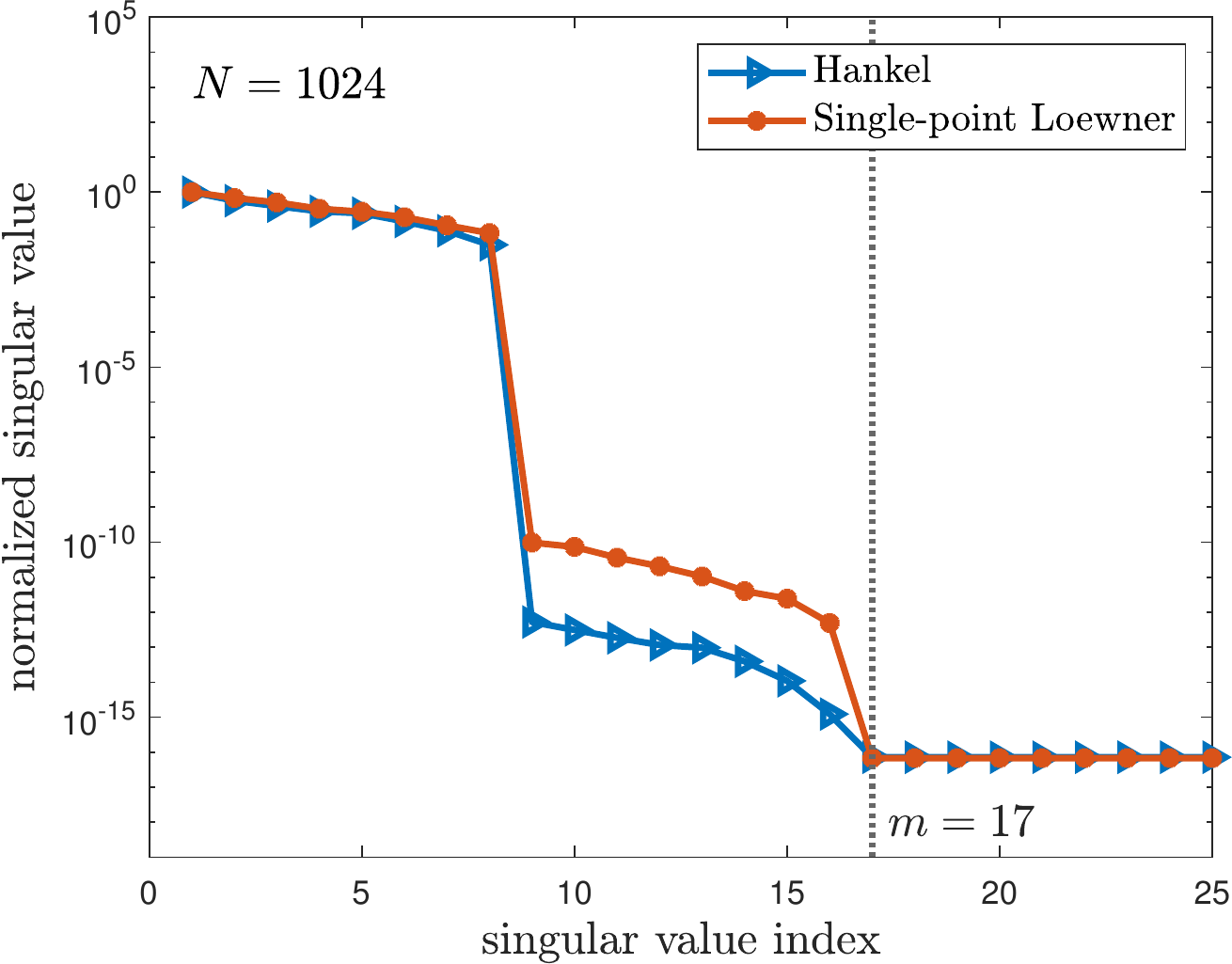}
	\end{center}
	
	\vspace*{-1em}
	\caption{\label{fig:sp_gun}
		{\tt gun} problem, $m=17$ eigenvalues in $\Omega$, Hankel and single-point Loewner methods with three pairs of $K$ and $\ell=r$, 
all giving $K r = 32$: the maximum relative residual error (top); the singular values of $\H$ and $\L$ with $N=128$ (middle) and $N=1024$ (bottom) quadrature points.}
\end{figure}

We next consider a common nonlinear eigenvalue benchmark, the {\tt gun} problem provided by the NLEVP collection \cite{BHMST13}, used to model a radio-frequency gun cavity. The problem takes the form
\[
\BT(z) = \BK - z \BM + \iop@\sqrt{z - \alpha_1^2}\,\BE_1 + \iop@\sqrt{z - \alpha_2^2}\,\BE_2.
\]  
The matrices $\BM, \BK, \BE_1,\BE_2 \in \R^{9956 \times 9956}$ are symmetric, and we use the common parameters $\alpha_1 = 0$ and $\alpha_2 = 108.8774$. We note several key differences in our numerical set-up compared to the results given for this example in~\cite{BK16,BP20}. We do not need one-sided probing, i.e., $\BL = \BI$ or $\BR=\BI$, as typically done 
to recover the eigenvectors.   \Cref{thm:singleshift} allows us to recover approximate eigenvectors even when we probe from both directions. 
Thus our final Hankel and Loewner matrices are $Kr\times Kr$ dimensional ($32 \times 32$ here), versus $n \times Kr$ ($9956 \times 32$) when $\BL = \BI$ or $\BR=\BI$. Our contour is a circle of radius $30000$ centered at $z = 141000$, whereas~\cite{BK16,BP20} center the contour at $z = 140000$. 
Our choice brings the contour closer to an external eigenvalue near the right edge of the contour, and farther from external eigenvalues left of the contour. 
We measure the accuracy of the methods using the relative residual error 
\[
\max_{j=1,\ldots,m}\|\BT(\widetilde{\lambda}_j)\wt{\Bv}_j\|_2/\|\BT(\widetilde{\lambda}_j)\|_F.
\]
The interpolation point was chosen directly under the contour at $\sigma = 141000 - 40000@\iop$. 
\Cref{fig:sp_gun} shows computational results.
First, we take $K = 1$ and set the probing dimensions to $\ell=r = 32$. We see similar results for each method in this case. 
Next we take $K = 2$ with smaller probing dimensions, $\ell = r = 16$. The single-point Loewner method remains accurate in this case, while the maximum residual error for the Hankel method apparently fails to converge as the number of quadrature points increases. Decreasing the probing dimensions further to $\ell = r = 8$ and computing four moments ($K = 4$), both methods fail to converge. 
In each case, the Hankel and Loewner matrices have size $32 \times 32$, but in the second and third cases the probing dimension $r$ is smaller than the number of eigenvalues in $\Omega$, $m = 17$. 
Both these methods should recover $m=17$ accurate eigenpairs, given exact data and computations. 
The middle and bottom rows of \Cref{fig:sp_gun} show the singular values of the Hankel and Loewner matrices for each of the three cases discussed above, with $N=128$ and $N=1024$ quadrature points.
These singular values suggest an explanation for the poor accuracy. 
The number of probing directions has a significant influence on the \emph{numerical} rank of $\H$ and $\L$,
indicating numerical error in computing higher moments. 
For $\ell=r=16$ and $\ell=r=8$, the requirement that the Hankel and Loewner matrices have rank equal 
to the number $m$ of eigenvalues in $\Omega$ is violated (as far as these finite precision computations 
can reveal, even for $N=1024$, an unusually large number of quadrature points). 
These results motivate an alternative method that will trade 
higher order interpolation at one point for first order interpolation at several points. 
We will see that the \emph{multi-point Loewner} method resolves this accuracy issue
for the {\tt gun} problem with probing dimensions as low as $\ell=r = 4$.

%%%%%%%%%%%%%%%%%%%%%%%%%%%%%%%%%%%%%%%%%%%%%%%%%%%%%%%%%%%%%%%%%%%%%%%%%%%%%%%%
%%%%%%%%%%%%%%%%%%%%%%%%%%%%%%%%%%%%%%%%%%%%%%%%%%%%%%%%%%%%%%%%%%%%%%%%%%%%%%%%

%% file: loewner.tex
%!TEX root = paper.tex

%%%%%%%%%%%%%%%%%%%%%%%%%%%%%%%%%%%%%%%%%%%%%%%%%%%%%%%%%%%%%%%%%%%%%%%%%%%%%%%%
\section{Multi-point Loewner algorithm for NLEVPs}    \label{sec:multipoint}
%%%%%%%%%%%%%%%%%%%%%%%%%%%%%%%%%%%%%%%%%%%%%%%%%%%%%%%%%%%%%%%%%%%%%%%%%%%%%%%%
In the last two sections we showed how to recover (realize) a rational transfer function
from high-order samples at a single point, using either $z=\infty$ (with $f(z)=z^k$)
or $z=\rp\in\C\setminus \overline{\Omega}$ (with $f(z) = (-1)^k/(\rp-z)^{(k+1)}$).  
{When the maximum order $K$ of these samples gets large
(as necessary to compute linearly dependent eigenvectors),
numerical difficulties can emerge.}
We now adapt these ideas to solve the NLEVP using 
\emph{low-order samples} from a variety of points in $\C\setminus \overline{\Omega}$,
paralleling developments in rational interpolation methods for 
data-driven reduced-order modeling~\cite{AntBG20,MA07}. 
(This multi-point approach can find linearly dependent eigenvectors 
without requiring use of larger powers.)
 
Suppose we have  $2@r$ points%
\footnote{The algorithm extends readily to the case of $\ell$ left and $r$ right interpolation points,
with $\ell\ne r$.  We take $\ell=r$ here to simplify the notation.}
 in $\C\setminus\overline{\Omega}$ 
at which we want to sample $\BH(z)$.
(These points are called \emph{interpolation points},
terminology we will justify in the next subsection.)
We group these points into two sets of $r$~points, 
called the \emph{left interpolation points} $\{\lp_1, \ldots, \lp_r\} \subset {\C\setminus \overline{\Omega}}$ 
and the \emph{right interpolation points} $\{\rp_{1}, \ldots, \rp_{r} \}\subset {\C\setminus \overline{\Omega}}$.
We assume these points are not poles of $\BH(z)$,
as ensured by taking them outside $\overline{\Omega}$.

For each point, we assign a nonzero probing (direction) vector:  Let 
\[
\Bell_1, \ldots, \Bell_r \in \C^n
\] 
denote the (left) directions associated with the (left) points $\lp_1, \ldots, \lp_r$, 
and let 
\[
 \Br_{1}, \ldots,\Br_{r} \in \C^n
\]
denote the (right) directions associated with the (right) points $\rp_1, \ldots, \rp_r$.
Now, assume we have one-sided measurements (interpolation data) of $\BH(z)$ at these interpolation points,
along the selected directions, i.e., we have the left and right data
\begin{equation} \label{eqn:intdata} 
\ld_i^* := \Bell_i^* \BH(\lp_i^{}) \in \C^{1\times n},\qquad  \rd_{j} := \BH(\rp_{j})\Br_{j} \in \C^{n\times 1}, \rlap{$\qquad i, j = 1,\dots r$.}
\end{equation}
As described in \Cref{thm:rationalfunctionCI}, these samples of $\BH(z)$ can be obtained
via contour integrals of $\BT(z)$, without directly accessing the unknown $\BH(z)$:
\begin{align} 
		\ld_i^* := \Bell_i^\ast  \BH(\lp_i) &= \frac{1}{2\pi@\iop} \int_{\partial \Omega} \frac{1}{\lp_i - z}\, \Bell_i^\ast  \BT(z)^{-1}\, \dop z,\label{eqn:leftdata}  \\[.25em]
		\rd_j :=  \BH(\rp_j) \Br_j&= \frac{1}{2\pi@\iop} \int_{\partial \Omega} \frac{1}{\rp_j - z}\, \BT(z)^{-1} \Br_j\, \dop z. \label{eqn:rightdata} 
		\end{align}
We can now pose the realization problem we seek to solve.
%%%%%%%%%%%%%%%%%%%%%%%%%%%%%%%%%%%%%%%%%%%%%%%%%%%%%%%%%%%%%%%%%%%%%%%%%%%%%%%%
\medskip
\begin{center}
\framebox{\smallskip\ \begin{minipage}{4.8in}
\vspace*{2pt}
\textbf{Realization problem: Data at multiple points, one-sided samples}\\
\textsl{Given left samples $\{\ld_i^* = \Bell_i^*\BH(\lp_i)\}_{i=1}^r$ 
and right samples $\{\rd_j = \BH(\rp_j)\Br_j\}_{j=1}^r$,
along the directions $\{\Bell_i\}_{i=1}^r$ and $\{\Br_j\}_{j=1}^r$,
construct the transfer function $\BH(z)$ for the full system~\cref{eqn:LTI}.}

\vspace*{2pt}
\end{minipage}\ }
\end{center}
\medskip
%%%%%%%%%%%%%%%%%%%%%%%%%%%%%%%%%%%%%%%%%%%%%%%%%%%%%%%%%%%%%%%%%%%%%%%%%%%%%%%%

\noindent
When the integrals~\cref{eqn:rightdata} are approximated via quadrature, 
we face a \emph{realization problem with inexact measurements}, a topic of 
interest in systems theory; in the Loewner setting, 
see, e.g., \cite{BGW12,DP,EI}.  The following derivation assumes exact measurements.

\smallskip
\begin{plainremark}
The choice of the interpolation points and directions is quite flexible,
beyond the requirement that the points fall outside $\overline{\Omega}$.
We assume at first that $\lp_i\ne \rp_j$ for all $i,j=1,\ldots,r$
(the $\lp_i=\rp_j$ case is treated in \cref{sec:derivdata}),
but otherwise permit repeated points $\lp_i = \lp_j$ or $\rp_i = \rp_j$, 
or repeated directions, $\Bell_i=\Bell_j$ or $\Br_i=\Br_j$.
Provided a rank condition stated in \cref{thm:multishift} holds, 
the interpolation points and directions are sufficiently rich to
fully recover $\BH(z)$.
\end{plainremark}

\smallskip
Organize the multi-point interpolation data \cref{eqn:intdata}
into $r\times r$ Loewner ($\L$) and shifted Loewner  ($\L_s$)
matrices defined entrywise by
\begin{align}  \label{eqn:Lmulti}
	[@@\L@]_{i,j} &=\frac{\Bell_i^\ast \left[\BH(\lp_i)  - \BH(\rp_j)\right] \Br_j }{\lp_i - \rp_j} 
   = \frac{\ld_i^*\Br_j - \Bell_i^\ast \rd_j}{\lp_i - \rp_j},
	\end{align}
and
\begin{align}
	[@@\L_s@]_{i,j} &= \frac{\Bell_i^\ast \left[\lp_i\BH(\lp_i)  - \rp_j\BH(\rp_j)\right] \Br_j }{\lp_i - \rp_j} =
	\frac{\lp_i@\ld_i^*\Br_j - \rp_j\Bell_i^\ast \rd_j}{\lp_i - \rp_j}
	, \label{eqn:Lsmulti}
\end{align}
 for $i,j=1,2,\ldots,r$.  Group the measurements into the matrices
 \begin{equation}  \label{eqn:YZmulti}
 \ldm = \left[ \begin{array}{c}
	\ld_1^* \\
	\ld_2^* \\
	\vdots \\
	\ld_r^* 
	\end{array} \right] \in \C^{r  \times n}
	\quad \mbox{and} \quad
\rdm = \left[\,\rd_1 ~~\rd_2 ~~ \cdots ~~ \rd_r \,\right]\in \C^{n \times r}.
 \end{equation}
Notice that $\L$, $\L_s$, $\ldm$, and $\rdm$ only contain
interpolation data.  
Following the template of the last two sections, 
use spectral quantities in $\BH(z) = \BV(z\BI-\BLambda)^{-1}\BW^*$
to define the \emph{generalized observability matrix} $\OO$
and the \emph{generalized reachability matrix} $\RR@$: 
\begin{align} 
\label{eqn:obs}
	\OO &= \left[ \begin{array}{c} \Bell_1^\ast\BV(\lp_1\BI - \BLambda)^{-1} \\
	\vdots\\
	\Bell_r^\ast\BV(\lp_r\BI - \BLambda)^{-1}
	\end{array}   \right] \in \C^{r \times m} \\[.5em]
 \label{eqn:reach}
	 \RR &= \left[(\rp_1\BI - \BLambda)^{-1}\BW^*\Br_1,\ \cdots \ ,(\rp_r\BI - \BLambda)^{-1}\BW^*\Br_r \right]\in \C^{m \times r}.
\end{align}
As in~\cite[sect.~4.1]{AntBG20}, we use the First Resolvent Identity
\[ (\lp_i\BI-\BLambda)^{-1}(\rp_j\BI-\BLambda)^{-1} 
    = \frac{(\lp_i\BI-\BLambda)^{-1}-(\rp_j\BI-\BLambda)^{-1}}{\rp_j-\lp_i} \]
to simplify the $(i,j)$ entry of the product $\OO\RR$:
\begin{align*}
[@\OO\RR@]_{i,j} 
   &= \Bell_i^*\BV(\lp_i\BI-\BLambda)^{-1}(\rp_j\BI-\BLambda)^{-1}\BW^*\Br_j \\[.5em]
   &= \frac{\Bell_i^*\BV(\lp_i\BI-\BLambda)^{-1}\BW^*\Br_j
            - \Bell_i^*\BV(\rp_j\BI-\BLambda)^{-1}\BW^*\Br_j}{\rp_j-\lp_i} \\[.5em]
   &= \frac{\Bell_i^*\BH(\lp_i)\Br_j - \Bell_i^*\BH(\rp_j)\Br_j}{\rp_j-\lp_i} 
    = \frac{\ld_i^*\Br_j - \Bell_i^*\rd_j}{\rp_j-\lp_i} 
    = -[@@\L@]_{i,j},
\end{align*}
and hence $\OO\RR = -\L$.\ \ 
Similarly, the identity 
\[ (\lp_i\BI-\BLambda)^{-1}\BLambda(\rp_j\BI-\BLambda)^{-1} 
    = \frac{\lp_i(\lp_i\BI-\BLambda)^{-1}-\rp_j(\rp_j\BI-\BLambda)^{-1}}{\rp_j-\lp_i} \]
implies that
\[ [@\OO\BLambda\RR@]_{i,j}
   = \frac{\lp_i \Bell_i^*\BV(\lp_i\BI-\BLambda)^{-1}\BW^*\Br_j 
             - \rp_j \Bell_i^*\BV(\rp_i\BI-\BLambda)^{-1}\BW^*\Br_j}{\rp_j-\lp_i} 
   = -[@@\L_s@]_{i,j},\]
and hence $\OO\BLambda\RR = -\L_s$. 
Moreover, $\OO\BW^* = \ldm$ and $\BV\RR=\rdm$.

\medskip
Given this set-up, the recovery of $\BH(z)$ proceeds exactly as expected from the
last two sections.
Assume that $\rank(\L)=m$ and take the reduced SVD
\begin{equation} \label{eqn:tsvd}
 \L = \BX\BSigma\BY^*.
\end{equation}
Since $\BSigma = \BX^*\L\BY = -(\BX^*\OO)(\RR\BY)$ has rank $m$,
so too must the $m\times m$ matrices $\BX^*\OO$ and $\RR\BY$.
Take the full system~\cref{eqn:LTI}, transform coordinates to 
$\Bx(t) = \RR\BY@\wh{\Bx}(t)$, and premultiply the resulting 
state equation by $\BX^*\OO$ to obtain
\begin{align*}
   (\BX^*\OO)(\RR\BY)@\wh{\Bx}@'(t) 
       &= (\BX^*\OO)\BLambda(\RR\BY)@\wh{\Bx}(t) + (\BX^*\OO) \BW^*\Bu(t) \\[.25em]
    \By(t) &= \BV\RR\BY @\wh{\Bx}(t).
\end{align*}
Since $\OO\RR=-\L$ and $\OO\BLambda\RR = -\L_s$, this system gives the realization
\begin{equation} \label{eqn:Hmultipoint}
 \BH(z) = (\rdm\BY)(\BX^*\L_s\BY - z@\BSigma)^{-1}(\BX^*\ldm).  
\end{equation}

%%%%%%%%%%%%%%%%%%%%%%%%%%%%%%%%%%%%%%%%%%%%%%%%%%%%%%%%%%%%%%%%%%%%%%%%%%%%%%%%
\begin{theorem}  	\label{thm:multishift}
Given the points $\lp_1, \ldots, \lp_r \in \C\setminus \overline{\Omega}$ 
and $\rp_{1}, \ldots, \rp_{r} \in \C\setminus \overline{\Omega}$, 
and associated probing vectors $\Bell_1, \ldots, \Bell_r \in \C^n$
and $\Br_{1}, \ldots, \Br_{r} \in \C^n$, 
compute the samples~\cref{eqn:intdata} 
using the contour integrals~\cref{eqn:leftdata} and~\cref{eqn:rightdata}.  
Assume $\lp_i \ne \rp_j$ for all $i,j \in \{1,\ldots,r\}$.
Let $\L$ and $\L_s$ be the Loewner and shifted Loewner matrices
given in \cref{eqn:Lmulti} and \cref{eqn:Lsmulti}, 
and suppose ${\rm rank}(\L) = m$.
Take the reduced SVD $\L = \BX\BSigma\BY^*$.
Arrange the samples into $\ldm$ and $\rdm$, as in \cref{eqn:YZmulti}.
	Define $\BB_{\theta,\sigma}  \in \C^{m \times m}$ and take its eigenvalue decomposition:
	\begin{equation} \label{eqn:Bsigmalti}
	\BB_{\theta,\sigma} :=   \BSigma^{-1} \BX^* \L_s\BY = \BS \BLambda \BS^{-1}.
	\end{equation}
The matrix $\BLambda = {\rm diag}(\lambda_1,\ldots, \lambda_m)$ reveals the 
$m$ eigenvalues of $\BT(z)$ in $\Omega$.\ Let $\Bs_j$ denote
the $j$th column of $\BS$.  Then $(\lambda_j,  \bC \BY \Bs_j)$ is an eigenpair of $\BT(z)$.
\end{theorem}
%%%%%%%%%%%%%%%%%%%%%%%%%%%%%%%%%%%%%%%%%%%%%%%%%%%%%%%%%%%%%%%%%%%%%%%%%%%%%%%%

%%%%%%%%%%%%%%%%%%%%%%%%%%%%%%%%%%%%%%%%%%%%%%%%%%%%%%%%%%%%%%%%%%%%%%%%%%%%%%%%
\begin{algorithm}[b!]
 \caption{Compute quadrature data for single-sided contour integral algorithms}
 \label{alg:quad_data}
 \begin{enumerate}
 \item[ ] Input:  $\{\qn_k\}_{k=1}^N$ and $\{\qw_k\}_{k=1}^N$, 
               quadrature nodes and weights 
 \item[ ] Input:  $r\ge m$, the number of left and right interpolation points.
 \item[ ] Input: left and right probing directions, 
           $\Bell_1, \ldots, \Bell_r \in \C^n$,
           $\Br_{1}, \ldots, \Br_{r} \in \C^n$.
\smallskip
 \item[ ] Output: Tensor $\QL\in \C^{r\times n\times N}$ with
                   $\QL(i,:,k) = \Bell_i^* \BT(z_k)^{-1}\in\C^n$.
 \item[ ] Output: Tensor $\QR\in \C^{n\times r\times N}$ with
                   $\QR(:,j,k) = \BT(z_k)^{-1}\Br_j \in \C^n$.

\medskip
\item Store the directions in $\BL := [\Bell_1, \ldots,\Bell_r]$, $\BR := [\Br_1,\ldots,\Br_r] \in \C^{n\times r}$.
\smallskip
\item For $k=1,\ldots,N$, \\
\ \hspace*{1.5em} Compute $\QL(:,:,k) := \BL^*\BT(\qn_k)^{-1}$ and $\QR(:,:,k) := \BT(\qn_k)^{-1}\BR$. 
 \end{enumerate}
\end{algorithm}
%%%%%%%%%%%%%%%%%%%%%%%%%%%%%%%%%%%%%%%%%%%%%%%%%%%%%%%%%%%%%%%%%%%%%%%%%%%%%%%%

%%%%%%%%%%%%%%%%%%%%%%%%%%%%%%%%%%%%%%%%%%%%%%%%%%%%%%%%%%%%%%%%%%%%%%%%%%%%%%%%
	\subsection{The multi-point Loewner algorithm}  \label{sec:contmultipoint}
%%%%%%%%%%%%%%%%%%%%%%%%%%%%%%%%%%%%%%%%%%%%%%%%%%%%%%%%%%%%%%%%%%%%%%%%%%%%%%%%
To develop a numerical method for the NLEVP using the multi-point Loewner approach 
we have just described, generate the data~\cref{eqn:intdata}
by approximating the contour integrals~\cref{eqn:leftdata} 
and~\cref{eqn:rightdata} via an $N$-point quadrature rule
with nodes~$\{\qn_k\}_{k=1}^N$ and weights~$\{\qw_k\}_{k=1}^N$:
for $i,j = 1,\ldots, r$,
\begin{align*}
\ld_i^* 
:= \Bell_i^\ast  \BH(\lp_i) 
   &= \frac{1}{2\pi \iop} \int_{\partial \Omega} 
      \frac{1}{\lp_i - z}\, \Bell_i^\ast  \BT(z)^{-1}\, \dop z
\ \approx\  
\sum_{k=1}^N \frac{\qw_k}{\lp_i - \qn_k} \Bell_i^*\BT(\qn_k)^{-1}, \\[.25em]
\rd_j
:= \BH(\rp_j) \Br_j
   &= \frac{1}{2\pi \iop} \int_{\partial \Omega} 
      \frac{1}{\rp_j - z}\, \BT(z)^{-1} \Br_j\, \dop z
\ \approx\  
\sum_{k=1}^N \frac{\qw_k}{\rp_j - \qn_k} \BT(\qn_k)^{-1}\Br_j.
\end{align*}
To expedite these computations, collect the left and right probing directions into
\[ \BL = \left[\begin{array}{ccc} \Bell_1 & \cdots & \Bell_r \end{array}\right]\in\C^{n\times r}, 
\qquad
   \BR = \left[\begin{array}{ccc} \Br_1 & \cdots & \Br_r \end{array}\right]\in\C^{n\times r},
\]
and for each quadrature node $\qn_k$, $k=1,\ldots, N$, solve systems to 
obtain the data for all probing directions at once:
\begin{equation} \label{eq:quaddata}
  \BL^*\BT(\qn_k)^{-1}, \qquad \BT(\qn_k)^{-1}\BR.
\end{equation}
\Cref{alg:quad_data} describes the computation of this quadrature data, 
storing the results in two tensors to facilitate reuse.
\emph{Note that these computations, the bulk of the work required for the multi-point 
Loewner contour integral algorithm, 
are independent of the choice of left and right interpolation points, 
$\{\lp_i\}_{i=1}^r$ and $\{\rp_j\}_{j=1}^r$.}
Quadrature approximations to the samples
$\ld_i^* = \Bell_i^*\BH(\lp_i)$ and $\rd_j = \BH(\rp_j)\Br_j$,
simply amount to different weighted sums of the data computed in~\cref{eq:quaddata}.
By combining the same data differently, one can also construct 
quadrature approximations to probed samples for the Hankel (\Cref{sec:hankelsys})
and single-point Loewner (\Cref{sec:singlepoint}) methods at little additional cost.
Similarly, one can readily experiment with different choices of 
$\{\lp_i\}_{i=1}^r$ and $\{\rp_j\}_{j=1}^r$,  provided the probing directions in
$\BL$ and $\BR$ remain fixed.
\Cref{alg:LNLEVP} describes how the quadrature data from \Cref{alg:quad_data}
can then be used in the multi-point Loewner method to arrive at estimates
for the $m$ eigenpairs of $\BT(z)$ in $\Omega$.

%%%%%%%%%%%%%%%%%%%%%%%%%%%%%%%%%%%%%%%%%%%%%%%%%%%%%%%%%%%%%%%%%%%%%%%%%%%%%%%%
\begin{algorithm}[t!]
 \caption{Multi-point Loewner contour integral algorithm for NLEVPs}
 \label{alg:LNLEVP}
 \begin{enumerate}
 \item[ ] Input:  $\Omega$, a target domain containing $m$ eigenvalues ({simple}).
 \item[ ] Input:  $\{\qn_k\}_{k=1}^N$ and $\{\qw_k\}_{k=1}^N$, 
               quadrature nodes and weights for $\Omega$.
 \item[ ] Input:  $r\ge m$, the number of left and right interpolation points.
 \item[ ] Input: left and right points, 
           $\lp_1, \ldots, \lp_r \in {\C\setminus \Omega}$ 
       and $\rp_1, \ldots, \rp_r \in {\C\setminus \Omega}$.
 \item[ ] Input: left and right probing directions, 
           $\Bell_1, \ldots, \Bell_r \in \C^n$
           and $ \Br_{1}, \ldots, \Br_{r} \in \C^n$,\\
           \phantom{Input: }\textsl{e.g., random vectors or approximate eigenvectors.}
\medskip
 \item[] Output: Approximations to $m$ eigenpairs of $\BT(z)$ in $\Omega$.

\medskip
\item Use~\Cref{alg:quad_data} to compute quadrature data $\QL\in \C^{r\times n\times N}$ , $\QR\in\C^{n\times r\times N}$.

\medskip
\item[ ] \hspace*{-17pt} \textsl{Compute quadrature approximations to the left and right data.}
\smallskip
\item For $i=1,\ldots,r$, \\
\ \hspace*{1.5em} Compute left samples:  $\displaystyle{\ld_i^* = \sum_{k=1}^N \frac{\qw_k}{\lp_i-\qn_k} {\QL(i,:,k)}}$.
\item For $j=1,\ldots,r$, \\
\ \hspace*{1.5em} Compute right samples:  $\displaystyle{\rd_j = \sum_{k=1}^N \frac{\qw_k}{\rp_j-\qn_k} {\QR(:,j,k)}}$.

\medskip
\item[ ] \hspace*{-17pt} \textsl{Form Loewner matrices and compute approximate eigenvalues, eigenvectors.}
\smallskip
\item Construct $\L$, $\L_s$, $\ldm$, and $\rdm$  
      as in \cref{eqn:Lmulti}, \cref{eqn:Lsmulti}, and~\cref{eqn:YZmulti}.\\
      (Use \cref{eq:Hermitedata} and \cref{eq:HermiteLoew} to handle any cases where  $\lp_i=\rp_j$.) 
\smallskip
\item Determine the rank $m$ of $\L$ and compute the reduced SVD $\L = \BX\BSigma\BY^*$.
\smallskip
\item Construct $\BB_{\theta,\sigma} :=  \BSigma^{-1} \BX^* \L_s\BY \in \C^{m\times m}$ and
compute its eigendecomposition $\BB_{\theta,\sigma} = \BS \BLambda \BS^{-1}$, 
where $\BLambda = {\rm diag}(\lambda_1,\ldots,\lambda_m)$ 
and $\BS = [\Bs_1\ \cdots\ \Bs_m]$.
\smallskip
\item For $j=1,\ldots, m$, $(\lambda_j,  \rdm \BY \Bs_j)$ is an approximate eigenpair of $\BT(z)$.
 \end{enumerate}
\end{algorithm}
%%%%%%%%%%%%%%%%%%%%%%%%%%%%%%%%%%%%%%%%%%%%%%%%%%%%%%%%%%%%%%%%%%%%%%%%%%%%%%%%

\smallskip
We briefly note that if one seeks only eigenvalue approximations, 
the multi-point Loewner method can be streamlined to require only two-sided samples 
(as in the initial version of Hankel algorithm in \Cref{sec:hankel}).  
One simply needs the Loewner and shifted Loewner matrices, 
with entries computed via the contour integrals
\begin{align*}
     (\L)_{i,j} &= {1\over 2\pi@\iop} \int_{\partial \Omega} \frac{-1}{(\lp_i-z)(\rp_j-z)}\Bell_i^*\BT(z)^{-1}\Br_j\,\dop z,\\[.5em]
   (\L_s)_{i,j} &= {1\over 2\pi@\iop} \int_{\partial \Omega} \frac{-z}{(\lp_i-z)(\rp_j-z)}\Bell_i^*\BT(z)^{-1}\Br_j\,\dop z.
\end{align*}

%%%%%%%%%%%%%%%%%%%%%%%%%%%%%%%%%%%%%%%%%%%%%%%%%%%%%%%%%%%%%%%%%%%%%%%%%%%%%%%%
\subsection{Interpolation perspective} \label{sec:loewinterp}
%%%%%%%%%%%%%%%%%%%%%%%%%%%%%%%%%%%%%%%%%%%%%%%%%%%%%%%%%%%%%%%%%%%%%%%%%%%%%%%%

Suppose we have underestimated the number of eigenvalues in $\Omega$, 
collecting sample data at $r<m$ left and right points.  
Does the resulting model we have just constructed bear any resemblance 
to the full-order transfer function $\BH(z) = \BV(z\BI-\BLambda)^{-1}\BW^*$?
Model reduction gives some insight.

Suppose we have contour integral data at $r$ left and right interpolation points,
resulting in an $r\times r$ Loewner matrix $\L$ having full rank, $r$.
In this case, we can skip the economy-sized SVD, and simply construct
from the data the transfer function
 \begin{equation} \label{eqn:Grom}
\BG(z) := \rdm\left(\L_s - z@\L \right)^{-1}\ldm.
\end{equation}
If $r<m$, $\BG(z)$ is a rational function generally with lower order than $\BH(z)$, 
and its poles cannot exactly match all poles of $\BH(z)$, 
the desired eigenvalues.
However, if $\L_s - \rp_j\L$ is invertible,%
\footnote{When $r<m$, the eigenvalues of the Loewner pencil $(\L_s,\L)$ 
can potentially fall outside $\Omega$.}
$\BG(z)$ \emph{tangential interpolates}
$\BH(z)$ at $z=\rp_j$ in the direction~$\Br_j$:
\[ \BG(\rp_j) \Br_j = \BH(\rp_j)\Br_j = \rd_j.\]
Similarly, if $\L_s - \lp_i\L$ is invertible, then
\[ \Bell_i^*\BG(\lp_i)  = \Bell_i^*\BH(\lp_i) = \ld_i^*.\]
To verify this interpolation property,
arrange the right and left interpolation points into the diagonal matrices
$\SigmaPts := {\rm diag}(\rp_1,\ldots, \rp_r)$ and $\ThetaPts := {\rm diag}(\lp_1,\ldots,\lp_r)$.
Then the structure of the Loewner and shifted Loewner matrices~\cref{eqn:Lmulti} 
and~\cref{eqn:Lsmulti} can be encoded in the Sylvester equations~\cite{MA07}
\[ \L\SigmaPts - \ThetaPts@ \L = \BL^*\rdm - \ldm@\BR, \qquad
   \L_s\SigmaPts - \ThetaPts@ \L_s = \BL^*\rdm\SigmaPts - \ThetaPts@\ldm@\BR,\]
which can be combined to show~\cite[prop.~3.1]{MA07}
\begin{equation}
\L_s - \L\SigmaPts = \ldm@\BR, \qquad
\L_s - \ThetaPts@\L = \BL^*\rdm.
\end{equation}
Letting $\Be_j$ denote the $j$th column of the $r\times r$ identity matrix,
\begin{align*}
   \BG(\rp_j)\Br_j &= \rdm(\L_s - \rp_j\L)^{-1}\ldm@\BR@\Be_j \\[.25em] 
                   &= \rdm(\L_s - \rp_j\L)^{-1}(\L_s-\L\SigmaPts)@\Be_j \\[.25em] 
                   &= \rdm(\L_s - \rp_j\L)^{-1}(\L_s-\rp_j\L)@\Be_j 
                   = \rdm@\Be_j = \rd_j = \BH(\rp_j)@\Br_j.
\end{align*}
The left interpolation property at $\lp_i$ follows similarly.

This framework is commonly used in data-driven reduced-order modeling, 
where the goal is to construct a reduced-order approximation $\BG(z)$ 
to a high-order transfer function $\BH(z)$ using only evaluations 
of $\BH(z)$, without access to a state-space representation. 
Using $\BG(z)$ in place of $\BH(z)$ can make tasks such as 
simulation and control computationally feasible for large-scale systems.
The poles of $\BG(z)$ will not precisely match those of $\BH(z)$, 
though one might expect certain poles to be good approximations
(depending on the interplay of interpolation points and system dynamics).
See~\cite{ABG10,ALI18,BG12,GAB08,MA07} and the references therein
for more details on rational interpolation in model reduction, 
the Loewner formulation, and optimal interpolation point selection.

%%%%%%%%%%%%%%%%%%%%%%%%%%%%%%%%%%%%%%%%%%%%%%%%%%%%%%%%%%%%%%%%%%%%%%%%%%%%%%%%
\subsection{Matching interpolation points utilize derivative data}  \label{sec:derivdata}
%%%%%%%%%%%%%%%%%%%%%%%%%%%%%%%%%%%%%%%%%%%%%%%%%%%%%%%%%%%%%%%%%%%%%%%%%%%%%%%%
Thus far we have assumed the left interpolation points are disjoint
from the right interpolation points, but this constraint can be relaxed.  
To illustrate how this changes the set-up we have just described, 
suppose $\lp_i=\rp_j$ for some $i,j \in \{1,\ldots,r\}$.
In the limit as $\lp_i \to \rp_j$, the definition of 
$\L$ in \cref{eqn:Lmulti} suggests that $\lp_i=\rp_j$ will require 
information about the derivative $\BH'(z)$ at $z=\lp_i=\rp_j$.
In addition to the interpolation data \cref{eqn:intdata}, 
we will require scalar measurements of the \emph{tangential derivative}
\begin{equation} \label{eqn:intdatahermite} 
 \Bell_i^* \BH'(\rp_j)\Br_{j} \in  \C. 
\end{equation}
\Cref{thm:rationalfunctionCI} shows that this data can also be obtained 
from a contour integral (and thus approximated by quadrature):
\begin{equation}  \label{eq:Hermitedata}
 \Bell_i^*  {\BH'(\rp_j)}\Br_j 
 = \frac{-1}{2\pi@\iop} 
    \int_{\partial \Omega} \frac{1}{(\rp_j-z)^{2}}\,\Bell_i^* \BT(z)^{-1}\Br_j\,\dop z
 \approx 
    \sum_{k=1}^N  \frac{-w_k}{(\rp_j-z_k)^{2}}\,\Bell_i^* \BT(z_k)^{-1}\Br_j.
\end{equation}
The matrices $\ldm$ and $\rdm$ are defined just as in~\cref{eqn:YZmulti} 
using the one-sided data~\cref{eqn:leftdata} and~\cref{eqn:rightdata}.
The Loewner matrix $\L$ in \cref{eqn:Lmulti} and the shifted-Loewner matrix $\L_s$ 
in~\cref{eqn:Lsmulti} are also defined as before, 
except the $(i,j)$ entries are replaced by
\begin{equation} \label{eq:HermiteLoew}
   [\L]_{i,j} = \Bell_i^*\BH'(\rp_j)\Br_j, \qquad 
   [\L_s]_{i,j} =
		\Bell_i^*\left[\rp_j\BH'(\rp_j)  + \BH(\rp_j)\right]\Br_j. 
\end{equation}
One then proceeds as before, either performing SVD truncation 
to recover $\BH(z)$ as in~\cref{eqn:Hmultipoint},	 
or constructing a reduced-order rational interpolant $\BG(z)$ as in~\cref{eqn:Grom}. 
Only in the latter case is there a change from the earlier discussion:
the reduced-order rational function $\BG(z)$ still satisfies the left and right
interpolation conditions $\Bell_i^*\BG(\lp_i) = \Bell_i^*\BH(\lp_i)$ 
and $\BG(\rp_j)\Br_j = \BH(\rp_j)\Br_j$,
but now also the Hermite interpolation condition
	\[
\Bell_i^*\BG'(\rp_j)\Br_j = \Bell_i^*\BH'(\rp_j)\Br_j.
	\]
See~\cite[sect.~6]{MA07} for additional details.
Such tangential Hermite conditions play an important role in  model reduction,
arising as  necessary conditions for $\BG(z)$ to optimally approximate $\BH(z)$ 
in the least-squares sense; see~\cite{ABG10,GAB08} for details.
However, in the context of contour algorithms for eigenvalues, 
we usually presume we have sufficient data to recover $\BH(z)$ completely.

%%%%%%%%%%%%%%%%%%%%%%%%%%%%%%%%%%%%%%%%%%%%%%%%%%%%%%%%%%%%%%%%%%%%%%%%%%%%%%%%
\subsection{Numerical illustration of the multi-point Loewner algorithm}
%%%%%%%%%%%%%%%%%%%%%%%%%%%%%%%%%%%%%%%%%%%%%%%%%%%%%%%%%%%%%%%%%%%%%%%%%%%%%%%%

As emphasized in \Cref{sec:contmultipoint}, 
the primary computational burden of contour based methods results 
from computing the quadrature data, $ \Bell^*\BT(\qn_k)^{-1}$ and $\BT(\qn_k)^{-1}\Br$. 
Once this data is computed, there is negligible additional cost to using 
additional expansion points $\rp_j$ for the same sampling  directions $\Bell$ and $\Br$.
To this end, in the numerical experiments below $r_\rp$  denotes the distinct number of probing 
(direction) vectors, and  we use all these vectors for each interpolation point $\rp_j$. 
This set-up is equivalent to repeating each interpolation point $\rp_j$ as many times 
as there are directions. Given $K_\rp$ distinct interpolation points and $r_\rp$ 
distinct probing vectors, the Loewner matrices then have size $r\times r = K_\rp r_\rp \times K_\rp r_\rp$.
We further take $\lp_j = \rp_j$ for $j=1,\ldots, r$, 
thus using the Hermite formulation of \Cref{sec:derivdata}.

\begin{figure}[h!]
	\begin{center}
		$\overset{K_\sigma = 2,\ r_\sigma = 16}{\includegraphics[width=.32\textwidth]{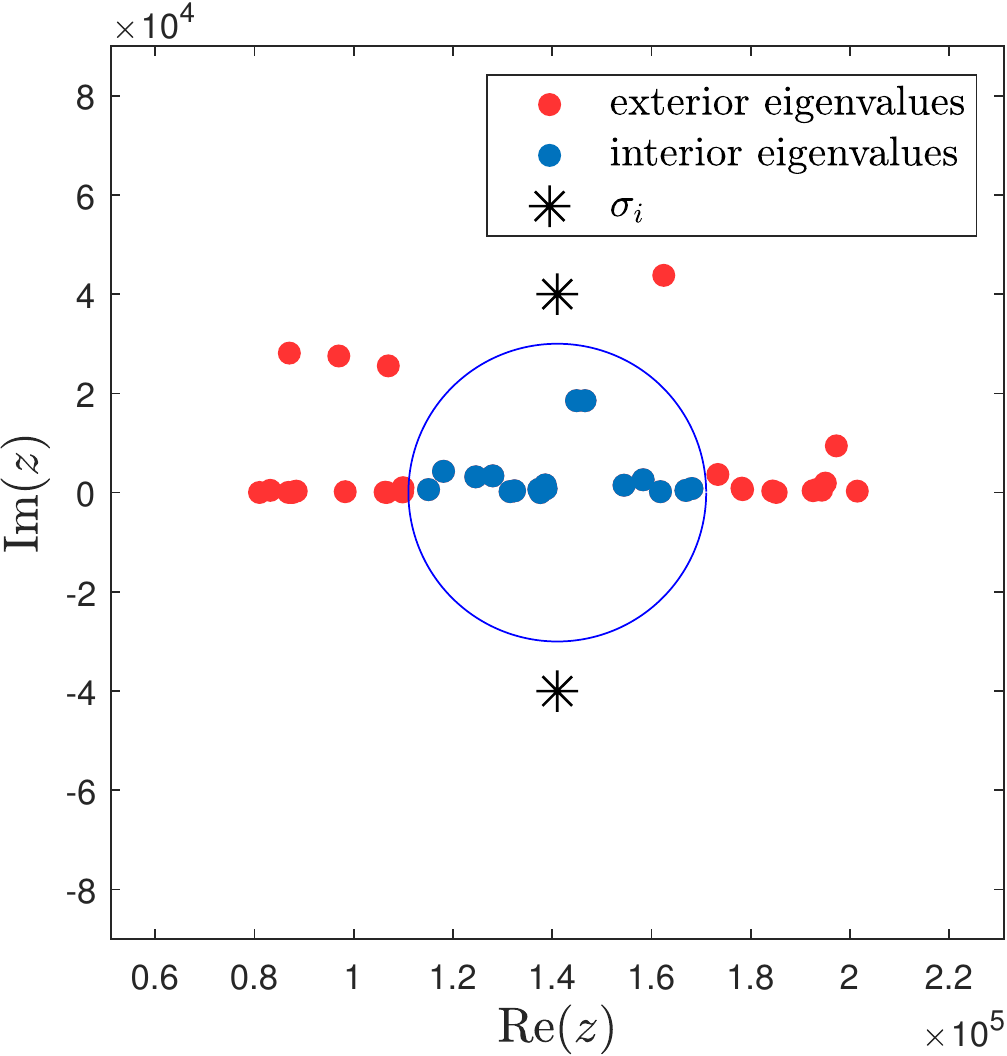}}$
		$\overset{K_\sigma = 4,\ r_\sigma = 8}{\includegraphics[width=.32\textwidth]{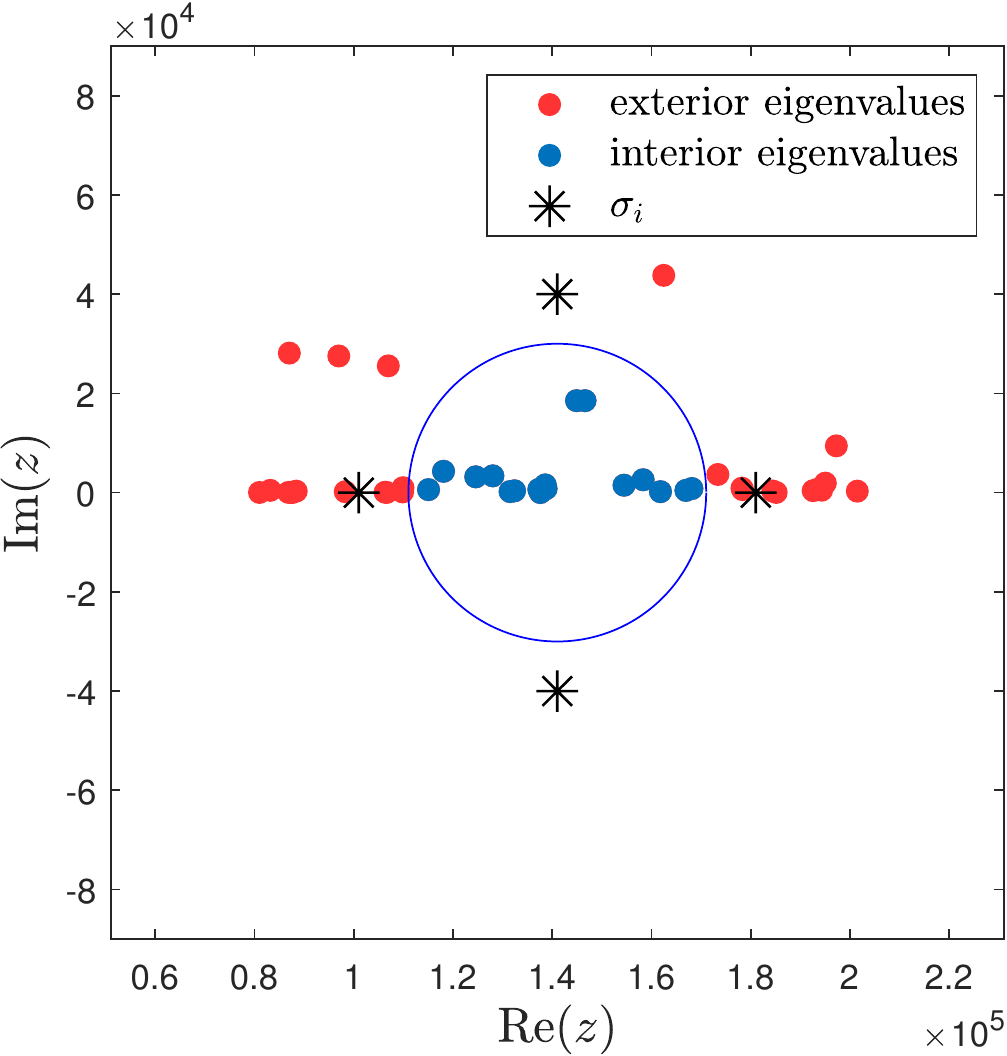}}$
		$\overset{K_\sigma = 8,\ r_\sigma = 4}{\includegraphics[width=.32\textwidth]{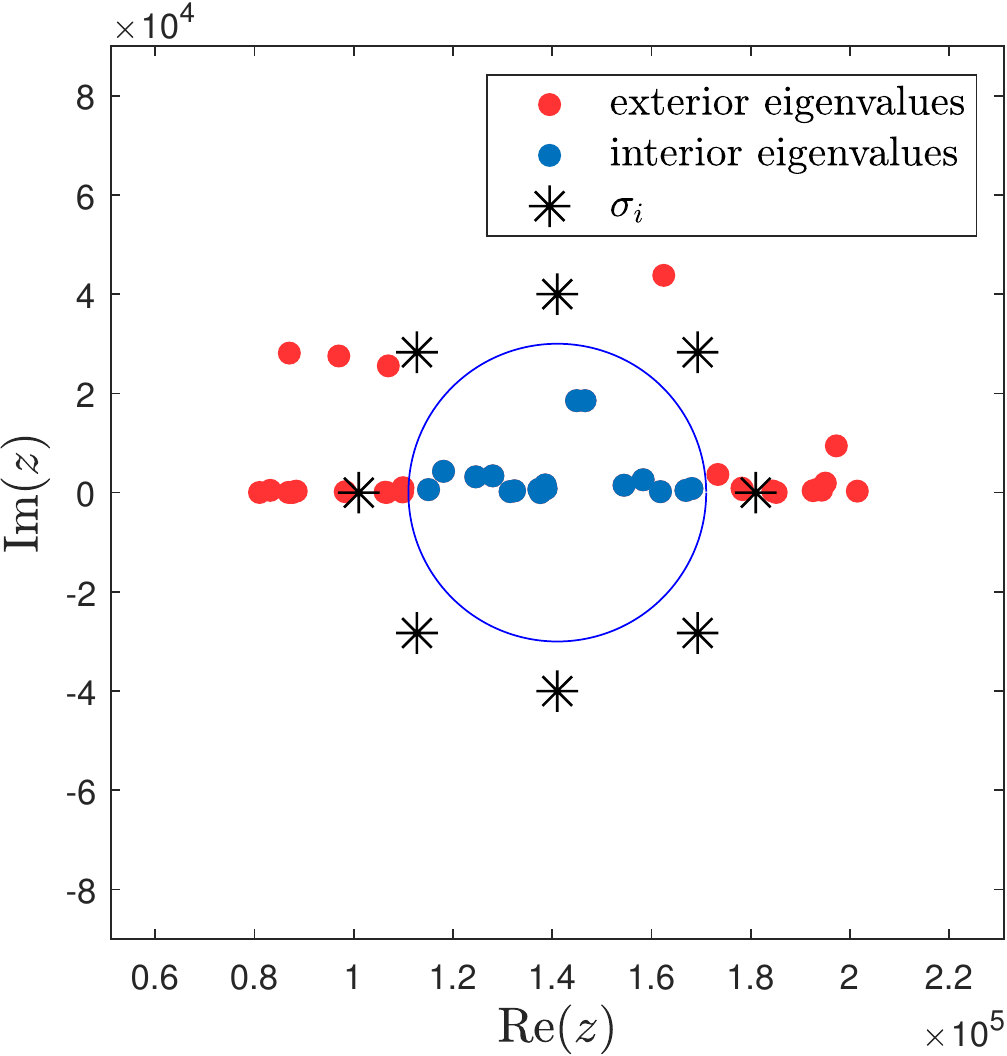}}$

		\includegraphics[width=0.32\textwidth]{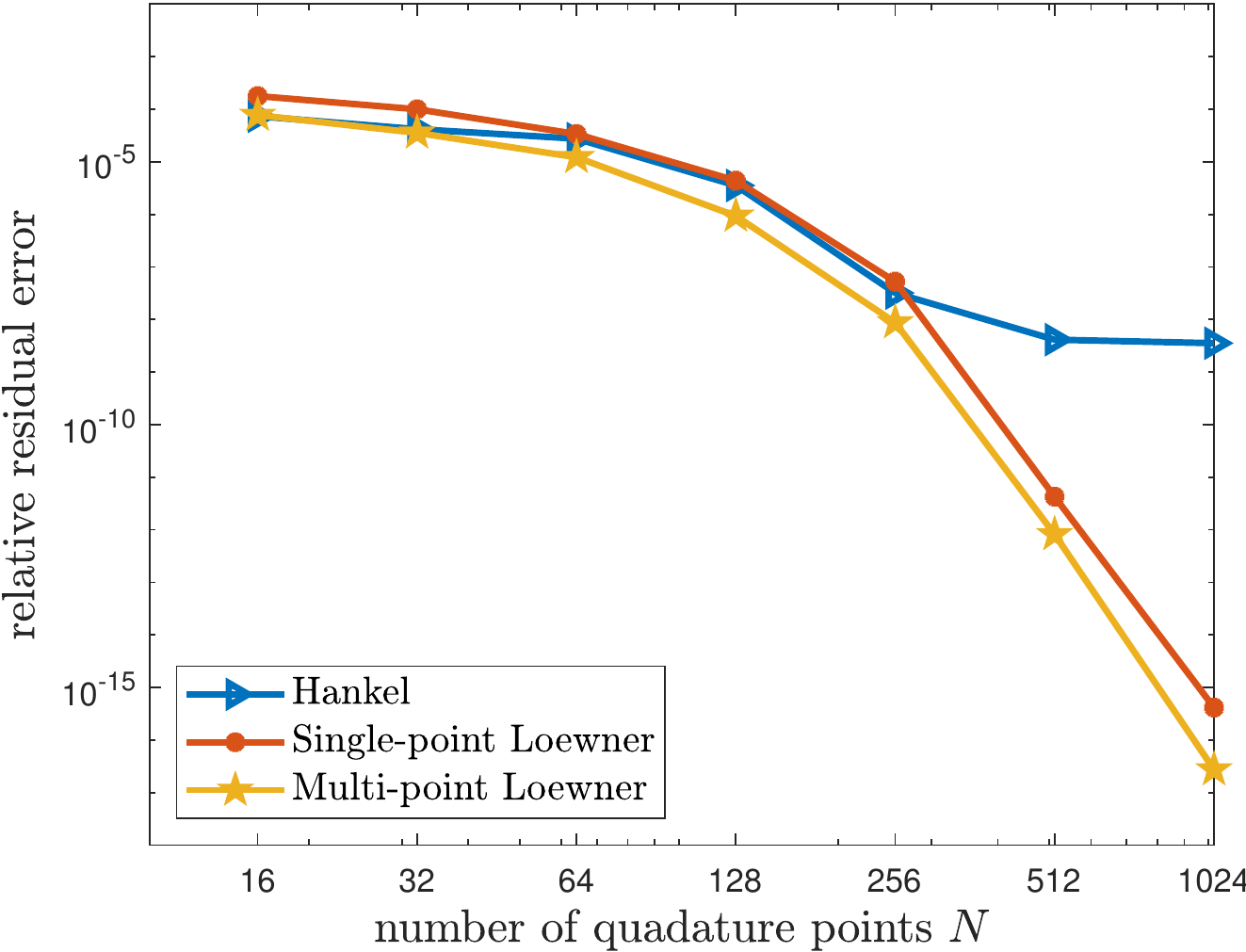}
		\includegraphics[width=0.32\textwidth]{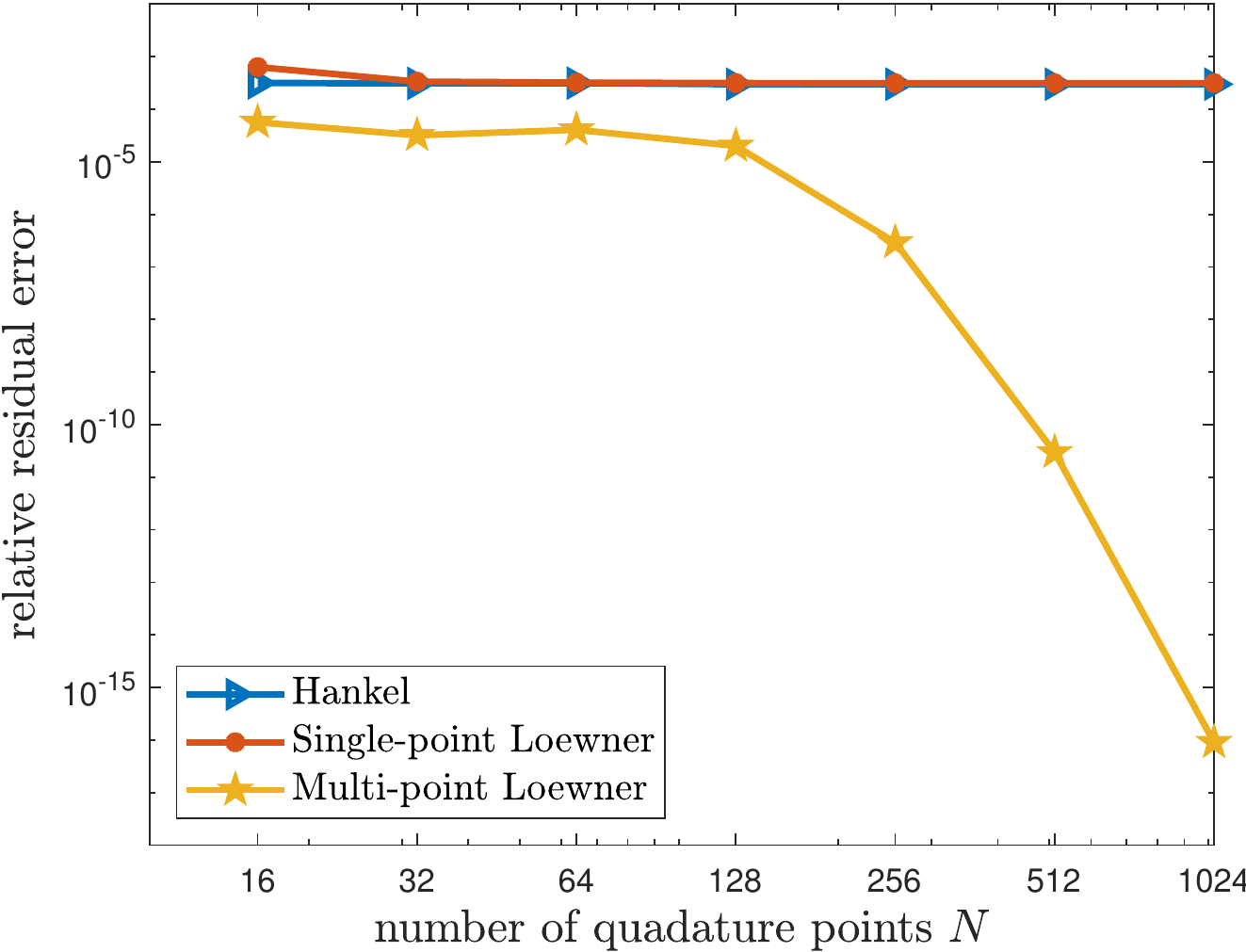}
		\includegraphics[width=0.32\textwidth]{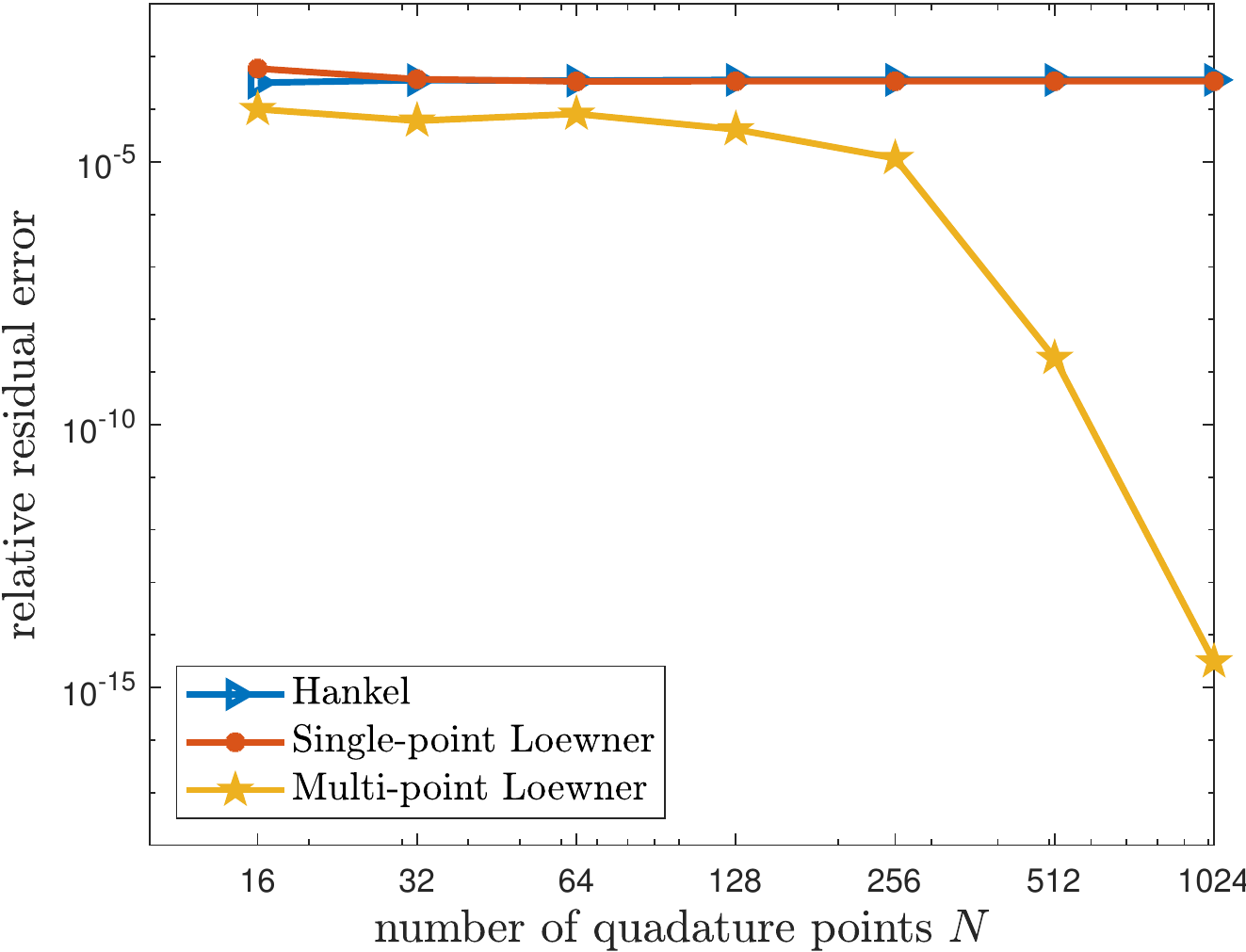}

		\includegraphics[width=0.32\textwidth]{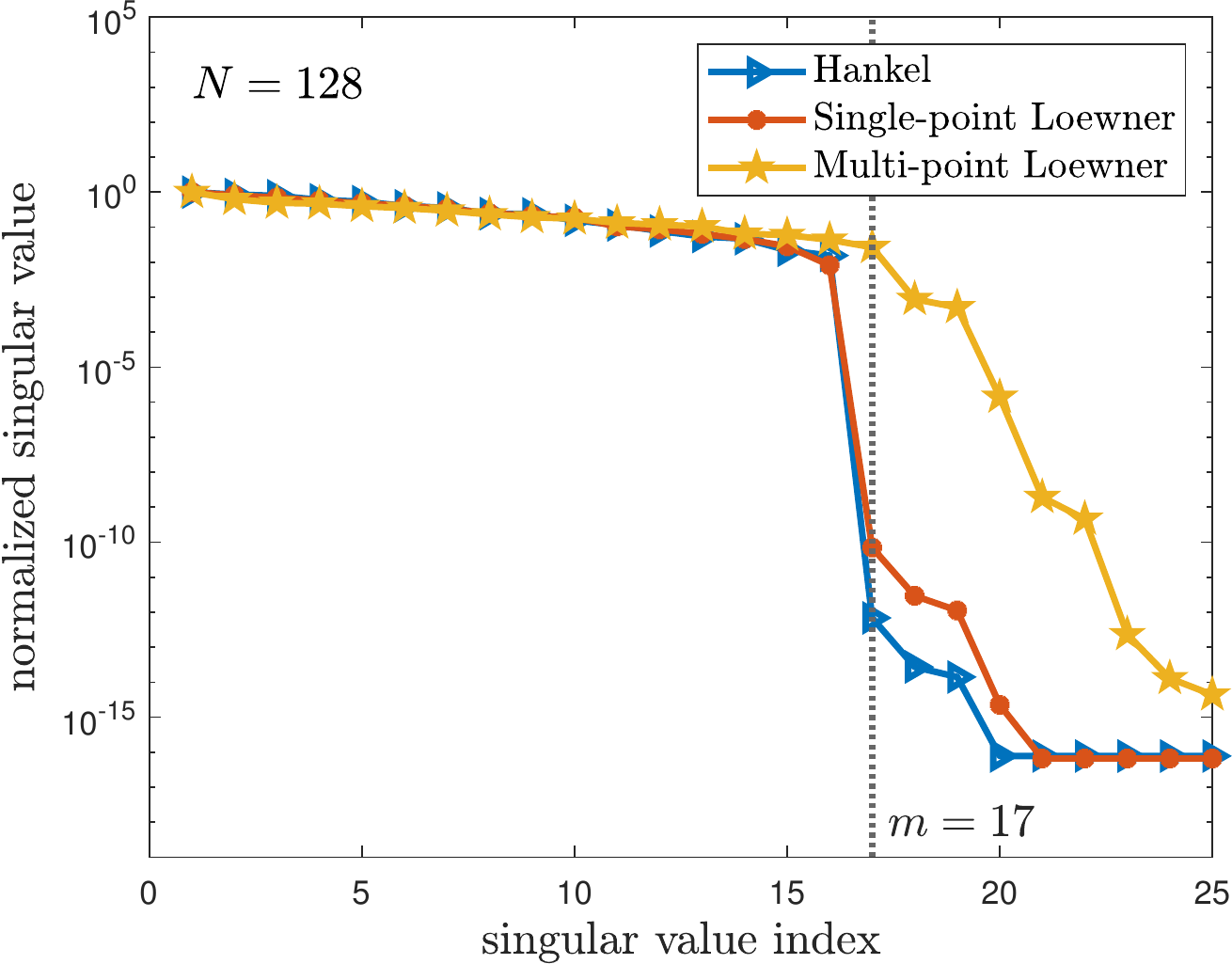}
		\includegraphics[width=0.32\textwidth]{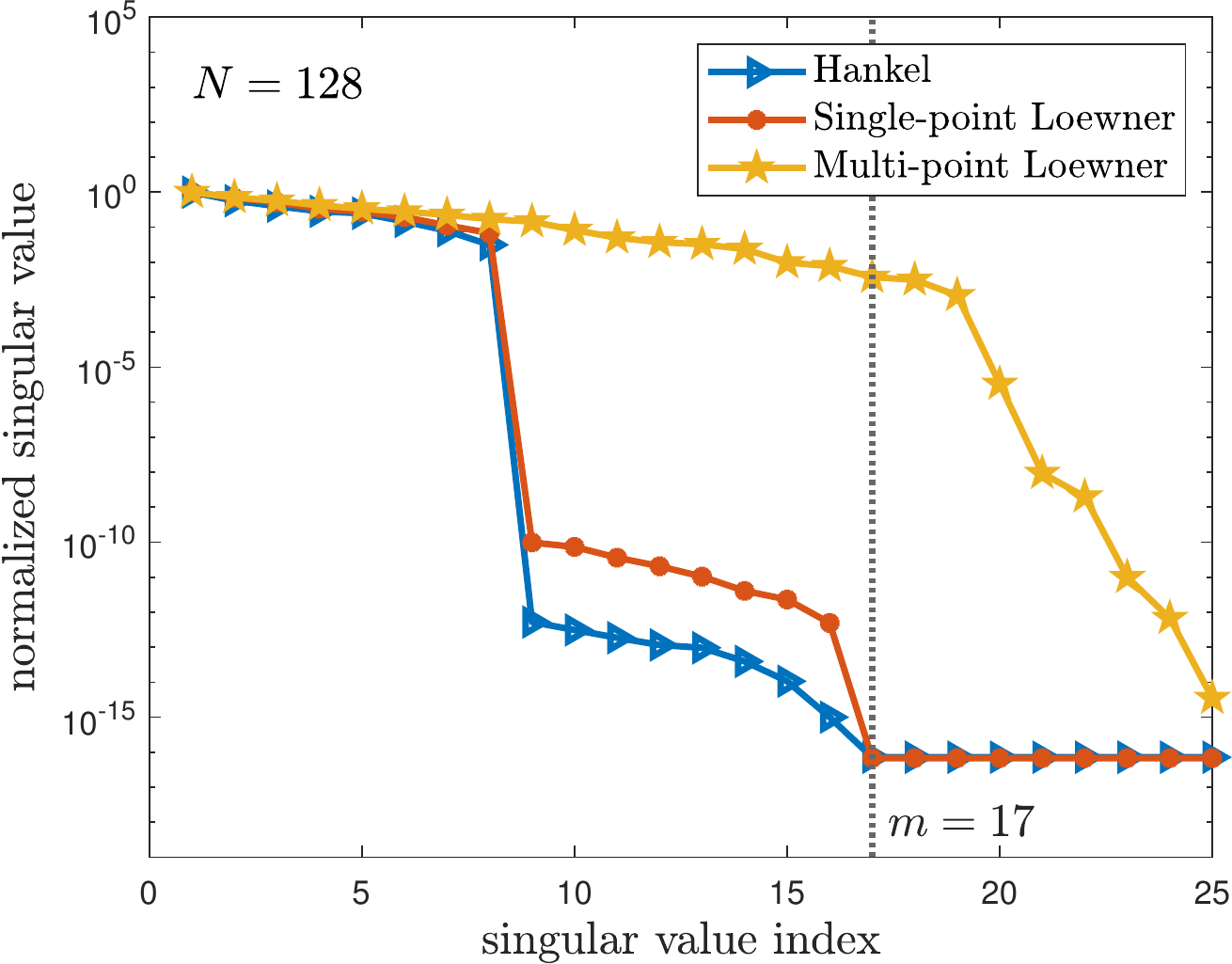}
		\includegraphics[width=0.32\textwidth]{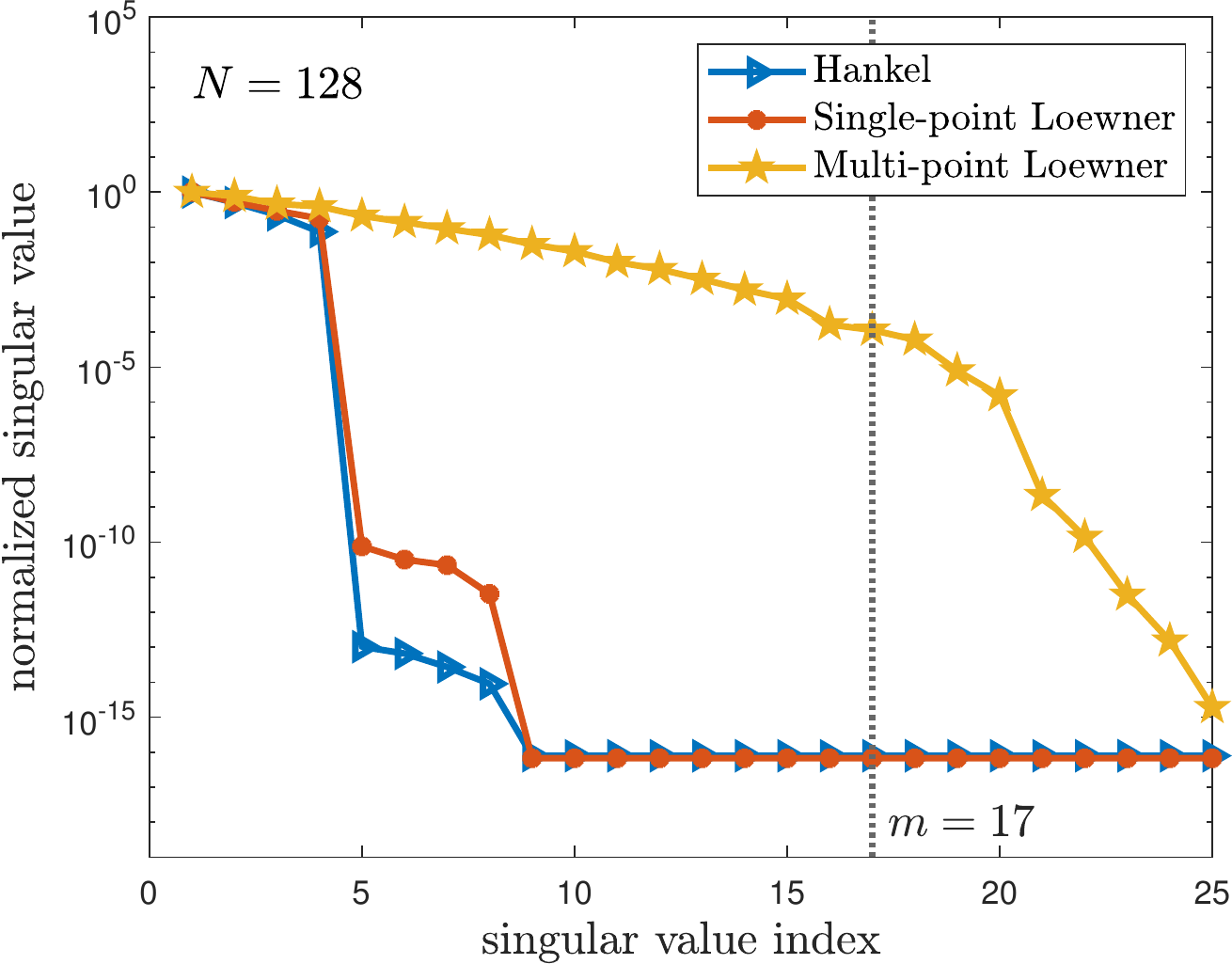}
		
		\includegraphics[width=0.32\textwidth]{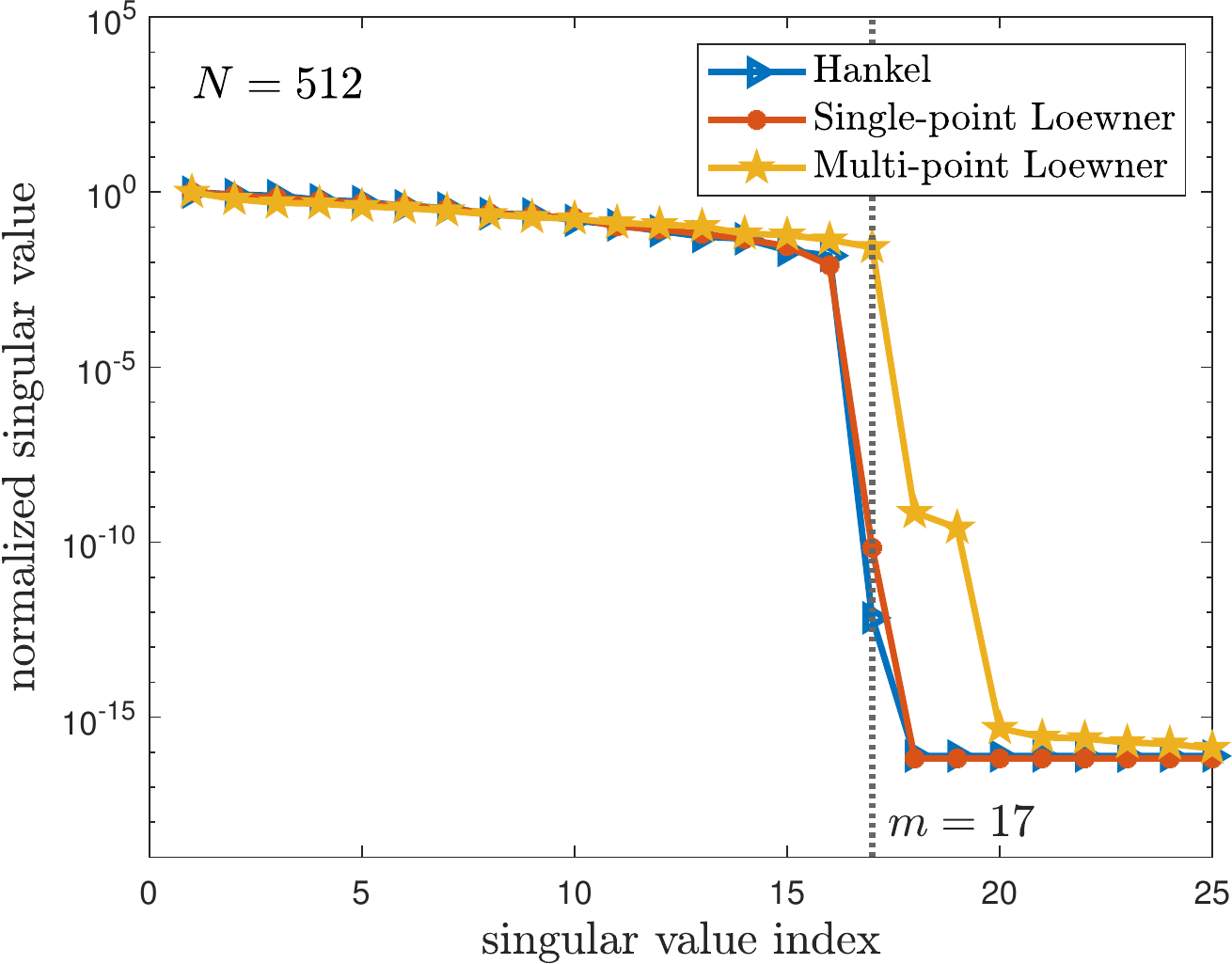}
		\includegraphics[width=0.32\textwidth]{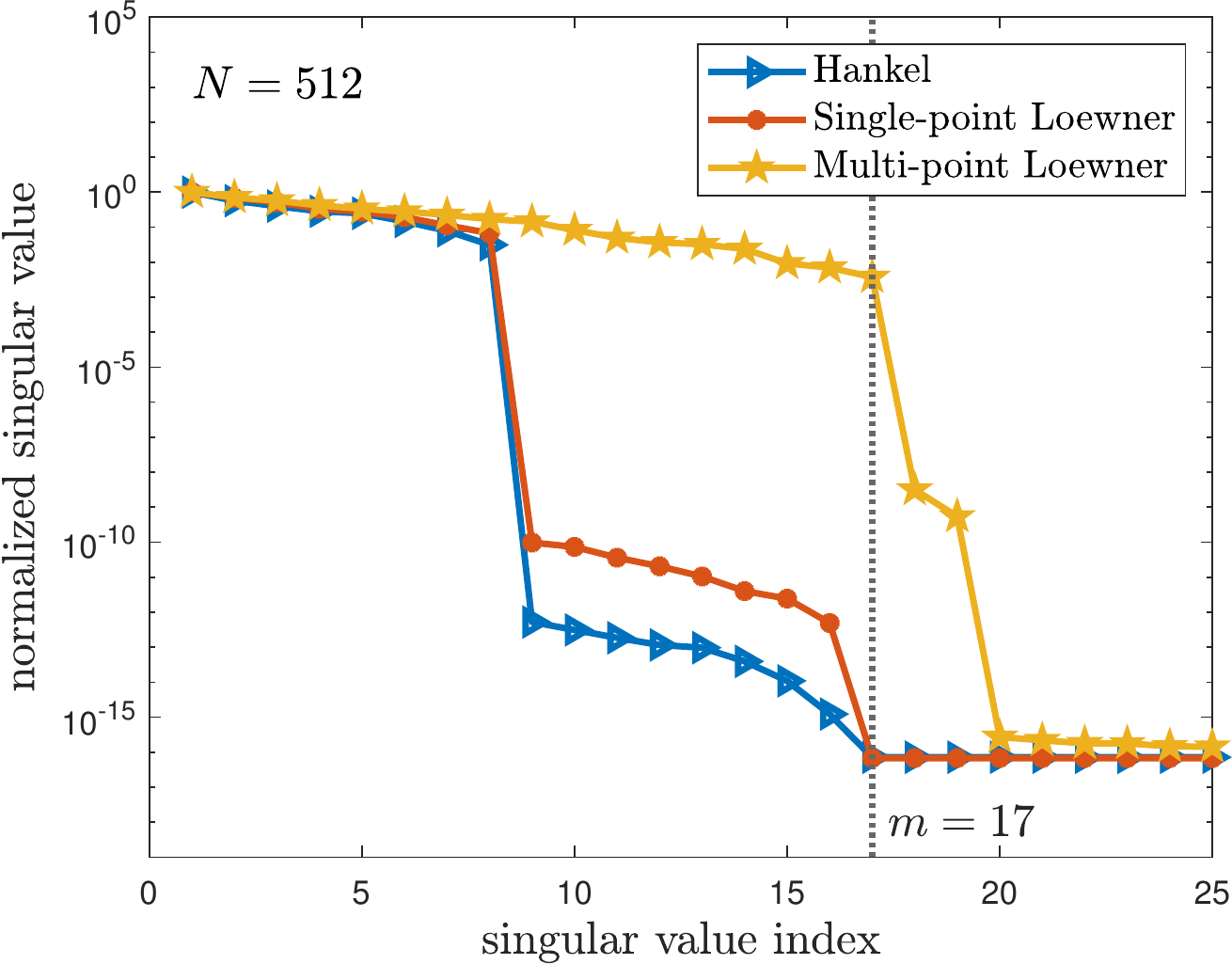}
		\includegraphics[width=0.32\textwidth]{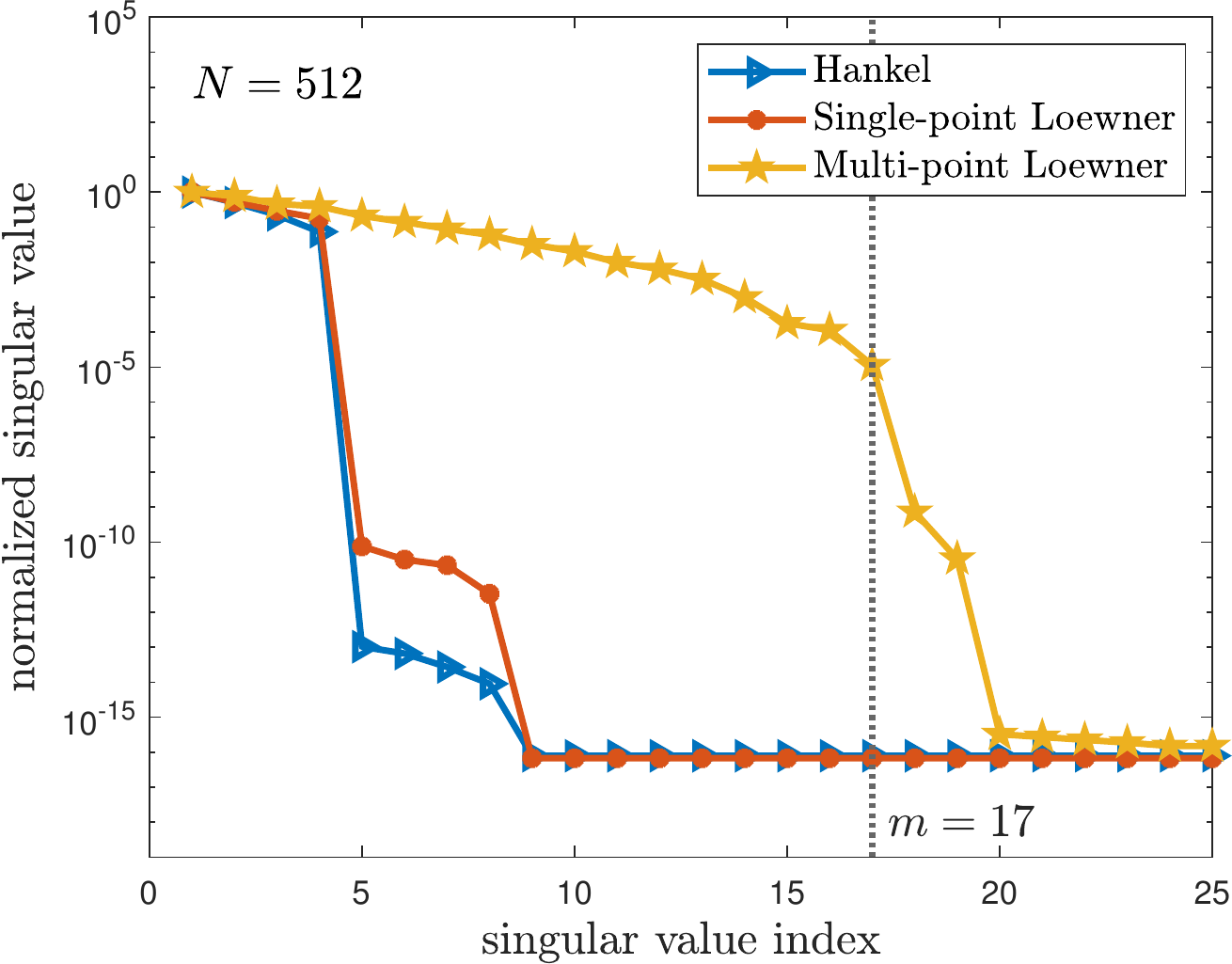}
	\end{center}
	
	\vspace*{-1em}
	\caption{\label{fig:mp_gun}
     {\tt gun} problem, $m=17$ eigenvalues in $\Omega$, Hankel, single-point Loewner, and multi-point Loewner experiments described in the text, all involving $32\times 32$ Hankel and Loewner matrices:
multi-point Loewner interpolation points (top row);
maximum relative residual error (second row);
singular values of $\H$ and $\L$ with $N=128$ (third row) and $N=512$ (bottom) quadrature points.
	}
\end{figure}

We now provide results for the multi-point Loewner method applied to the \verb|gun| problem described 
in \Cref{sec:singlept:numerics}, using the same contour $\partial \Omega$ and interpolation point 
$\sigma$ for single-point Loewner as before. 
For multi-point Loewner, we consider interpolation points in a circle concentric to $\partial \Omega$, 
but of larger radius $40000$.
\Cref{fig:mp_gun} shows the placement of the interpolation points compared to the contour and eigenvalues. 
We take the number of interpolation points $K_\rp = K$ and the number of probing directions $r_\rp$ 
so that both the Hankel and Loewner matrices are  $32 \times 32$ for each method. 
In the cases where $r_\rp = 4$ ($K_\rp = 8)$ and $r_\rp = 8$ ($K_\rp = 4$), 
we see that as the number of quadrature points, $N$, increases, the approximate eigenvalues from the multi-point Loewner approach converge to the true eigenvalues. On the other hand, the eigenvalue estimates from both the Hankel and single-point Loewner methods fail to converge. The singular values of the Hankel and Loewner matrices shed light on this difference: unlike the Hankel and single-point Loewner methods,
the multi-point approach does not underestimate the rank $m$; 
as the number $N$ of interpolation points increases, the correct rank $m=17$ is revealed for the multi-point case, whereas the other methods give a numerical rank approximately equal to $r_\sigma < m$.

%%%%%%%%%%%%%%%%%%%%%%%%%%%%%%%%%%%%%%%%%%%%%%%%%%%%%%%%%%%%%%%%%%%%%%%%%%%%%%%%

%% file: rationalforTinv.tex
%!TEX root = paper.tex

\section{Approximate eigenvalues via rational approximation of {\boldmath$\BT(z)^{-1}$}}
\label{sec:ratTinv}

To this point we have addressed contour integral methods that compute 
all eigenvalues in a prescribed domain, $\Omega\subset \C$.\ \  
With exact contour integral evaluations, the eigenvalues would be found exactly.
Quadrature gives noisy data. If $N$ quadrature points $\{\qn_k\}_{k=1}^N$ are used,  
one must evaluate $N$ quantities of the form $\Bell^*\BT(\qn_k)^{-1}$ and $\BT(\qn_k)^{-1}\Br$.

Here, we briefly introduce another approach to NLEVPs that uses similar 
systems theory techniques, but seeks to use fewer evaluations involving $\BT(z)^{-1}$ 
to develop (potentially crude) approximations to eigenvalues of $\BT(z)$.
This alternative method approximates $\BT(z)^{-1}$ with a rational function $\BG(z)$,
then uses the  poles and residues of $\BG(z)$ to approximate the eigenvalues and 
eigenvectors of $\BT(z)$.  In contrast to many existing methods, we do not 
approximate $\BT(z)$ via a rational or polynomial function, 
then solve the resulting (rational or polynomial) eigenvalue problem. 
Instead, we directly approximate $\BT(z)^{-1}$ in a manner that
delivers a linear eigenvalue problem. 

First consider the linear (generalized) eigenvalue problem: $\BT(z) = z \BE - \BA$ 
with $\BE,\BA \in \C^{n\times n}$.
The resolvent $\BT(z)^{-1} = (z \BE - \BA)^{-1}$ is the transfer function 
for a degree-$n$ dynamical system with $n$-inputs and $n$-outputs.
Now construct an order $r \leq n$ rational approximant
$\BG(z) := \Cred(z \Ered - \Ared)^{-1}\Bred$,
where $\Ered,\Ared \in \C^{r \times r}$, $\Bred \in \C^{r \times n}$, 
and $\Cred \in \C^{n \times r}$.
Take the poles of $\BG(z)$ as approximations to eigenvalues of $\BT(z)$.

This same approach can be applied to NLEVPs, 
using only evaluations of $\BT(z)^{-1}$ at selected points. 
While in \Cref{sec:multipoint} we used the data-driven Loewner framework 
to sample and recover
the \emph{rational} transfer function $\BH(z) = \BV(z\BI-\BLambda)^{-1}\BW^*$,
that rational interpolation methodology does not require the sampled function to be rational,
if the goal is only approximation and not exact recovery. 
Choose {left interpolation points} $\lp_1, \ldots, \lp_r \in \C$ 
with left direction vectors
$\Bell_1, \ldots, \Bell_r \in \C^n$,
and {right interpolation points} $\rp_{1}, \ldots, \rp_{r}\in \C$ 
with right direction vectors
$\Br_{1}, \ldots,\Br_{r} \in \C^n$. 
Assume here that $\rp_i \ne \lp_j$ for $i,j \in\{ 1,\dots, r\}$, 
and suppose none of these points is an eigenvalue of $\BT(z)$.
Compute the probed (tangential) samples of $\BT(z)^{-1}$ at these points: 
\begin{equation} \label{eqn:intdataforGr} 
\hspace*{-2em}\ld_i^* = \Bell_i^* \BT(\lp_i)^{-1} \in \C^{1\times n},\qquad
\rd_j = \BT(\rp_{j})^{-1} \Br_{j} \in \C^n,
\rlap{\qquad $i,j = 1,\ldots,r$.}
\end{equation}
Construct  $\BG(s)$ as in \cref{eqn:Lmulti}--\cref{eqn:YZmulti} and~\cref{eqn:Grom},
but now using probed samples of $\BT(z)^{-1}$ in \cref{eqn:intdataforGr} instead of $\BH(z)$, 
i.e., $\BG(z) = \Cred(z \Ered - \Ared)^{-1}\Bred$, where, for $i,j=1,\ldots,r$,
\[
	[\Ered]_{i,j} = -[\L]_{i,j}
                 = -\,\frac{\ld_i^*\Br_j - \Bell_i^\ast \rd_j}{\lp_i - \rp_j},
   \qquad
	[\Ared]_{i,j} = -[\L_s]_{i,j} = 
	-\,\frac{\lp_i@\ld_i^*\Br_j - \rp_j\Bell_i^\ast \rd_j}{\lp_i - \rp_j},
\]
\[ \Bred = \left[\,\ld_1 ~~\ld_2~~ \cdots ~~ \ld_r \,\right]^* \in \C^{r  \times n}, 
   \qquad
   \Cred = \left[\,\rd_1 ~~\rd_2 ~~ \cdots ~~ \rd_r \,\right]\in \C^{n \times r}.
\]
Assuming $\lp_i@\L-\L_s$ and $\rp_j@\L-\L_s$ are all invertible,
$\BG(z) $ is a low-order rational interpolant to $\BT(z)^{-1}$. 
The poles of $\BG(z)$ approximate poles of $\BT(z)^{-1}$, 
eigenvalues of $\BT(z)$:
if $(\lambda,\Bz)$ is an eigenpair of the pencil $(\Ared,\Ered)$,
then $(\lambda,\Cred@\Bz)$ is an \emph{approximate} eigenpair of $\BT(z)$.

Notice that this approach requires $r$ evaluations each of quantities of the form
$\Bell^*\BT(\lp)^{-1}$ and $\BT(\rp)^{-1}\Br$, compared to $N$ evaluations
of similar quantities $\Bell^*\BT(\qn_k)^{-1}$ and $\BT(\qn_k)^{-1}\Br$ required
for contour integral methods.  When $r \ll N$, the direct
interpolation approach described in this section will be much cheaper to execute.

As stated, there are no guarantees on the accuracy of the approximate eigenpairs
extracted from this method. 
Their accuracy will depend strongly on the selection of interpolation data. Interpolation points should be chosen in the region of interest. 
The direction vectors are also crucial, 
and the closer they come to an actual eigenvector, 
the better the approximation should be. 
Since  the matrix $\BT(\lambda)$ is singular at a true eigenvalue $\lambda$, 
for a given right interpolation $\rp_j$, the smallest right singular vector of $\BT(\rp_j)$ 
is an appealing choice for the right direction $\Br_j$. 
One can perform a few steps of an iterative algorithm (such as inverse iteration) 
to approximate this vector. 
A better strategy will be to iteratively correct the selection of interpolation points 
and directions, akin to the similar procedure in the Dominant Pole Algorithm \cite{RN06,RS08} 
and the Iterative Rational Krylov Algorithm \cite{GAB08}. 
These approximate eigenpairs can also be used as a preprocessing step for 
contour integral methods, to estimate the location and number of eigenvalues;
corresponding approximate eigenvectors could then be used as probing directions.  
We will pursue these consideration in future work, and refer the reader 
to \cite{Brennan18} for an initial investigation.

%% file: filterfunctions.tex
%!TEX root = paper.tex

\section{Filter functions in numerical approximation of  contour integrals}
\label{section:filterfunctions}
Contour integral methods (both established Hankel methods
and the Loewner approaches we have described) require integrals of the form~\cref{eq:fprobe}. 
In practice these integrals are approximated via numerical quadrature as in~\cref{eq:trap},
leading to inexact data.
How do these quadrature errors interact with the terms in the Keldysh decomposition~\cref{eq:keldysh}?
Van Barel and Kravanja~\cite{BK16,Bare16} have studied this question for the Hankel approach using the
concept of \emph{rational filter functions}.  Here we briefly summarize the analysis 
from~\cite{BK16,Bare16}, and indicate how it could be extended to the Loewner setting.

Approximate the 
Markov parameter $\BA_\pow = \BL^*\BV\BLambdas{\pow}\BW^*\BR$ in \cref{eq:ApLR} via quadrature: 
\[
\BA_\pow \approx \widetilde{\BA}_\pow := \sum_{j=1}^{N} \qw_j @ \qn_j^\pow \,\BL^*\BT(\qn_j)^{-1} \BR,
\]
where $\{\qw_j\}_{j=1}^N$ and $\{\qn_j\}_{j=1}^N$  are the quadrature weights and nodes.
In the Keldysh decomposition $\BT(z)^{-1} =\BH(z) + \BN(z)$ in~\cref{eq:keldysh}, 
\[
\BH(z) = \BV(z\BI-\BLambda)^{-1}\BW^* = \sum_{i = 1}^{m}  \dfrac{\Bv_i^{} \Bw_i^* }{z - \lambda_{i}}.
\]
The approximate Markov parameter $\widetilde{\BA}_\pow$ is then 
\begin{align}
\widetilde{\BA}_\pow 
    =& \sum_{j=1}^{N} \qw_j \qn_j^\pow @@@\BL^*\BH(\qn_j) \BR 
        + \sum_{j=1}^{N} \qw_j \qn_j^\pow @@@\BL^* \BN(\qn_j) \BR \nonumber \\ 
    =& \sum_{i=1}^{m}\bigg(\BL^* \Bv_i^{} \Bw_i^* \BR \sum_{j=1}^{N} \dfrac{\qw_j \qn_j^\pow}{\qn_j - \lambda_{i}} \bigg)
        + \sum_{j=1}^{N} \qw_j \qn_j^\pow @@@\BL^*\BN(\qn_j) \BR.  \label{Akapprox}
\end{align}
In the last expression, the rational function
\begin{align} \label{hfilterk} 
b_\pow(z) := \sum_{j=1}^{N} \dfrac{\qw_j \qn_j^\pow}{\qn_j - z},
\rlap{\qquad $\pow=0,1,\ldots.$}
\end{align}
is called a \emph{rational filter function.} 
To distinguish it from the Loewner case, we call it the {\it Hankel filter function}.%
\footnote{We emphasize that the use of \emph{rational filter functions} for Hankel methods 
is entirely distinct from our use of \emph{rational interpolation methods} for system realization,
which develops rational approximations and realizations of the function $\BH(z)$. 
}

To simplify the presentation, let the domain $\Omega$ be the unit circle, 
$\Omega = \{z: |z| < 1  \}$, and assume the trapezoidal rule is applied with 
quadrature nodes $\qn_j = e^{2\pi\iop j/N}$ and quadrature weights $\qw_j = \qn_j/N$. 
Then one can show that 
\[
b_0(z) = \frac{1}{N} \sum_{j=1}^{N}\frac{\qn_j}{\qn_j -z} 
       = \frac{1}{1 - z^N},
\qquad\quad
b_\pow(z) = z^\pow @b_0(z), 
\quad \pow=1,2,\ldots.
\]
The filter function  $b_0(z)$ is in fact the trapezoidal rule approximation of the ideal filter, the indicator function of the unit circle as defined by a Cauchy integral:
\[
b_0(\lambda) \approx \frac{1}{2\pi@\iop} \int_{\partial \Omega} \frac{1}{z-\lambda} @@\dop z =
	\left\{\begin{array}{ll}
       1, & |\lambda| < 1; \\[.25em]
       0, & |\lambda| > 1.
   \end{array}\right.
\]
As discussed in \cite{BK16}, even in the case of $\Omega = \{z: |z| < 1  \}$
the choice for the quadrature is not restricted to the trapezoidal rule. 
By viewing  $b_\pow(z)$ in~\cref{hfilterk} as a rational function characterized 
by the poles $\{\qn_j\}_{j=1}^N$ and residues $\{\qw_j\}_{j=1}^N$, 
one can design new rational filter functions to achieve specific goals.
Based on~\cref{Akapprox} and the structure from the trapezoidal rule approximation, 
\cite{BK16,Bare16} propose three goals for rational filter design
(taken from~\cite[p.~349]{Bare16}, with an adjustment to the third condition):
\smallskip
\begin{enumerate}
\item $b_{\pow}(z) = z^\pow b_0(z)$ for $\pow=1,2,\ldots;$
\smallskip
    \item $\| \sum_{j=1}^{N} \qw_j \qn_j^\pow \BL^*\BN(\qn_j)\BR\|$ should be small; 
\smallskip
     \item {$b_0(z)$ should approximate one} inside $\Omega$ 
and be small  {in magnitude} outside $\Omega$.
\end{enumerate}
\smallskip
Van Barel and Kravanja then convert these design goals into an optimization problem;
see~\cite{BK16} for details. 
Regardless of how the filter functions are designed, the guiding objective is to obtain better
approximations of the Markov parameters $\BA_\pow$, thus to improve the performance of 
Hankel contour methods. 

To show how this idea can be extended to the Loewner setting, we generalize
the three filter design goals to the single-point Loewner method from \Cref{sec:singlepoint},
with interpolation point $\sigma\not\in\overline{\Omega}$.\ \  
The quantity $\BL^*\BM_\pow\BR = (-1)^k\BL^*\BV(\rp@\BI-\BLambda)^{-(k+1)}\BW^*\BR$ 
is replaced by a quadrature approximation 
of the contour integral~\cref{eqn:Hviacontour}:
\begin{align*}
	&\BL^*\BM_\pow\BR  
    \approx (-1)^{\pow}
	\sum_{j=1}^{N}
		 \frac{\qw_j}{(\rp-\qn_j)^{\pow+1}} \BL^* \BT(\qn_j)^{-1} \BR\  \\ 
     &=  (-1)^{\pow} \left(\sum_{i=1}^{m}\BL^* \Bv_i^{} \Bw_i^* \BR \sum_{j=1}^{N} 
               \dfrac{\qw_j }{(\qn_j - \lambda_{i})(\rp-\qn_j)^{\pow+1}} + \sum_{j=1}^{N}  \dfrac{\qw_j }{(\rp-\qn_j)^{\pow+1}} \BL^*\BN(\qn_j) \BR\right), \nonumber
		\end{align*}
where, as before, $\{\qw_j\}_{j=1}^N$ and $\{\qn_j\}_{j=1}^N$  denote the quadrature weights  and nodes. 
The Hankel filter function $b_\pow(z)$ in \cref{hfilterk} is now replaced with the \emph{Loewner filter function}
\begin{align}
b_{\rp,\pow}(z) := \sum_{j=1}^{N} \dfrac{\qw_j }{(\qn_j - z)(\rp-\qn_j)^{\pow+1}}, \rlap{\qquad $\pow=0,1,\ldots.$}
\end{align}
Consider, for example, the case $\pow=0$.  Then, 
\[
\BL^*\BM_0\BR 
    \approx \sum_{i=1}^{m} \bigg(\BL^* \Bv_i^{} \Bw_i^*\BR \sum_{j=0}^{N} \dfrac{\qw_j}{(\qn_j - \lambda_i)(\sigma - \qn_j)}\bigg) 
    + \sum_{j=1}^{N}  \frac{\qw_j }{\sigma  - \qn_j}\BL^*\BN(\qn_j)\BR,
\]
yielding the zeroth-order Loewner filter function
\begin{align}
b_{\rp,0}(z) =   \sum_{j=1}^{N} \dfrac{\qw_j }{(\qn_j - z)(\rp-\qn_j)}.
\end{align}
By analogy with the Hankel case {(again using $\Omega = \{z:|z|<1\}$)}, 
$b_{\rp,0}(z)$ approximates the ideal filter 
\[
b_{\rp,0}(\lambda) 
   \approx \frac{1}{2\pi@\iop} \int_{\partial \Omega} \cfrac{1}{(z-\lambda)(\rp-z)} \,\dop z
= 	\left\{\begin{array}{cl}
1/(\rp-z), & |\lambda| < 1; \\[3pt]
  0,       & |\lambda| > 1.
   \end{array}\right.
\]
The ideas developed in~\cite{Bare16,BK16} can be expanded here to design \emph{Loewner filter functions},
requiring modification of the design objectives stated above.
One might consider, for example, the following goals:
\smallskip
\begin{enumerate}
\item ${\displaystyle b_{\rp,\pow}(z) = \dfrac{1}{(\rp-z)^\pow}b_{\rp,0}(z)}$ for $\pow=1,2,\ldots;$
\smallskip
\item ${\displaystyle \bigg\| @@\sum_{j=1}^{N}  \dfrac{\qw_j}{(\rp-\qn_j)^{\pow+1}} \BL^*\BN(\qn_j) \BR@@\bigg\|}$ should be small;
\smallskip
\item {$b_{\rp,0}(z)$ should approximate $1/(\rp-z)$ inside $\Omega$} and be small {in magnitude} outside~$\Omega$.
\end{enumerate}
\smallskip

We have sought to briefly show how rational filter design could be 
generalized to Loewner-based contour methods. 
Thus the improvements (in accuracy and computational speed) 
that filters provide for Hankel methods can be anticipated to similarly
benefit Loewner approaches. 
These issues will be considered in future work. 

%% file: modaltruncation.tex
%!TEX root = paper.tex

\section{Contour-integration methods for data-driven modal truncation}  \label{sec:modal}
So far, we have used tools from systems theory and rational interpolation 
to cast contour integral methods for NLEVPs in the framework of 
data-driven realization, exploiting this perspective to propose new methods 
for solving eigenvalue problems. 
In this section we do the opposite: the machinery behind contour integral methods 
suggests a new data-driven approach for computing reduced order models 
of dynamical systems using modal truncation.

Modal truncation constructs a reduced model for a linear time-invariant dynamical
system by restricting the dynamics to evolve within the span of selected eigenvectors, 
giving a reduced model that only contains the corresponding eigenvalues. 
(For example, one might omit leftmost eigenvalues, or those with large imaginary parts.)
More precisely, let
\begin{equation}
\TF(z) = \Cmod( z\BI - \Amod)^{-1} \Bmod,~~\mbox{where}~~ \Cmod \in \C^{\nout \times \ns},~\Amod \in \C^{\ns \times \ns},~ \mbox{and}~\Bmod \in \C^{\ns \times \nin},
\end{equation}
be the transfer function of a linear dynamical system of degree-$\ns$, with $\nin$ inputs  and $\nout$ outputs.  
For simplicity, 
suppose all the eigenvalues $\pole_1,\ldots,\pole_\ns$ of $\Amod$ are distinct,
permitting us to write $\TF(z)$ in the pole-residue form 
\begin{equation} \label{poleresidue}
\TF(z) = \sum_{j=1}^\ns \frac{\lrd_j^{} \rrd_j^*}{z-\pole_j},~~\mbox{where}~~\lrd_j \in \C^{\nout}~\mbox{and}~\rrd_j \in \C^{\nin}.
 \end{equation}
Modal truncation constructs a reduced  transfer function of order $m<n$
of the form
\begin{equation}  \label{redpoleresidue}
\TFr(z) = \sum_{j=1}^\nred \frac{\lrd_j^{} \rrd_j^*}{z-\pole_j}.
 \end{equation}
In other words, the reduced model is constructed by retaining only the terms
corresponding to the first $\nred$ poles,\footnote{We assume the poles are ordered such that the retained poles are the leading ones.} 
$\pole_1,\ldots,\pole_\nred$, with residues $\lrd_j^{} \rrd_j^*$. 
Construction of  $\TFr(z)$ is typically achieved by computing a spectral decomposition 
of $\Amod$, then truncating the pole-residue form \cref{poleresidue}  to obtain the approximant \cref{redpoleresidue}; thus modal truncation typically requires access to the system's state-space representation, 
i.e., the matrices $\BA$, $\BB$, and $\BC$.\ \ However, 
contour integration tools enable one to perform modal truncation 
using only evaluations of $\TF(z)$,  without access to state-space quantities. 

Let $\TFt(z) $ denote the truncated part  (tail) of the modal decomposition, i.e.,
\[
\TFt(z) = \sum_{j=\nred+1}^\ns \frac{\lrd_j^{} \rrd_j^*}{z-\pole_j},
\]
so that 
\[
\TF(z)  =  \TFr(z) + \TFt(z). 
\]
By sampling $\TFr(z)$ at enough points in the complex domain, then using the Loewner modeling framework 
of \Cref{sec:multipoint}, we could exactly recover $\TFr(z)$, thus performing a data-driven modal truncation. 
However, we have direct access to $\TF(z)$, not $\TFr(z)$: 
precisely the problem the contour integration resolves. 
In the language of \Cref{fig:schematic}, 
$\BT(z)^{-1}$, $\BH(z)$, and $\BN(z)$ are now replaced by $\TF(z)$, $\TFr(z)$, and $\TFt(z)$, respectively. 
Contour integration of $\TF(z)$ enables us to sample $\TFr(z)$ at selected points in the complex plane. 
Let $\Omega$ be a domain containing the  poles (eigenvalues) $\pole_1,\ldots,\pole_\nred$ 
to be retained in the modal truncation,
and let $\cS$ denote the set $\Omega \cup \{\pole_{\nred+1},\ldots,\pole_\ns \}$.
Akin to \Cref{sec:multipoint}, pick interpolation points 
$\lp_1, \ldots, \lp_\nred \in {\C\setminus\cS }$ and 
\ $\rp_{1}, \ldots, \rp_{\nred} \in {\C\setminus \cS}$, 
and probing (direction) vectors $\Bell_1, \ldots, \Bell_\nred \in \C^{\nout}$ 
and $ \Br_{1}, \ldots, \Br_{\nred} \in \C^{\nin}$, 
then construct the interpolation data \cref{eqn:intdata} using the contour integrals
\begin{align} 
\label{eqn:modtrunleftdata}  
		\Bell_i^\ast  \TFr(\lp_i) &= \frac{1}{2\pi \iop} \int_{\partial \Omega} \frac{1}{\lp_i - z}\, \Bell_i^\ast  \TF(z)^{-1}\, \dop z, 
\rlap{\qquad$i=1,\ldots,\nred$,}\\[3pt]
		 \TFr(\rp_j) \Br_j&= \frac{1}{2\pi \iop} \int_{\partial \Omega} \frac{1}{\rp_j - z}\, \TF(z)^{-1} \Br_j\, \dop z,
\rlap{\qquad$j=1,\ldots,\nred$.}
\label{eqn:modtrunrightdata}  
		\end{align}
Given enough interpolation data,  one can use the tangential samples  {\cref{eqn:modtrunleftdata} and  \cref{eqn:modtrunrightdata}} 
in the Loewner framework to recover the desired modal truncation approximant $\TFr(z)$
\emph{without access to $\Amod$ and its spectral decomposition}. 
This approach requires that the interpolation point cannot 
coincide with the other poles (eigenvalues) outside $\Omega$.
(One can revert to the Hankel-based approach to guarantee this, 
thus sampling the Markov parameters of  $\TFr(z)$.  
Such an approach amounts to applying Sakurai--Sugiura~\cite{SS03}
or FEAST type methods~\cite{Pol09} for the standard eigenvalue problem.)

%% file: conclude.tex
\section{Conclusions}

Contour integral methods provide an effective tool for computing
eigenvalues in a bounded region of the complex plane.
By casting these algorithms in the framework of systems theory, 
we have proposed several new Loewner matrix methods inspired by
rational interpolation for system realization and model reduction.
Since these new techniques use the same quadrature data as existing
Hankel matrix methods, they can be implemented at little additional 
cost, yet can yield eigenvalue estimates with considerably improved 
accuracy.
Our primary intention has been to open up a broad family of 
algorithms for exploration and refinement.

The rational interpolation perspective suggests another, lower-cost
approach to the nonlinear eigenvalue problem $\BT(z)\Bv=\Bzero$:
use the Loewner framework to construct a linear matrix pencil that 
tangentially interpolates the nonlinear problem, in a sense detailed
in \Cref{sec:ratTinv}.
This method can be adapted into an
iterative method to refine the set of interpolation points.

Just as systems theory can inspire new eigenvalue algorithms,
contour integral methods provide a new approach to model reduction:
a way to perform modal truncation using only samples of the
transfer function without requiring access to a state-space representation, 
as described in \Cref{sec:modal}.

\begin{figure}[t!]
	\begin{center}
		$\overset{K = 1, \ r = 11}{\includegraphics[width=.32\textwidth]{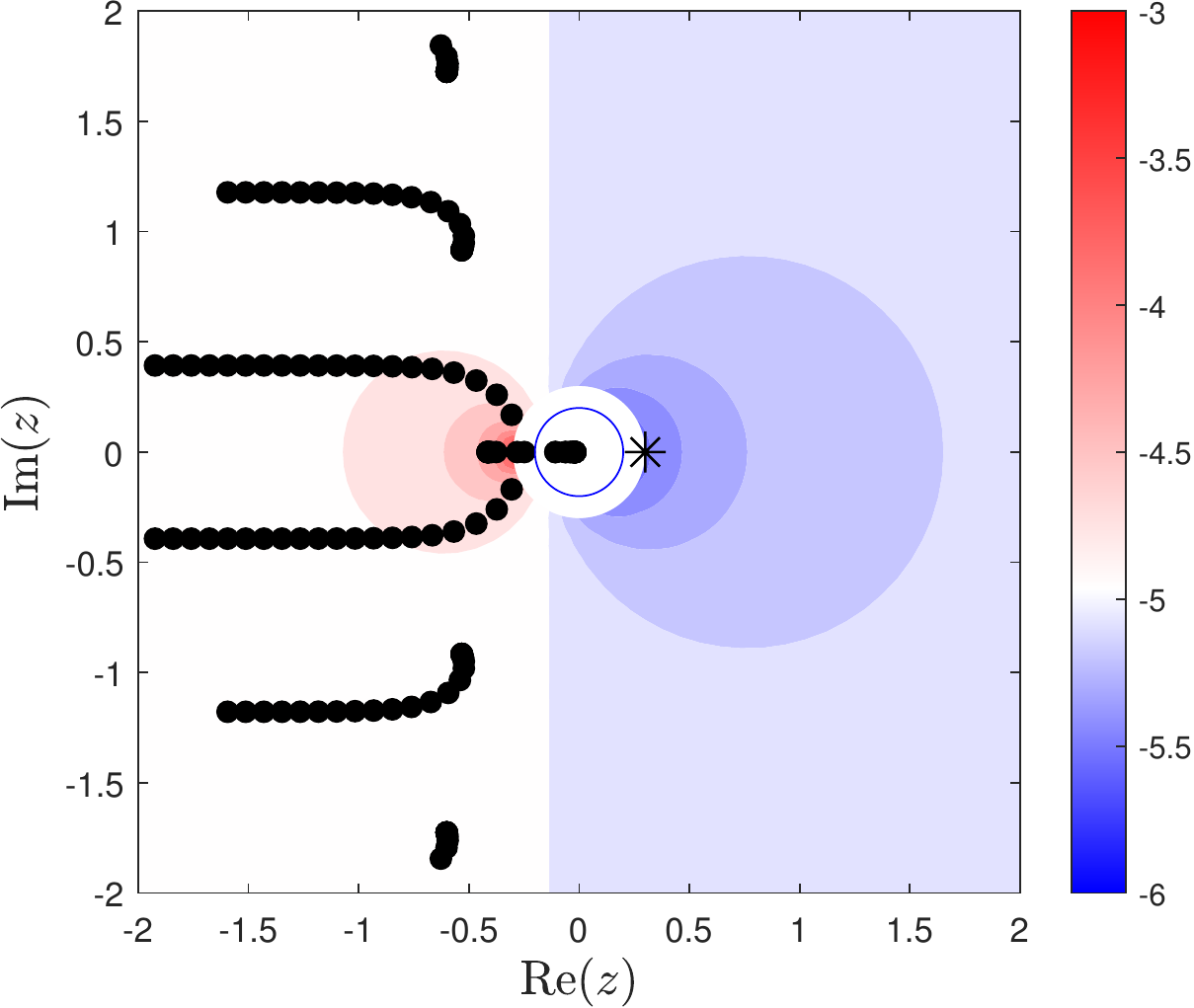}}$
		$\overset{K = 3, \ r = 11}{\includegraphics[width=.32\textwidth]{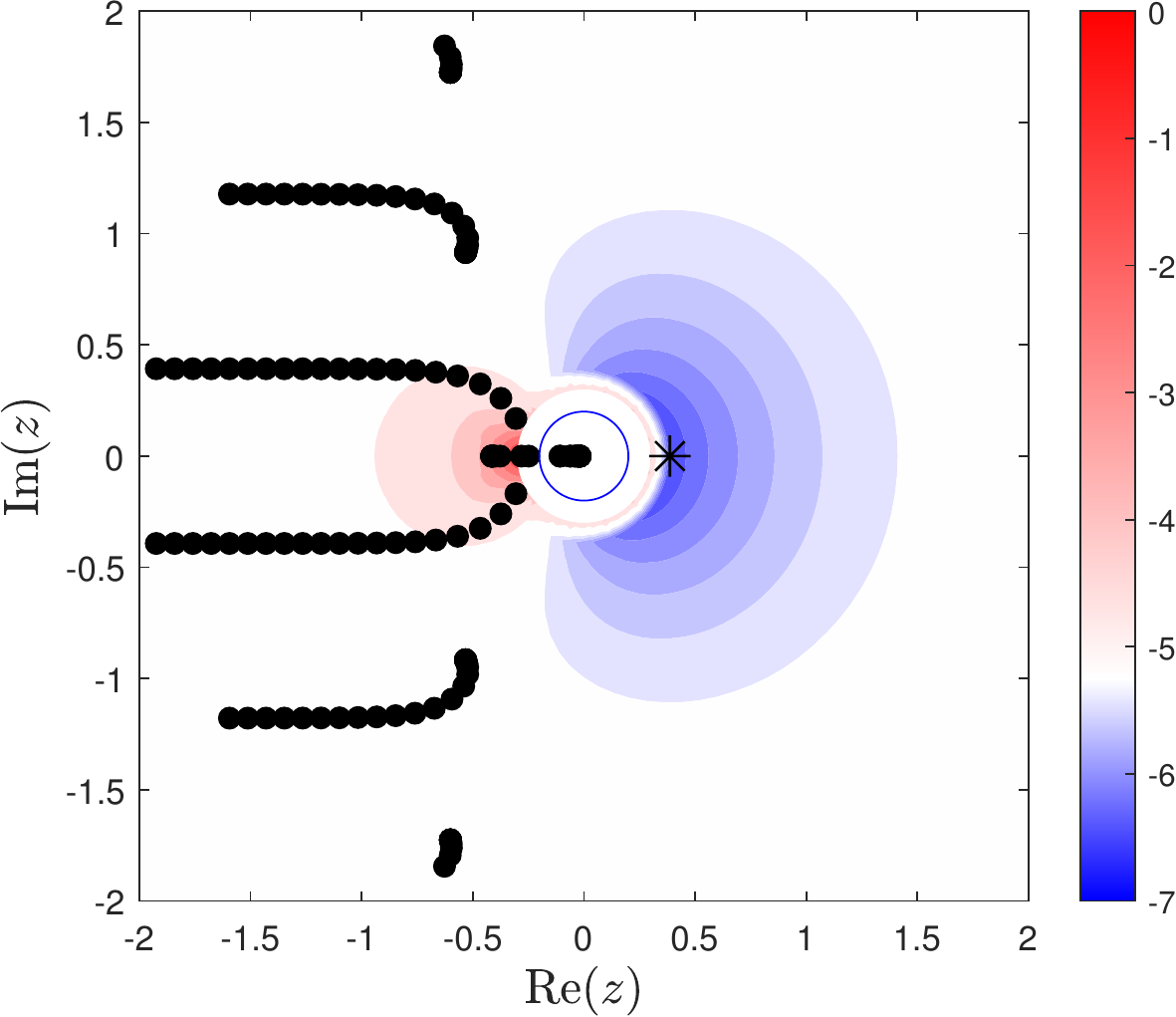}}$
		$\overset{K = 5, \ r = 11}{\includegraphics[width=.32\textwidth]{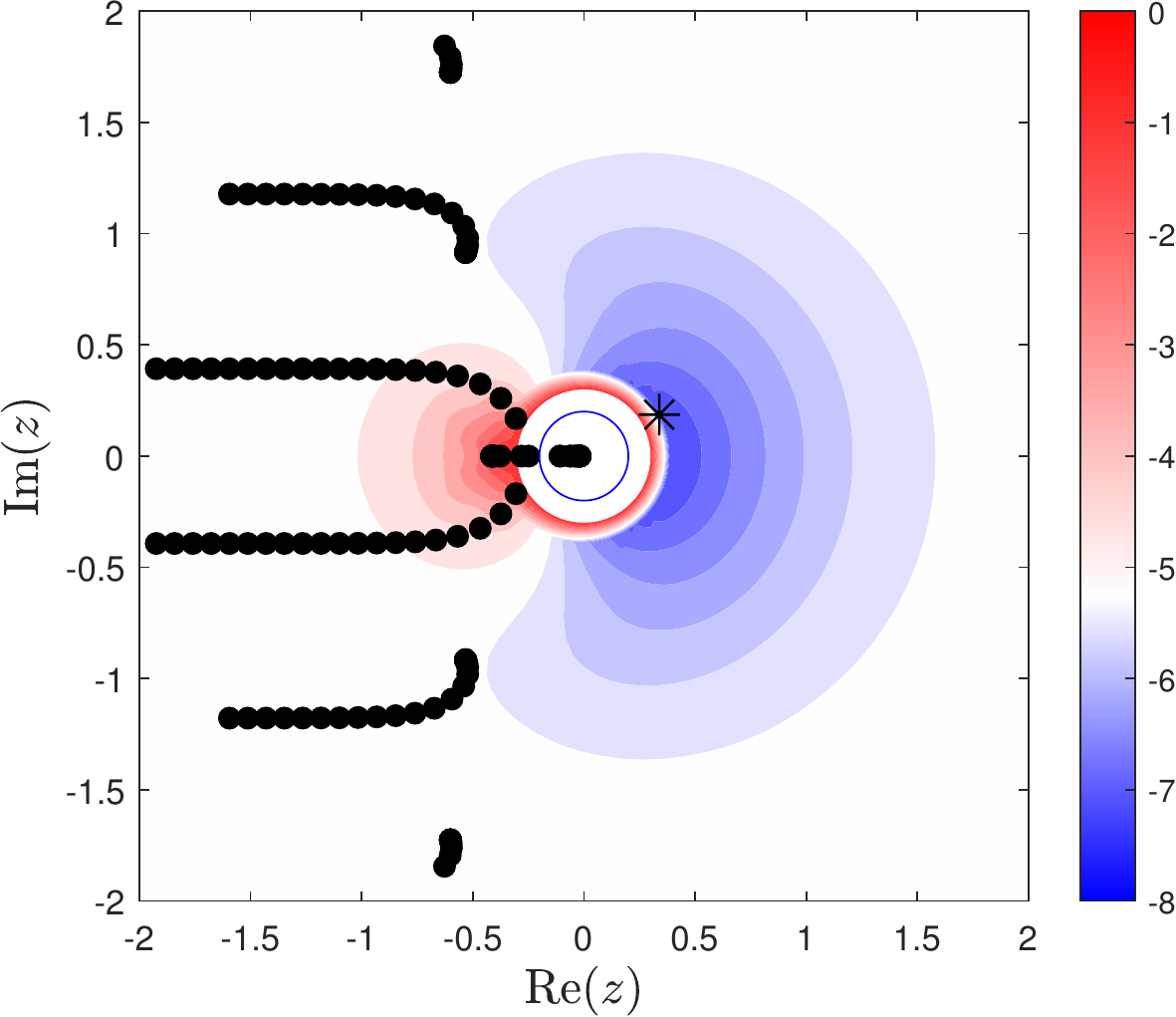}}$

\smallskip		
		$\overset{N = 32}{\includegraphics[width=.32\textwidth]{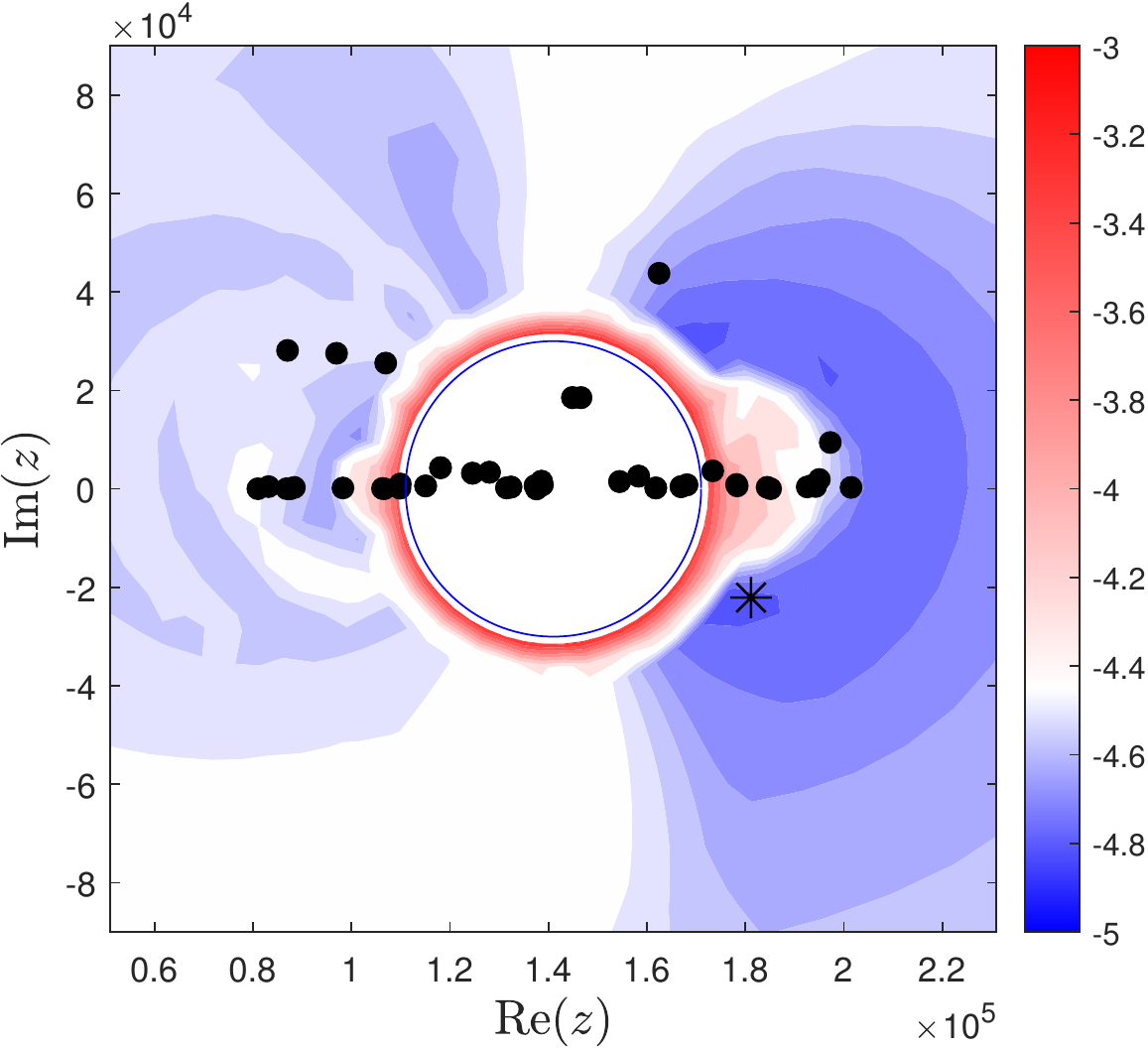}}$
		$\overset{N = 64}{\includegraphics[width=.32\textwidth]{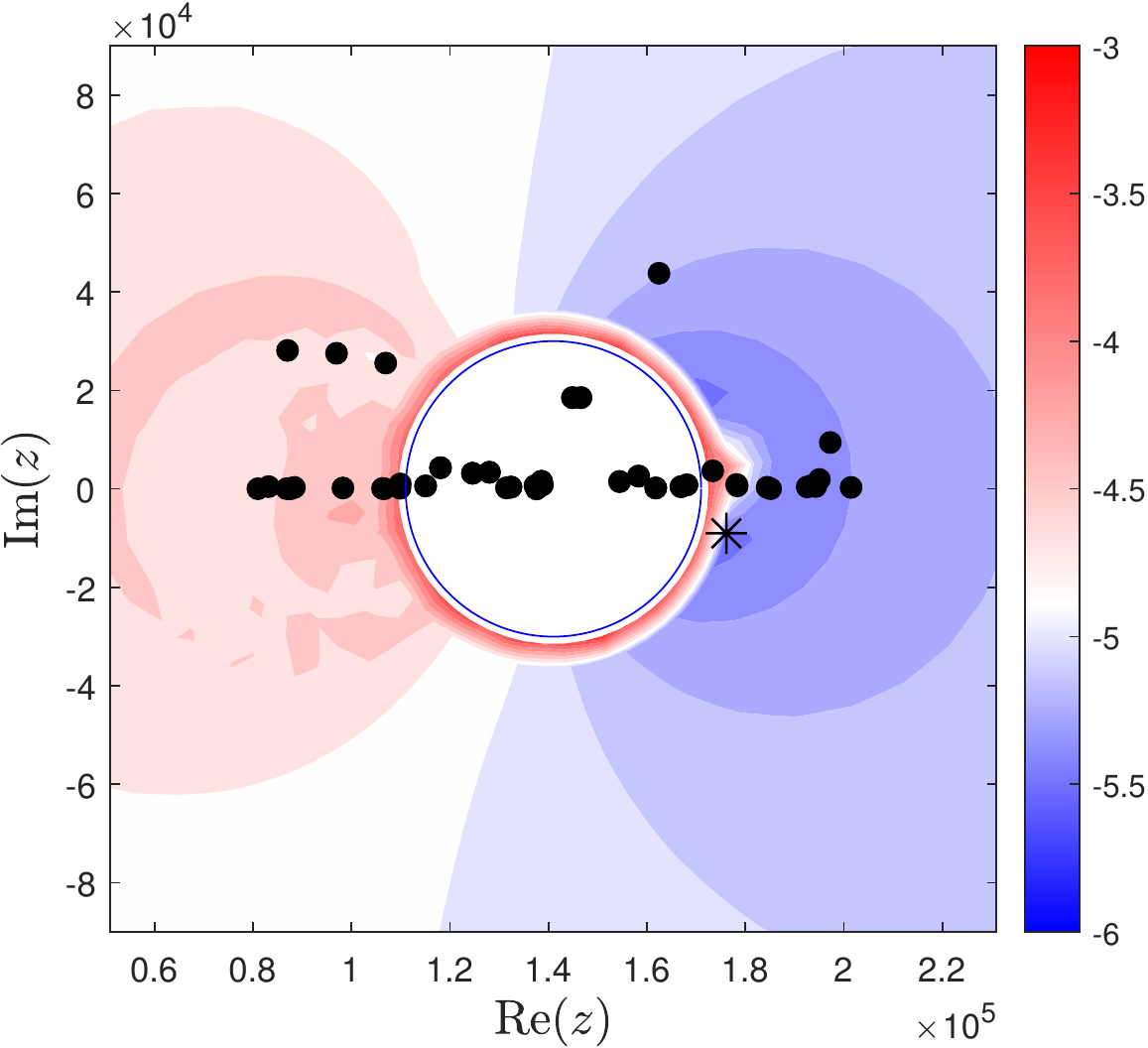}}$
		$\overset{N = 128}{\includegraphics[width=.32\textwidth]{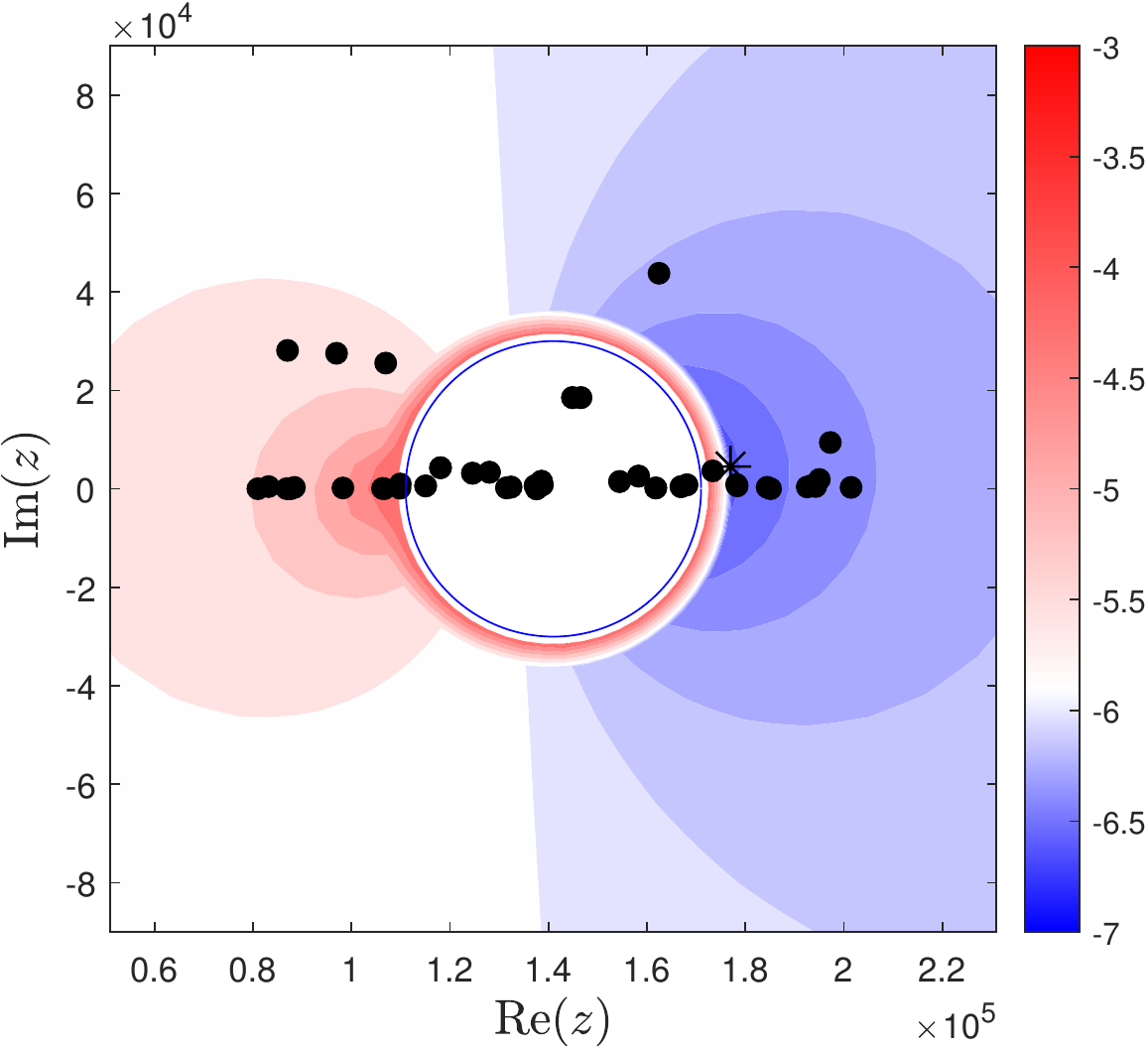}}$
	\end{center}
	
	\vspace*{-1em}
	\caption{\label{fig:delay_concl_contours}
		Filled contour plots showing the maximum eigenvalue residual error for the single-point Loewner method
applied to the delay (top) and {\tt gun} (bottom) problems, as a function of the interpolation point $\sigma\not\in\overline{\Omega}$;
$*$ denotes the best $\sigma$ in each plot. For the delay problem, we fix $N = 64$ quadrature points: as $K$ increases, 
the contours subtly change, but it is consistently better to place the interpolation point opposite the closest exterior eigenvalue. 
For the {\tt gun} problem, we fix $K = 1$ and vary $N$. For both problems, the color maps are defined so that white denotes 
the maximum residual error for the Hankel method, and the color depicts the $\log_{10}$ of the error for the single-point Loewner method:
the Loewner method outperforms the Hankel method for values of $\sigma$ in blue regions.
	}
\end{figure}

Many avenues for additional research remain open.  
The selection of optimal interpolation points and directions for Loewner methods, 
and the influence of the conditioning of the eigenvalues of the
associated Loewner pencil, remain important areas for investigation;
see~\cite{EI} for preliminary results in the context of system
realization. \Cref{fig:delay_concl_contours} shows contour plots 
of the maximum eigenvalue residual error for the single-point Loewner method 
applied to the delay and {\tt gun} problems, 
as a function of the interpolation point $\sigma\not\in\overline{\Omega}$. 
We see how the accuracy of the method indeed depends on the location of the interpolation point. 
(Our earlier experiments used \emph{good, but not optimal} values of $\sigma$.) 
These contour plots hint at the potential benefit that can come from theoretical 
insight into optimal interpolation point selection, as well as algorithms that can identify good candidates for $\sigma$.

%%%%%%%%%%%%%%%%%%%%%%%%%%%%%%%%%%%%%%%%%%%%%%%%%%%%%%%%%%%%%%%%%%%%%%%%%%%%%%%%
%%%%%%%%%%%%%%%%%%%%%%%%%%%%%%%%%%%%%%%%%%%%%%%%%%%%%%%%%%%%%%%%%%%%%%%%%%%%%%%%

One might also study the design of optimal filter functions to improve 
the convergence of the quadrature rules, the affect of quadrature errors 
on the computed eigenvalues,
how eigenvalues outside $\Omega$ affect convergence,
and how interpolation points should be placed relative 
to those exterior eigenvalues to give optimal accuracy.
Options abound for algorithm designers to explore and exploit.

%% file: paper.bbl
\begin{thebibliography}{10}

\bibitem{Ant05b}
{\sc A.~C. Antoulas}, {\em Approximation of Large-Scale Dynamical Systems},
  SIAM, Philadelphia, 2005.

\bibitem{ABG10}
{\sc A.~C. Antoulas, C.~A. Beattie, and S.~Gugercin}, {\em Interpolatory model
  reduction of large-scale dynamical systems}, in Efficient Modeling and
  Control of Large-Scale Systems, Springer, 2010, pp.~3--58.

\bibitem{AntBG20}
{\sc A.~C. Antoulas, C.~A. Beattie, and S.~G\"u\u{g}erc\.{\sc i}n}, {\em
  Interpolatory Methods for Model Reduction}, SIAM, Philadelphia, 2020.

\bibitem{ALI18}
{\sc A.~C. Antoulas, S.~Lefteriu, and A.~C. Ionita}, {\em A tutorial
  introduction to the {Loewner} framework for model reduction}, in Model
  Reduction and Approximation, SIAM, Philadelphia, 2017, pp.~335--376.

\bibitem{ASTIK09}
{\sc J.~Asakura, T.~Sakurai, H.~Tadano, T.~Ikegami, and K.~Kimura}, {\em A
  numerical method for nonlinear eigenvalue problems using contour integrals},
  JSIAM Letters, 1 (2009), pp.~52--55.

\bibitem{BG12}
{\sc C.~Beattie and S.~Gugercin}, {\em Realization-independent
  $\mathcal{H}_2$-approximation}, in 51st IEEE Conference on Decision and
  Control (CDC), IEEE, 2012, pp.~4953--4958.

\bibitem{BGW12}
{\sc C.~Beattie, S.~Gugercin, and S.~Wyatt}, {\em Inexact solves in
  interpolatory model reduction}, Linear Algebra Appl., 436 (2012),
  pp.~2916--2943.

\bibitem{BHMST13}
{\sc T.~Betcke, N.~J. Higham, V.~Mehrmann, C.~Schr{\"o}der, and F.~Tisseur},
  {\em {NLEVP}: A collection of nonlinear eigenvalue problems}, ACM
  Transactions on Mathematical Software (TOMS), 39 (2013), pp.~1--28.

\bibitem{Beyn12}
{\sc W.-J. Beyn}, {\em An integral method for solving nonlinear eigenvalue
  problems}, Linear Algebra Appl., 436 (2012), pp.~3839--3863.

\bibitem{BEGM05}
{\sc G.~Boutry, M.~Elad, G.~H. Golub, and P.~Milanfar}, {\em The generalized
  eigenvalue problem for non-square pencils using a minimal perturbation
  approach}, SIAM J.~Matrix Anal.\ Appl., 27 (2005), pp.~582--601.

\bibitem{Brennan18}
{\sc M.~C. Brennan}, {\em Rational interpolation methods for nonlinear
  eigenvalue problems}, master's thesis, Virginia Tech, 2018.

\bibitem{BP20}
{\sc J.~Brenneck and E.~Polizzi}, {\em An iterative method for contour-based
  nonlinear egensolvers}, arXiv preprint arXiv:2007.03000,  (2020).

\bibitem{CD05}
{\sc Z.~Chen and J.~J. Dongarra}, {\em Condition numbers of {Gaussian} random
  matrices}, SIAM J.~Matrix Anal.\ Appl., 27 (2005), pp.~603--620.

\bibitem{DeS00}
{\sc B.~De~Schutter}, {\em Minimal state-space realization in linear system
  theory: an overview}, J.~Comput.\ Appl.\ Math., 121 (2000), pp.~331--354.

\bibitem{DP}
{\sc Z.~Drma\v{c} and B.~Peherstorfer}, {\em Learning low-dimensional
  dynamical-system models from noisy frequency-response data with {Loewner}
  rational interpolation}.
\newblock arXiv:1910.00110; To appear in \emph{Realization and Model Reduction
  of Dynamical Systems: A Festschrift in Honor of the 70th Birthday of Thanos
  Antoulas}, C. A. Beattie, P. Benner, M. Embree, S. Gugergin, S. Lefteriu,
  eds.

\bibitem{EI}
{\sc M.~Embree and A.~C. Ionita}, {\em Pseudospectra of {Loewner} pencils}.
\newblock arXiv:1910.12153; To appear in \emph{Realization and Model Reduction
  of Dynamical Systems: A Festschrift in Honor of the 70th Birthday of Thanos
  Antoulas}, C. A. Beattie, P. Benner, M. Embree, S. Gugergin, S. Lefteriu,
  eds.

\bibitem{GAB08}
{\sc S.~Gugercin, A.~C. Antoulas, and C.~Beattie}, {\em {${\cal H}_2$} model
  reduction for large-scale linear dynamical systems}, SIAM J.~Matrix Anal.\
  Appl., 30 (2008), pp.~609--638.

\bibitem{GT17}
{\sc S.~G\"uttel and F.~Tisseur}, {\em The nonlinear eigenvalue problem}, Acta
  Numerica,  (2017), pp.~1--94.

\bibitem{Hig08}
{\sc N.~J. Higham}, {\em Functions of Matrices: Theory and Computation}, SIAM,
  Philadelphia, 2008.

\bibitem{HK66}
{\sc B.~L. Ho and R.~E. Kalman}, {\em Effective construction of linear
  state-variable models from input/output functions}, Regelungstechnik, 12
  (1966), pp.~545--548.

\bibitem{HMP}
{\sc M.~E. Hochstenbach, C.~Mehl, and B.~Plestenjak}, {\em Solving singular
  generalized eigenvalue problems by a rank-completing perturbation}, 2018,
  \href{http://arxiv.org/abs/1805.07657}{arXiv:1805.07657}.

\bibitem{HJ85}
{\sc R.~A. Horn and C.~R. Johnson}, {\em Matrix Analysis}, Cambridge University
  Press, Cambridge, 1985.

\bibitem{Kai80}
{\sc T.~Kailath}, {\em Linear Systems}, Prentice-Hall, Englewood Cliffs, NJ,
  1980.

\bibitem{Kel51}
{\sc M.~V. Keldysh}, {\em On the characteristic values and characteristic
  functions of certain classes of non-self-adjoint equations}, Doklady Akad.\
  Nauk SSSR (NS), 77 (1951), pp.~11--14.

\bibitem{Keld71}
{\sc M.~V. Keldysh}, {\em On the completeness of the eigenfunctions of some
  classes of non-selfadjoint linear operators}, Russian Math.\ Surveys, 26
  (1971), pp.~15--44.

\bibitem{MA07}
{\sc A.~J. Mayo and A.~C. Antoulas}, {\em A framework for the solution of the
  generalized realization problem}, Linear Algebra Appl., 425 (2007),
  pp.~634--662.

\bibitem{MV04}
{\sc V.~Mehrmann and H.~Voss}, {\em Nonlinear eigenvalue problems: a challenge
  for modern eigenvalue methods}, GAMM-Mitt., 27 (2004), pp.~121--152.

\bibitem{MN14}
{\sc W.~Michiels and S.-I. Niculescu}, {\em Stability and Stabilization of
  Time-Delay Systems: An Eigenvalue-Based Approach}, SIAM, Philadelphia,
  second~ed., 2014.

\bibitem{Pol09}
{\sc E.~Polizzi}, {\em Density-matrix-based algorithm for solving eigenvalue
  problems}, Phys.\ Rev.~B, 79 (2009), p.~115112.

\bibitem{RN06}
{\sc J.~Rommes and N.~Martins}, {\em Efficient computation of transfer function
  dominant poles using subspace acceleration}, IEEE Trans. Power Systems, 21
  (2006), pp.~1218--1226.

\bibitem{RS08}
{\sc J.~Rommes and G.~L.~G. Sleijpen}, {\em Convergence of the dominant pole
  algorithm and {Rayleigh} quotient iteration}, SIAM J.~Matrix Anal.\ Appl., 30
  (2008), pp.~346--363.

\bibitem{SS03}
{\sc T.~Sakurai and H.~Sugiura}, {\em A projection method for generalized
  eigenvalue problems using numerical integration}, J.~Comput.\ Appl.\ Math.,
  159 (2003), pp.~119--128.

\bibitem{Silv71}
{\sc L.~Silverman}, {\em Realization of linear dynamical systems}, IEEE Trans.\
  Auto.\ Control, 16 (1971), pp.~554--567.

\bibitem{Ste94}
{\sc G.~W. Stewart}, {\em Perturbation theory for rectangular matrix pencils},
  Linear Algebra Appl., 208/209 (1994), pp.~297--301.

\bibitem{TP19}
{\sc F.~Tisseur and G.~Porzio}, {\em An algorithm for dense nonlinear
  eigenvalue problems}, July 2019.
\newblock Presentation, 9th International Congress on Industrial and Applied
  Mathematics, Valencia, Spain.

\bibitem{TW14}
{\sc L.~N. Trefethen and J.~A.~C. Weideman}, {\em The exponentially convergent
  trapezoidal rule}, SIAM Review, 56 (2014), pp.~385--458.

\bibitem{TYUC17}
{\sc J.~A. Tropp, A.~Yurtsever, M.~Udell, and V.~Cevher}, {\em Practical
  sketching algorithms for low-rank matrix approximation}, SIAM J.~Matrix
  Anal.\ Appl., 38 (2017), pp.~1454--1485.

\bibitem{Bare16}
{\sc M.~Van~Barel}, {\em Designing rational filter functions for solving
  eigenvalue problems by contour integration}, Linear Algebra Appl., 502
  (2016), pp.~346--365.

\bibitem{BK16}
{\sc M.~Van~Barel and P.~Kravanja}, {\em Nonlinear eigenvalue problems and
  contour integrals}, J.~Comput.\ Appl.\ Math., 292 (2016), pp.~526--540.

\bibitem{Vos14}
{\sc H.~Voss}, {\em Nonlinear eigenvalue problems}, in Handbook of Linear
  Algebra, L.~Hogben, ed., CRC/Taylor \& Francis, Boca Raton, FL, second~ed.,
  2014, ch.~60.

\bibitem{Woo14}
{\sc D.~P. Woodruff}, {\em Sketching as a tool for numerical linear algebra},
  Found.\ Trends Theoret.\ Comput.\ Sci., 10 (2014), pp.~1--157.

\bibitem{WT02}
{\sc T.~G. Wright and L.~N. Trefethen}, {\em Pseudospectra of rectangular
  matrices}, IMA J.~Numer.\ Anal., 22 (2002), pp.~501--519.

\end{thebibliography}
